\documentclass[bj,authoryear, doublespacing, noshowframe]{imsart}

\usepackage{epsfig}
\usepackage{amssymb,amsmath, amsthm,url,graphicx,color}
\usepackage{multirow}
\usepackage{booktabs}
\usepackage{comment}

\usepackage{array}
\usepackage{placeins}

\usepackage{xcolor}

\startlocaldefs
\newtheorem{theorem}{Theorem}[section]

\newtheorem{proposition}{Proposition}[section]

\newtheorem*{assumption*}{Assumptions}

\theoremstyle{definition}

\numberwithin{equation}{section}
\numberwithin{theorem}{section}
\numberwithin{lemma}{section}
\numberwithin{remark}{section}
\numberwithin{example}{section}
\numberwithin{table}{section}
\numberwithin{figure}{section}
\numberwithin{definition}{section}

\newcommand{\mbf}{\mathbf}
\newcommand{\mbb}{\mathbb}
\newcommand{\mcal}{\mathcal}
\newcommand{\bs}{\boldsymbol}

\usepackage{float}
\floatstyle{ruled}
\newfloat{algo}{thp}{lop}
\floatname{algo}{Algorithm}

\allowdisplaybreaks

\DeclareMathOperator*{\argmin}{argmin}

\endlocaldefs

\begin{document}

\begin{frontmatter}
\title{Spectral estimation for high-dimensional linear processes}
\runtitle{Spectral estimation for high-dimensional linear processes}

\begin{aug}
\author[A]{\inits{J.}\fnms{Jamshid}~\snm{Namdari}\ead[label=e1]{jamshid.namdari@emory.edu}}
\author[B]{\inits{A.}\fnms{Alexander}~\snm{Aue}\ead[label=e2]{aaue@ucdavis.edu}}
\author[B,C]{\inits{D.}\fnms{Debashis}~\snm{Paul}\ead[label=e3]{debpaul@ucdavis.edu}}
\address[A]{Department of Biostatistics \& Bioinformatics, Emory University \printead[presep={,\ }]{e1}}

\address[B]{Department of Statistics, University of California, Davis \printead[presep={,\ }]{e2}}

\address[C]{Applied Statistics Unit, Indian Statistical Institute, Kolkata\printead[presep={,\ }]{e3}}
\end{aug}

\begin{abstract}
We develop a random-matrix-based framework for estimating the latent joint spectral distribution of the coefficient matrices and innovation covariance of a class of high-dimensional linear processes from the eigenvalues of
weighted sample periodograms. The model assumes that these matrices share a common orthonormal eigenbasis, while the dimension $p$ and the sample size $n$ satisfy $p/n \to c \in (0,\infty)$. We derive the almost-sure limiting spectral distribution of frequency-weighted integrals of the sample periodogram and characterize its Stieltjes transform by a system of nonlinear integral equations involving the unknown joint spectrum. For a discrete spectrum supported on a prescribed grid, we estimate the mixture weights by matching empirical and limiting Stieltjes transforms with respect to the general $L_\kappa$ discrepancy, $\kappa\geq 1$, and establish consistency for $\kappa=2$. We further propose an estimator of the spectral density matrix and a bootstrap procedure to select among the candidate ARMA models. Simulations demonstrate accurate recovery of the underlying spectral distributions across several autoregressive and moving-average models, and an application to S\&P~500 returns reveals temporal dependence beyond the dominant market factor.
\end{abstract}


\end{frontmatter}

\section{Introduction}
\label{sec:int}

High-dimensional time series arise in many modern applications, including finance, climatology, genomics, and neuroscience, where the dimension of the observations is comparable to or larger than the sample size. In this regime, classical multivariate estimators are no longer consistent, and random matrix theory has become an important tool for understanding the behavior of covariance-type matrices and developing statistically valid inference procedures; see \citet{pa13}. Although these developments are now well established for independent observations, analogous estimation problems for dependent high-dimensional data remain considerably less understood.

For stationary multivariate linear processes, the spectral density matrix and the autocovariance sequence play a central role in estimation, prediction, and model selection. In high dimensions, however, sample covariance matrices, sample autocovariances, and sample periodograms exhibit non-negligible random matrix effects, implying that their empirical spectral distributions do not converge to their population counterparts. Consequently, classical procedures for estimating spectral density matrices or identifying autoregressive and moving-average models become unreliable when the dimension is of the same order as the sample size.

Recent work has established limiting spectral distributions for various high-dimensional autocovariance matrices. \cite{jwbnh14} derived the limiting empirical spectral distribution of symmetrized sample autocovariances under independent noise, while \citet{lap15} extended these results to a substantially richer class of linear processes whose coefficient matrices are simultaneously diagonalizable. \citet{bb16} further relaxed the simultaneous diagonalization assumption using tools from free probability. These papers provide a detailed understanding of the forward asymptotic behavior of sample autocovariance matrices and related quantities. However, unlike the corresponding literature for covariance matrices of independent observations, they do not address the inverse statistical problem of estimating the underlying spectral distribution of the model parameters.
The present paper addresses this problem. We study high-dimensional linear processes $(X_t\colon t\in\mathbb{Z})$ given by
\begin{equation}
\label{eq:linproc}
X_t = \sum_{\ell=0}^\infty \mathbf{A}_\ell Z_{t-\ell},\qquad t\in\mathbb{Z},
\end{equation}
under the assumption that the coefficient matrices and the innovation covariance matrix are simultaneously diagonalizable. This assumption allows the multivariate process to be represented, after a common orthogonal transformation, as a collection of uncorrelated univariate linear processes, making it possible to characterize the limiting spectral behavior of weighted sample periodograms through the joint spectral distribution of the coefficient matrices and the innovation covariance. 

The starting point is a new characterization of the limiting spectral distribution of weighted integrals of the sample periodogram under the high-dimensional asymptotic regime 
\begin{equation}\label{eq:pn}
p,n\to\infty\quad\mbox{such that}\quad\frac pn\to c\in(0,\infty).
\end{equation} 
The resulting system of nonlinear equations links the empirical spectral distributions of weighted sample periodograms to the unknown joint spectral distribution of the coefficient matrices and the innovation covariance. This relationship forms the basis of a new estimation procedure that solves an inverse random matrix problem for high-dimensional time series.

While spectrum estimation is well developed for high-dimensional independent data \citep{e08,lw15}, analogous methods for high-dimensional time series remain limited. We address this problem for the class of linear processes studied by \cite{lap15}, in which the coefficient matrices are simultaneously diagonalizable. Although restrictive, this assumption allows for nontrivial cross-sectional dependence and arises naturally in several multivariate settings, including blind source separation \citep{minott16,mnot14}, common principal component analysis \citep{flury84}, and, more generally, joint diagonalization methods in multivariate statistics \citep{nr22}. It has also proved useful in the random matrix analysis of high-dimensional linear processes and sample autocovariance matrices \citep{lap15,wap17}. Within this framework, we develop an inverse random matrix approach for estimating the underlying joint spectral distribution from weighted sample periodograms. The joint spectral distribution provides a basis-invariant summary of the heterogeneous dynamic behavior across coordinates. It describes the limiting proportions of directions exhibiting different autoregressive, moving-average, and innovation-variance characteristics, without requiring recovery of the corresponding individual eigenvectors.

The main contributions of the paper are fourfold.
\begin{itemize}
\item First, we establish the limiting spectral distribution of weighted sample periodograms and characterize its Stieltjes transform through a system of nonlinear integral equations.
\item Second, assuming that the joint spectral distribution is supported on a finite grid, we estimate the corresponding mixture weights by matching empirical and theoretical Stieltjes transforms and prove the consistency of the estimator under a $L_2$-discrepancy criterion.
\item Third, we develop an estimator of the spectral density matrix together with a bootstrap-based procedure for selecting among competing high-dimensional ARMA models.
\item
Finally, simulation studies demonstrate that the proposed procedure accurately recovers the underlying spectral distributions across a range of autoregressive and moving-average models, while an application to S\&P 500 returns reveals temporal dependence beyond the dominant market factor.
\end{itemize}

The remainder of the paper is organized as follows. Section \ref{sec:Main_Results} derives the limiting spectral distribution of weighted sample periodograms. Sections 3 and 4 introduce the estimation procedures for the joint spectral distribution and the spectral density matrix, respectively. Section 5 establishes consistency, Section 6 develops model selection, Sections 7 and 8 present simulation studies and a data analysis, and Section 9 concludes.

\section{Spectral distributions of weighted integrals of the sample periodogram} \label{sec:Main_Results}

Let the data matrix $\mathbf{X_n}=[X_1,\dots,X_n]$ be obtained from the $p$-dimensional linear process $(X_t\colon t\in\mathbb{Z})$ in \eqref{eq:linproc}, with
iid innovations $(Z_t\colon t\in\mathbb{Z)}$ having zero mean and covariance matrix $\Sigma$. Throughout the paper we consider the high-dimensional asymptotic regime \eqref{eq:pn}. Our objective is to estimate the joint spectral distribution of the coefficient matrices $(\mathbf{A}_\ell\colon\ell\geq 0)$ together with the innovation covariance matrix $\Sigma$. Rather than working directly with sample autocovariance matrices, we consider weighted integrals of the multivariate sample periodogram. Their empirical spectral distributions admit deterministic limits that depend only on the unknown joint spectrum and therefore provide the basis for statistical estimation developed in Section 3.

To obtain tractable limiting equations, we impose the following structural assumptions. Let $\mathbb{N}$ and $\mathbb{N}_0$ denote the set of positive and nonnegative  integers, respectively.

\begin{itemize}
\item[\textbf{A.1}]
The process $(Z_t\colon t\in\mathbb{Z})$ is represented as $Z_t=\Sigma^{1/2}\tilde{Z}_t$ where $\Sigma^{1/2}$ is a square-root of $\Sigma$, and the $p$-dimensional vectors $\tilde{Z}_t$ have iid\ entries $\tilde{Z}_{jt}$ satisfying one 
of the following conditions:
\begin{itemize}
\item
(for complex-valued processes): $\tilde{Z}_{jt}$ is complex valued,
with zero mean, and
$\mathbb{E}[\Re(\tilde{Z}_{jt})^2] =\mathbb{E}[\Im(\tilde{Z}_{jt})^2]=1/2$ where $\Re(\tilde{Z}_{jt})$ and $\Im(\tilde{Z}_{jt})$ are the real and imaginary parts of $\tilde{Z}_{jt}$, which are independent, and $\mathbb{E}[|\tilde{Z}_{jt}|^4] < \infty$; 
\item
(for real-valued processes):
$\tilde{Z}_{jt}$ is real-valued, with zero mean, $\mathbb{E}[|\tilde{Z}_{jt}|^2]=1$, and $\mathbb{E}[|\tilde{Z}_{jt}|^4] < \infty$.
\end{itemize}

\item[\textbf{A.2}]
	The matrices $(\mathbf{A}_\ell\colon \ell \in \mathbb{N}_0)$ are simultaneously diagonalizable, random Hermitian matrices, independent of $(Z_t\colon t \in \mathbb{Z})$, satisfying $\|\mathbf{A}_\ell\| \leq \bar\lambda_\ell$ for all $\ell \in \mathbb{N}_0$ and large $p$ with
	\begin{align*}
		\sum_{\ell=0}^{\infty}\bar\lambda_\ell < \infty \qquad \mbox{and} \qquad \sum_{\ell=0}^{\infty}\ell\bar\lambda_\ell  < \infty. 
	\end{align*} 
\item[\textbf{A.3}] 
There is a compact subset $\mathbb{K}$ of $\mathbb{R}^m$ for some $m < \infty$, and there are differentiable functions $f_\ell\colon\mathbb{K}\subset\mathbb{R}^m\rightarrow\mathbb{R}$, $\ell \in \mathbb{N}_0$, such that $\sum_{\ell=0}^{\infty}\|f_\ell^\prime\|_\infty < \infty$,
	where $f_\ell^\prime$ denotes the derivative of $f_\ell$. For each $p$, there is a set of points $\lambda_{p,1}, \ldots, \lambda_{p,p} \in \mathbb{K}$, not necessarily distinct, and a unitary $p\times p$ matrix $\mathbf{U}$ such that 
	\begin{align*}
		\mathbf{U}^*\mathbf{A}_\ell \mathbf{U}=\mathrm{diag}(f_\ell(\lambda_{p,1}),\dots, f_\ell(\lambda_{p,p})), \qquad \ell \in \mathbb{N},
	\end{align*}
	and $f_0(\lambda)=1$. 
	In addition, we assume that $\mathbf{U}^*\Sigma \mathbf{U}=\mathrm{diag}\{\sigma_{p,1},\dots,\sigma_{p,p}\}$, where $\sigma_{p,1},\dots,\sigma_{p,p}\in\mathbb{K}_\sigma\subset\mathbb{R}^+$, $\mathbb{K}_\sigma$ compact.
	
\item[\textbf{A.4}]
	Let  $\boldsymbol{\lambda}_k=(\lambda_{p,k},\sigma_{p,k})$ for $k=1,\dots,p$, where $\lambda_{p,k}=(\lambda_{p,k}^{(1)},\ldots,\lambda_{p,k}^{(m)})\in\mathbb{K}$ are as above. Then 
	the empirical distribution of $\{\boldsymbol{\lambda}_1, \ldots , \boldsymbol{\lambda}_p\}$, denoted by $F_p^{A,\Sigma}$,  converges weakly to a probability distribution function $F^{A,\Sigma}$ on $\mathbb{R}^{m+1}$ as $p \to \infty$.
\end{itemize}

The simultaneous diagonalizability assumption implies that there exists a common orthogonal (or unitary) basis $\mathbf{U}$ in which the transformed process
$\mathbf{U}^*X_t$ has uncorrelated coordinates. Consequently, each coordinate is a univariate, stationary linear process whose coefficients are determined by the corresponding eigenvalues of the matrices $\mathbf{A}_\ell$ and $\Sigma$. The limiting distribution $F^{A,\Sigma}$ introduced in Assumption A.4 therefore plays the role of a population distribution describing the random coefficients of these coordinate processes.

Assumptions A.2 and A.3 guarantee short-range dependence and regularity of the coefficient functions, while Assumption A.4 ensures the existence of a limiting joint spectral distribution as $p\to\infty$.

For $\boldsymbol{\lambda}=(\lambda,\sigma)\in\mathbb{K}\times\mathbb{K}_\sigma$ and $\theta\in[0,2\pi]$ define
\begin{equation}\label{eq:def_psi_and_h}
	\psi(\boldsymbol{\lambda}, \theta)=\sum_{l=0}^{\infty}e^{il\theta}\sigma f_l(\bs\lambda) \qquad \mbox{and} \qquad h(\boldsymbol{\lambda}, \theta)=|\psi(\boldsymbol{\lambda}, \theta)|^2.
\end{equation}
The quantity $h(\bs\lambda_k,\theta)$ is proportional to the spectral density of the
$k$-th coordinate process. Throughout the paper it serves as the fundamental building block linking the unknown population spectrum with the empirical spectral distributions of the sample periodogram. 

For a bounded nonnegative weight function $g$, define
\begin{equation} \label{eq:Weighted_Sum_Sample_Periodogram}
\mathbf{S}_g^{(n)} 
= \frac{1}{n}\sum_{t=1}^{n} g(\theta_t)\tilde{X}_t\tilde{X}_t^* 
\end{equation}
where $\theta_t=2\pi t/n$ are the Fourier frequencies and {\color{black} $\tilde{X_j} = n^{-1/2}\sum_{t=1}^{n}X_je^{-it\theta_j}$ denotes the discrete Fourier transform of the observations.} 
The matrix $\mathbf{S}_g^{(n)}$ is a weighted integral of the sample periodogram. At the population level, the corresponding weighted spectral quantity is
\[
\frac{1}{2\pi}\int_0^{2\pi} g(\theta) f(\theta)\,d\theta,
\]
where
\[
f(\theta)
=
\sum_{h=-\infty}^{\infty}
\boldsymbol{\Gamma}(h)e^{-ih\theta},
\qquad 0\leq\theta\leq2\pi,
\]
and $\boldsymbol{\Gamma}(h)$ denotes the lag-$h$ autocovariance matrix. By choosing weight functions concentrated on different regions of the frequency domain, one can extract complementary information about the temporal dependence structure. This observation will play an important role for the estimation of the joint spectral distribution proposed in Section \ref{sub:SpecEst_Algm}.

Define $\mathbf{W}_g=\mathrm{diag}\{\sqrt{g(\theta_1)},\dots,\sqrt{g(\theta_n)}\}$ and  $\tilde{\bs X}_n = [\tilde{X}_1,\dots,\tilde{X}_{n}]$ and note that $\mathbf{S}_g^{(n)} = \frac{1}{n}\boldsymbol{\tilde{X}}_n\mathbf{W}_g^2\boldsymbol{\tilde{X}}_n^*$. To characterize the limiting empirical spectral distribution of $\mathbf{S}_g^{(n)}$ it is convenient to work with its \textit{dual} matrix
\begin{equation}\label{eq:S_n_tilde_def}
\tilde{\mathbf{S}}_g^{(n)}=\frac{1}{n}\mathbf{W}_g\boldsymbol{\tilde{X}}_n^*\boldsymbol{\tilde{X}}_n\mathbf{W}_g,
\end{equation}
which has the same nonzero eigenvalues as $\mathbf{S}_g^{(n)}$. Denote its empirical spectral distribution (ESD) by $F_g^{(n)}$.
We need the following additional assumption for the function $g$. 
\begin{itemize}
	\item[\textbf{A.5}]
	Let $g\colon [0,2\pi]\to\mathbb{R}^+$ be a bounded function such that $|g(\theta_2)-g(\theta_1)| \leq c_g |\theta_2 - \theta_1|$ for some constant $c_g >0$.
\end{itemize}

The following theorem establishes the almost sure limit of $F_g^{(n)}$. Recall that, for a probability distribution function $F$ on $\mathbb{R}$, its \textit{Stieltjes transform} is defined by
\[
\mathcal{S}_F\colon\mathbb{C}^+\to\mathbb{C}^+,\qquad
z\mapsto\mathcal{S}_F(z)=\int\frac{1}{\lambda-z}F(dx),
\]
where $\mathbb{C}^+=\{x+iy\colon x\in\mathbb{R},y>0\}$. The Stieltjes transform uniquely determines $F$ and will be used to characterize the limiting spectral distribution.

\begin{theorem} \label{theorem1}
Consider the linear process $(X_t\colon t \in \mathbb{Z})$ satisfying 
assumptions \textbf{A.1}--\textbf{A.5}, and suppose that $p/n \to c \in (0,\infty)$ as $p,n\to \infty$. Then,
$F_g^{(n)}$ converges almost surely to a nonrandom probability distribution $F_g$ with Stieltjes transform $\mathcal{S}_g(z)=\mathcal{S}_{F_g}(z)$ determined by the equations
\begin{equation} \label{eq:s}
\mathcal{S}_g(z)=\frac{1}{2\pi}\int_{0}^{2\pi}\frac{1}{cM_g(z,\theta)-z}d\theta,
\end{equation}
\begin{equation} \label{eq:K}
K_g(z,\boldsymbol{\lambda})=\frac{1}{2\pi}\int_{0}^{2\pi}\frac{g(\theta)h(\boldsymbol{\lambda},\theta)}{cM_g(z,\theta)-z}d\theta,
\end{equation}
\begin{equation} \label{eq:M}
M_g(z,\theta)=\int\frac{g(\theta)h(\boldsymbol{\lambda},\theta)}{1+K_g(\boldsymbol{\lambda},z)}dF^{A,\Sigma}(\boldsymbol{\lambda}),
\end{equation}
where $K_g\colon \mathbb{K}\times\mathbb{C}^+ \to \mathbb{C}^+$ is the unique solution to (\ref{eq:K}) subject to the restriction that $K_g(\boldsymbol{\lambda},z)$ is a Stieltjes transform  with total mass $m_{\boldsymbol{\lambda}}^{(g)}=\frac{1}{2\pi}\int_{0}^{2\pi}g(\theta)h(\boldsymbol{\lambda}, \theta)d\theta$.
\end{theorem}

Theorem~\ref{theorem1} provides a deterministic characterization of the limiting spectral distribution through a coupled system of nonlinear equations for its Stieltjes transform. The unknown joint spectral distribution $F^{A,\Sigma}$ enters only through equation \eqref{eq:M}, making these equations suitable for statistical inversion. This observation forms the basis of the estimation procedure proposed in the next section.

The proof follows the deterministic-equivalent approach developed in random matrix theory. A key step consists of deriving a deterministic approximation to the resolvent of the dual matrix, from which the limiting equations for the Stieltjes transform follow. Details are deferred to the Appendix.

\section{Estimation procedure} \label{sub:SpecEst_Algm}

We now use the limiting spectral characterization developed in Section 2 to estimate the joint spectral distribution of the coefficient matrices and the innovation covariance matrix. To distinguish the population distributions considered in the estimation procedure, throughout this section we denote the true spectral distribution by $F^0=F^{A^0,\Sigma^0}$. The corresponding quantities in the limiting equations \eqref{eq:s}--\eqref{eq:M} are denoted by $\mathcal{S}_g^0$, $K_g^0$ and $M_g^0$. The estimation strategy is based on the observation that Theorem~\ref{theorem1} determines the limiting spectral distribution of each weighted sample periodogram as a function of the unknown population distribution $F^0$. We consider a class of candidate distributions $F$ and select the one whose implied limiting spectral distributions most closely match the empirical spectral distributions obtained from the data. The details of the estimation procedure are detailed in the following.

\medskip



Motivated by spectrum estimation for high-dimensional covariance matrices, see for example \cite{e08}, we represent candidate joint spectral distributions $F=F^{A,\Sigma}$ by a mixture of point masses supported on a prespecified grid $\Lambda^0=\{\boldsymbol{\lambda}_1^0,\ldots,\boldsymbol{\lambda}_J^0\}
\subset \mathbb{K}\times \mathbb{K}_\sigma$. Any such distribution can be written as
\begin{equation}\label{eq:F}
F(\boldsymbol{\lambda})
=F^{A,\Sigma}(\boldsymbol{\lambda})
=\sum_{j=1}^J\omega_j
\mathbb{I}\{\boldsymbol{\lambda}_j^0
\leq\boldsymbol{\lambda}\}
\end{equation}
and is therefore characterized by the weights
$\boldsymbol{\omega}=(\omega_1,\ldots,\omega_J)$. To be admissible, these weights have to be chosen from the $(J-1)$-dimensional simplex $\Delta_J=\{\boldsymbol{\omega}\in\mathbb{R}^J\colon\omega_j\geq 0,\,\|\boldsymbol{\omega}\|_1=1\}$. Note that under this discretized model $F^0$ is represented by replacing $\boldsymbol{\omega}$ with the true mixture weights, say $\boldsymbol{\omega}^0$. Given $\Lambda^0$, estimation of $F^0$ reduces to estimation of the mixture weights $\boldsymbol{\omega}^0$. 

To see how Theorem~\ref{theorem1} leads to an estimator of these weights, we first specialize its limiting equations to the discrete model above. Under the discrete representation \eqref{eq:F}, the integral equations in \eqref{eq:s}--\eqref{eq:M} reduce to
\begin{align}
\label{s0}
{\color{black} 
\mathcal{S}_g(z|\bs\omega)}&=
\frac{1}{2\pi}\int_{0}^{2\pi}\frac{1}{cM_g(z,\theta|\bs\omega)-z}d\theta,\\
\label{K0}
K_g(z,\boldsymbol{\lambda}^0_j|\bs\omega)&=
\frac{1}{2\pi}\int_{0}^{2\pi}\frac{g(\theta)h(\boldsymbol{\lambda}_j^0,\theta)}{cM_g(z,\theta|\bs\omega)-z}d\theta, \qquad j=1,\ldots,J, \\
\label{M0}
M_g(z,\theta|\bs\omega)&=\sum_{j=1}^{J}\omega_j\;\frac{g(\theta)h(\boldsymbol{\lambda}_j^0,\theta)}{1+K_g(\boldsymbol{\lambda}_j,z|\bs\omega)}. 
\end{align}
where the notation highlights the dependence on the mixture weights $\boldsymbol{\omega}$. These equations provide the link between the joint spectral distribution $F$ 
and the limiting spectral distribution of the weighted sample periodogram. Specifically, for each $g$, the weights $\bs\omega$ determine $M_g$ and $K_g$, and hence the limiting Stieltjes transform $\mathcal{S}_g$. 

\medskip

Theorem \ref{theorem1} together with intermediate results used in its proof also provides the corresponding empirical quantities. Let $\tilde{\mathbf{R}}_g(z)=(\tilde{\mathbf{S}}_g^{(n)}(z)-z\mathbf{I})^{-1}$ be the resolvent matrix and write the Stieltjes transform of $F_g^{(n)}$ as 
\begin{equation} \label{eq:s_hat}
\hat{\mathcal{S}}_g^{(n)}(z)=\frac{1}{n}\mathrm{trace}\big[\tilde{\mathbf{R}}_g(z)\big], \qquad z \in \mathbb{C}^+.
\end{equation}
Also define the Stieltjes kernel
\begin{equation} \label{eq:K_hat}
\hat{K}_g^{(n)}(z,\boldsymbol{\lambda})=\frac{1}{n}\mathrm{trace}\big[\tilde{\mathbf{R}}_g(z)\mathbf{W}_g\Psi^*(\boldsymbol{\lambda})\Psi(\boldsymbol{\lambda})\mathbf{W}_g\big],
\end{equation}
where $\Psi(\boldsymbol{\lambda})=\mathrm{diag}\{\psi(\boldsymbol{\lambda},\theta_1),\dots,\psi(\boldsymbol{\lambda},\theta_n)\}$ and $\psi$ is defined in \eqref{eq:def_psi_and_h}. Under the assumptions of Theorem~\ref{theorem1}, $\hat{\mathcal{S}}_g^{(n)}(z)\approx \mathcal{S}_g(z|\bs\omega^0)$ and $\hat{K}_g^{(n)}(z,\boldsymbol{\lambda})\approx K_g(z,\boldsymbol{\lambda}^0_j|\bs\omega^0)$ for all $\boldsymbol{\lambda}_j^0\in\Lambda^0$, with the exact meaning of $\approx$ given in the proof. These relations suggest estimating $\bs\omega^0$ by the weights $\hat{\bs\omega}$ producing the Stieltjes transform $\mathcal{S}_g(z|\hat{\bs\omega})$ closest to the empirical counterpart $\hat{\mathcal{S}}_g^{(n)}(z)$. 
To formalize, let 
\[
\mathcal{Z}=\{z_1,\ldots,z_D\}\subset\mathbb{C}^+
\]
be a prespecified finite set of evaluation points. 
Since the most informative frequencies for the estimation are generally unknown, select a finite collection $\mathcal{G}$ of weight functions $g$ satisfying Assumption~\textbf{A.5}. The estimated weights $\hat{\bs\omega}=(\hat\omega_1,\ldots,\hat\omega_J)$ are then obtained as solution to the following optimization problem:
\begin{align} \label{eq:optimization_problem}
\hat{\bs\omega}
=\argmin_{\bs\omega\in\Delta^J}
\sum_{z\in\mathcal{Z}}\sum_{g\in\mathcal{G}}\mathcal{L}\big(\hat{\mathcal{S}}^{(n)}_g(z), \mathcal{S}_g(z|\bs\omega)\big), 
\end{align}
where $\mathcal{L}$ is a suitable discrepancy measure. In this paper, we consider the $L_\kappa$-loss for $\kappa\geq 1$, and use $\kappa=1$ or $\kappa=2$ in the numerical and theoretical analyses. Putting the above together, the resulting estimator of the joint spectral distribution $F^0$ is 
\begin{equation}\label{eq:F_hat_A&Sigma}
    \hat F(\bs\lambda)= 
    \sum_{j=1}^{J}\hat{\omega}_j\mathbb{I}\{\boldsymbol{\lambda}_j^0\leq\boldsymbol{\lambda}\}.
\end{equation}
The estimator in \eqref{eq:F_hat_A&Sigma} can also be employed when the true joint spectral distribution is not a finite mixture of point masses. 
The fixed-grid formulation considered here leads to a parametric estimation problem and is used for the consistency theory developed in Section \ref{sec:consistency}. It should not be interpreted as requiring the underlying joint spectral distribution to have finite support. For a general $F^{A,\Sigma}$, the same procedure can be applied using a discrete approximation on a sufficiently fine grid. A formal treatment of this setting would require allowing the number of grid points to increase with $p$ and $n$. Developing such a sieve method is, however, beyond the scope of the present paper.

\medskip

The evaluation of the objective function in \eqref{eq:optimization_problem} requires the calculation of $\mathcal{S}_g(z|\bs\omega)$ for every $\bs\omega\in\Delta_J$, $g\in\mathcal{G}$ and $z\in\mathcal{Z}$, which in turn depends on the kernel values $K_g(z,\bs\lambda^0_j)$ for $j=1,\ldots,J$. We compute these quantities using the fixed-point equations in Theorem~\ref{theorem1}. The empirical kernel $\hat K_g^{(n)}(z,\bs\lambda^0_j)$ defined in \eqref{eq:K_hat} provides a natural initialization, since it converges pointwise to $K_g(z,\bs\lambda_j^0)$ under the given assumptions. This leads to the following algorithm.

\medskip

\hrule
\medskip

\noindent\textbf{Algorithm}

{\it Set $K_g^{(1)}(z,\boldsymbol{\lambda_j^0})=\hat{K}_g^{(n)}(z,\boldsymbol{\lambda_j^0})$. \\
\indent For  $i = 1,\ldots, I$ compute\\ 
$\indent \indent  \displaystyle M_g^{(i)}(z,\theta|\bs\omega)=\sum_{j=1}^{J}\omega_j\;\frac{g(\theta)h(\boldsymbol{\lambda}_j^0,\theta)}{1+K_g^{(i)}(\boldsymbol{\lambda}_j^0,z)}$, \\[.2cm]
$\indent \indent\displaystyle K_g^{(i+1)}(z,\boldsymbol{\lambda_j^0}) = \frac{1}{2\pi}\int_{0}^{2\pi}\frac{g(\theta)h(\boldsymbol{\lambda_j^0},\theta)}{cM_g^{(i)}(z,\theta|\bs\omega)-z}d\theta$. \\
\indent end.} \medskip

\hrule

\bigskip

A practical difficulty arises because the size of the grid $\Lambda^0$ can increase rapidly with the order of the time series under consideration. Consider, for example, an AR($q$) model
\[
X_t=A_1X_{t-1}+\cdots+A_qX_{t-q}+Z_t.
\]
If the marginal spectral distribution of each coefficient matrix and that of the innovation covariance are supported on $r$ points, then the unrestricted Cartesian grid has $r^{q+1}$ points. Direct optimization over the corresponding simplex can therefore become computationally demanding. One way to reduce the dimension is to assume that the joint spectral distribution factorizes into its marginal distributions. Under this product assumption, $F^{A,\Sigma}=F^{A_1}\cdots F^{A_q}F^{\Sigma}$ with joint support $\Lambda^0=\Lambda^{(1)}\times\cdots\times\Lambda^{(q+1)}$. The optimization can then be carried out over the marginal weight vectors
\[
\bs\omega^{(k)}=(\omega_1^{(k)},\ldots,\omega_{r_k}^{(k)}),\qquad k=1,\ldots,q+1,
\]
subject to $\sum_{i=1}^{r_k}\omega_i^{(k)}=1$ and $\omega_i^{(k)}\geq 0$. This reduces the number of mixture weights from $r_1\cdots r_{q+1}$ to $r_1+\cdots+r_{q+1}$. This reduction can be substantial for higher-order models, although it comes at the cost of imposing independence among the eigenvalue coordinates of the join spectral distribution.

\section{Estimation of the spectral density matrix}\label{sec:estimation_SDM}

Estimation of the joint spectral distribution $F^{A,\Sigma}$ does not by itself yield estimates of the individual eigenvalues of the coefficient matrices and the innovation covariance matrix. This differs from the iid setting, where the only relevant population matrix is the covariance matrix $\Sigma$. If $F^\Sigma$ denotes its limiting spectral distribution and $\hat F^\Sigma$ is an estimator thereof, the eigenvalues $\sigma_1,\ldots,\sigma_p$ of $\Sigma$ can, for example be estimated by the corresponding quantiles:
\[
\hat\sigma_j=(\hat F^\Sigma)^{-1}\bigg(1-\frac j{p+1}\bigg),\qquad j=1,\ldots,p. 
\]
Such a construction relies on the natural ordering of the eigenvalues of a single matrix. In the present setting, however, $F^{A,\Sigma}$ is a multivariate distribution describing the joint behavior of the eigenvalue vectors
\[
\bs\lambda_j = (\lambda_{p,j}, \sigma_j),
\qquad j=1,\ldots, p.
\]
There is no natural ordering of these vectors. Consequently, estimation of ${F}^{A,\Sigma}$ alone does not determine how its estimated support points should be paired with the eigenvectors of the coefficent matrices and $\Sigma$. To construct an estimator of the spectral density matrix, we therefore supplement $\hat F$ with an estimate $\hat{\mbf{U}}$ of the common eigenbasis and introduce an ordering of the support points based on a scalar summary of their spectral densities. The estimated eigenbasis is then used to represent the spectral density matrix at all frequencies. 

Note that the estimation of $F^{A\Sigma}$ does not require the estimation of the common eigenbasis. Reconstruction of the full spectral density matrix, however, is a more demanding problem because it additionally requires associating the estimated spectral components with directions in the observation space. We therefore regard the following procedure as a reconstruction method conditional on an adequate estimate of the common eigenspaces being available.
\medskip

\noindent\textbf{Estimation procedure: }
\vspace{-5pt}
\begin{enumerate}
	\item 
	Choose a positive-valued weight function $g_0$ on $[0,2\pi]$ and compute the weighted sample periodogram matrix $\mbf{S}_{g_0}^{(n)}$. Its population counterpart is the weighted integrated spectral density matrix. In particular, when $g_0 \equiv 1$, the latter reduces to the integrated spectral density matrix. 
	
	\item
	Let $\hat{\mbf{U}}$ be the matrix of eigenvectors of $\mbf{S}_{g_0}^{(n)}$ arranged according to decreasing eigenvalues. We use $\hat{\mbf{U}}$ as an estimate of the common eigenbasis. 
	
	\item 
	Since $\hat{F}$ is a discrete mixture, we represent it by the pairs 
    \[
    (\bs\lambda_j^0, \hat\omega_j),
    \qquad j=1,\ldots,J,
    \]
    where $\bs\lambda_j^0$ and $\hat\omega_j$, $j=1,\ldots,J$, are the grid points and the corresponding estimated weights, respectively.
	
	\item
	For each grid point define
    \[
	m_j(g_0) =	\frac{1}{2\pi}\int g_0(\theta)h(\bs\lambda_j^0,\theta)d\theta.
	\]
    Order the grid points by decreasing values of $m_j(g_0)$ and denote the resulting ordering by 
    \[
    \bs\lambda_{(j)}^0,
    \qquad j=1,\ldots,J.
    \]
	The quantities $m_j(g_0)$ are the population eigenvalues associated with the weighted integrated spectral density matrix and therefore provide a natural scalar criterion for pairing the support points with the ordered empirical eigenvectors. Ties may be broken arbitrarily.
	
	\item
    Convert the estimated mixture weights into multiplicities by setting
	\[
    \hat{p}_j = \lfloor p \hat\omega_{j}\rfloor,
    \qquad j=1,\ldots,J,
    \]
    where $\lfloor\cdot\rfloor$ denotes integer part. The multiplicities are arranged according to the ordering in the preceding step and adjusted, if necessary, so that $\sum_j \hat p_j=p$. This adjustment can be made, for example, by pruning sufficiently small estimated weights and reallocating the resulting rounding discrepancy. 
	
	\item
	Define the estimator of the spectral density matrix at frequency $\theta$ by
	\[
	\hat{\mbf{H}}(\theta) = \hat{\mbf{U}}~\mbox{diag}\Big(\big(h(\bs\lambda_j^0,\theta)\bs I_{\hat{p}_{(j)}}\big)_{j=1}^J\Big) \hat{\mbf{U}}^*. 
	\]
	{\color{black} where $h(\bs\lambda_j^0,\theta)$ is as defined in (\ref{eq:def_psi_and_h}) and $\bs I_{\hat{p}_{(j)}}$ denotes the $\hat{p}_{(j)}\times\hat{p}_{(j)}$ identity matrix.}
\end{enumerate}

By construction, $\mbf{H}(\theta)$ is symmetric (Hermitian) for every $\theta$, and the estimated spectral density matrices share the common eigenbasis $\hat{\mbf{U}}$ across frequencies. The estimator depends on the choice of the target weight function $g_0$, since $g_0$ determines both the estimated eigenbasis and the ordering of the support points. One could therefore construct estimators using several choices of $g_0$ and combine them trhough model selection or model averaging. Alternatively, the ordering of the support points could be chose by searching over plausible permutations and selecting the one providing the best fit to the observed data. We do not pursue these extensions here.

\section{Consistency under the squared loss}
\label{sec:consistency}

In this section, we establish consistency of the estimator introduced in Section~\ref{sub:SpecEst_Algm} for the $L^2$-discrepancy criterion. We assume that the true joint spectral distribution $F^0$ is supported on the known grid $\Lambda^0=\{\bs\lambda_1^0,\ldots,\bs\lambda_J^0\}$, so that the unknown parameter is the vector of population weights $\bs\omega^0=(\omega_1^0,\ldots,\omega_J^0)$. Note that $\bs\omega^0$ may lie on the boundary of $\Delta_J$.
%
%
For $\kappa=2$, the estimator defined in Section~\ref{sub:SpecEst_Algm} minimizes the squared discrepancy between the empirical Stieltjes transforms and those associated with the candidate distributions supported on $\bs\Lambda^0$. The optimization problem in \eqref{eq:optimization_problem} becomes
\begin{equation}\label{eq:omega_hat_L2_estimator}
\hat{\bs\omega}=\arg\min_{\bs\omega\in\Delta_J}\sum_{g\in\mcal G}\sum_{z\in\mcal Z}\left|\hat{\mcal S}_g^{(n)}(z)-
\mathcal{S}_g(z|\omega)\right|^2.
\end{equation}
Our goal is to give conditions under which this $L^2$-estimator satisfies $\hat{\bs\omega} \stackrel{a.s.}{\longrightarrow} \bs\omega^0$.


\subsection{Sufficient condition for consistency}
\label{subsec:consistency_condition}


The essential point for the consistency argument is that the empirical criterion converges to its population counterpart, which is minimized at the true weight vector. Consistency therefore reduces to establishing sufficient local identification of $\bs\omega^0$.
To formulate this identification condition, we parameterize perturbations of $\bs\omega^0$ within the simplex. Let $\mbf B$ be any fixed $J\times(J-1)$ matrix of rank $J-1$ satisfying \[ \mbf B^T\mbf 1_J=0. \] Since any two elements of $\Delta_J$ have coordinates summing to one, $\bs\omega-\bs\omega^0$ belongs to the $(J-1)$-dimensional subspace orthogonal to $\mbf 1_J$. The columns of $\mbf B$ form a basis for this subspace, and hence every $\bs\omega\in\Delta_J$ can be represented as \[ \bs\omega = \bs\omega^0+\mbf B\bs\eta \] for some $\bs\eta\in\mbb R^{J-1}$, subject to $\bs\omega^0+\mbf B\bs\eta\in\Delta_J$.

For $g\in\mcal G$ and $z\in\mcal Z$, define 
\[ 
\mbf a_g^0(z) = \frac{1}{2\pi} \int_0^{2\pi} \frac{ g(\theta)\mbf v_g^0(z,\theta) }{ \{cM_g^0(z,\theta)-z\}^2 }\,d\theta, 
\] 
where 
\[ 
\mbf v_g^0(z,\theta) = \left( \frac{ h(\bs\lambda_j^0,\theta) }{ 1+K_g^0(\bs\lambda_j^0,z) } \right)_{j=1}^J. 
\]
We then define the $J\times J$ matrix \begin{equation}\label{eq:M_G_Z_Lambda} 
\mcal M_{\mcal G,\mcal Z} (\bs\Lambda^0,\bs\omega^0) = \frac{1}{|\mcal G||\mcal Z|} \sum_{g\in\mcal G} \sum_{z\in\mcal Z} \mbf a_g^0(z)\mbf a_g^0(z)^*. 
\end{equation} 
Thus, $\mcal M_{\mcal G,\mcal Z}$ is a Gram-type matrix and is automatically positive semidefinite. Its restriction to the tangent space of the simplex describes whether local perturbations of the population weights can be distinguished through the collection of limiting Stieltjes transforms. This leads to the following sufficient condition for consistency.

\begin{theorem}
\label{thm:sufficient_cond_L2_consistency} 
Suppose that Assumptions~\textbf{A.1}--\textbf{A.5} hold. Provided that $|\mathcal{Z}|\geq J$ and the elements of $\mathcal{Z}$ are distinct, a sufficient condition for the consistency of $\hat{\bs\omega}$ is that 
\[ 
\mbf B^T \mcal M_{\mcal G,\mcal Z} (\bs\Lambda_J^0,\bs\omega^0) \mbf B 
\]
is positive definite. 
\end{theorem}

The positive-definiteness condition in Theorem~\ref{thm:sufficient_cond_L2_consistency} has a natural identifiability interpretation. For $\bs\eta\in\mbb R^{J-1}$, 
\[ 
\bs\eta^T \mbf B^T \mcal M_{\mcal G,\mcal Z} \mbf B \bs\eta = \frac{1}{|\mcal G||\mcal Z|} \sum_{g\in\mcal G} \sum_{z\in\mcal Z} \left| \mbf a_g^0(z)^*\mbf B\bs\eta \right|^2.
\]
Hence, positive definiteness requires that every nonzero admissible perturbation of the population weights changes at least one of the quantities used in the $L^2$-criterion. Equivalently, the collection of weight functions in $\mcal G$ and evaluation points in $\mcal Z$ must provide enough information to locally distinguish $\bs\omega^0$ from neighboring candidate distributions. 
The proof makes this connection explicit. If $D(\bs\omega)$ denotes the population $L^2$ criterion, then a Taylor expansion around $\bs\omega^0$ yields, for $\bs\omega=\bs\omega^0+\mbf B\bs\eta$, \[ D(\bs\omega) = \bs\eta^T \mbf B^T \mcal M_{\mcal G,\mcal Z} (\bs\Lambda_J^0,\bs\omega^0) \mbf B \bs\eta + O(\|\bs\eta\|^3). \] Thus, the condition in Theorem~\ref{thm:sufficient_cond_L2_consistency} ensures that the population criterion has a strict local minimum at $\bs\omega^0$. Together with uniform convergence of the empirical criterion to its population counterpart, this yields consistency of the estimator.

\subsection{Special cases}\label{subsec:special_cases} We next illustrate the identifiability condition in Theorem~\ref{thm:sufficient_cond_L2_consistency} in several important special cases. The examples also clarify the role of the weight functions in $\mcal G$. For independent observations, a single constant weight function is sufficient to distinguish the population weights. For dependent processes, by contrast, information from different frequency bands may be needed to distinguish the spectral behavior associated with different support points.

\subsubsection{Independent observations}
\label{subsec:indpendence} 
Suppose first that the observations are independent, so that 
\[ 
X_t=\mbf A_0 Z_t, \qquad \Sigma=\mbf A_0^2=\operatorname{Var}(X_t).
\] 
In this case, the joint spectral distribution reduces to the limiting spectral distribution of the covariance matrix, 
\[ 
F^0\equiv F^\Sigma = \sum_{j=1}^J \omega_j^0\delta_{\lambda_j^0}, 
\] 
where $\lambda_1^0,\ldots,\lambda_J^0$ are distinct nonnegative real numbers. Since there is no temporal dependence, frequency localization provides no additional information, and it is sufficient to take 
\[ 
\mcal G=\{g_0\}, \qquad g_0(\theta)\equiv1. 
\] 
The following result shows that, in this case, the population weights are locally identifiable from the limiting Stieltjes transforms.

\begin{proposition}
\label{prop:L2_consistency_independent_observations}
Let $\mcal G=\{g_0\}$ with $g_0(\theta)\equiv1$. Provided that $|\mathcal{Z}|\geq J$ and the elements of $\mathcal{Z}$ are distinct, 
\[ 
\mbf B^T \mcal M_{\mcal G,\mcal Z} (\bs\Lambda^0,\bs\omega^0) \mbf B 
\] 
is positive definite. 
\end{proposition}

Thus, for independent observations, distinct support points together with sufficiently many distinct evaluation points in $\mcal Z$ are enough to identify perturbations of the population weights. By Theorem~\ref{thm:sufficient_cond_L2_consistency}, the $L^2$-estimator is therefore consistent under the stated conditions.

\subsubsection{ARMA processes}\label{subsec:ARMA} For a dependent linear process, a constant weight function need not provide enough information to distinguish the spectral contributions of different support points. The relevant information is now distributed across frequencies, which motivates the use of a collection of localized weight functions. To formulate a sufficient condition, define \begin{equation}\label{eq:Gram_spectral_density} \mbf G(\bs\Lambda^0) = \frac{1}{2\pi} \int_0^{2\pi} \mbf h^0(\theta) \{\mbf h^0(\theta)\}^T\,d\theta, \end{equation} where 
\[
\mbf h^0(\theta) = \bigl( h(\bs\lambda_1^0,\theta), \ldots, h(\bs\lambda_J^0,\theta) \bigr)^T. 
\]
The matrix $\mbf G(\bs\Lambda^0)$ is the Gram matrix of the spectral profiles 
\[ 
h(\bs\lambda_1^0,\cdot),\ldots, h(\bs\lambda_J^0,\cdot) 
\] 
in $L^2([0,2\pi])$. Consequently, its positive definiteness means that these spectral profiles are linearly independent. In particular, the different support points generate sufficiently distinct frequency-domain behavior to be distinguished through appropriately chosen weight functions.

This condition also has a time-domain interpretation. Let 
\[
\bs\gamma_\ell^0 = \bigl( \gamma_\ell(\bs\lambda_1^0), \ldots, \gamma_\ell(\bs\lambda_J^0) \bigr)^T, \qquad \ell\geq0, 
\] 
where $\gamma_\ell(\bs\lambda_j^0)$ denotes the lag-$\ell$ autocovariance of the one-dimensional linear process associated with $\bs\lambda_j^0$. By the spectral representation, 
\[ 
h(\bs\lambda_j^0,\theta) = \gamma_0(\bs\lambda_j^0) + 2\sum_{\ell=1}^\infty \gamma_\ell(\bs\lambda_j^0)\cos(\ell\theta), \qquad \theta\in[0,2\pi], 
\]
and hence \begin{equation}
\label{eq:Gram_spectral_representation} 
\mbf G(\bs\Lambda^0) = \bs\gamma_0^0(\bs\gamma_0^0)^T + 2\sum_{\ell=1}^\infty \bs\gamma_\ell^0(\bs\gamma_\ell^0)^T. 
\end{equation} 
Thus, positive definiteness of $\mbf G(\bs\Lambda^0)$ can equivalently be interpreted as requiring the support points to generate sufficiently different autocovariance sequences. The following proposition connects this condition to the local identifiability condition in Theorem~\ref{thm:sufficient_cond_L2_consistency}. 

\begin{proposition}
\label{prop:ARMA_L2_consistency_condition} Suppose that $(X_t)$ satisfies the assumptions of Theorem~\ref{thm:sufficient_cond_L2_consistency}. If $\mbf G(\bs\Lambda^0)$ is positive definite, then there exists a collection $\mcal G$ of weight functions such that 
\[ 
\mbf B^T \mcal M_{\mcal G,\mcal Z} (\bs\Lambda^0,\bs\omega^0) \mbf B 
\] 
is positive definite. 
\end{proposition}

The proposition shows why several frequency-localized weight functions are useful for dependent processes. Linear independence of the spectral profiles ensures that there are frequency regions at which different support points can be distinguished. By choosing the functions in $\mcal G$ to concentrate on such regions, the resulting collection of limiting Stieltjes transforms provides sufficient local information to identify the population weights.

\subsubsection{AR(1) processes} 

We finally specialize the preceding condition to an AR$(1)$ process, where positive definiteness of $\mbf G(\bs\Lambda^0)$ can be verified directly. Let 
\[ 
X_t = \mbf A X_{t-1} + \Sigma^{1/2}Z_t, 
\] 
where $\mbf A$ and $\Sigma$ are symmetric/Hermitian, $\|\mbf A\|<1$, and there exists an orthogonal/unitary matrix $\mbf U$ such that 
$
\mbf U^*\mbf A\mbf U = 
\operatorname{diag}(\alpha_1,\ldots,\alpha_p)$ and $\mbf U^*\Sigma\mbf U 
= \operatorname{diag}(\sigma_1^2,\ldots,\sigma_p^2)$. Let $\bs\lambda=(\alpha,\sigma)$ and suppose that the true joint spectral distribution is 
\[ 
F^0 = \sum_{j=1}^J \omega_j^0\delta_{\bs\lambda_j^0}, 
\qquad \bs\lambda_j^0 = (\alpha_j^0,\sigma_j^0), 
\] 
where $\max_{j}|\alpha_j^0|<1$ and $\min_{j}\sigma_j^0>0$. For the univariate AR$(1)$ process corresponding to $\bs\lambda_j^0$, the autocovariance function is 
\[ 
\gamma_\ell(\bs\lambda_j^0) 
= \frac{(\sigma_j^0)^2}{1-(\alpha_j^0)^2} (\alpha_j^0)^\ell, 
\qquad \ell\geq0. 
\] 
Thus, different autoregressive coefficients generate different geometric autocovariance sequences. This leads to the following simple sufficient condition. 

\begin{proposition}
\label{prop:L2_consistency_AR1} If $\alpha_1^0,\ldots,\alpha_J^0$ are distinct, then $\mbf G(\bs\Lambda^0)$ defined in (\ref{eq:Gram_spectral_density}) is positive definite. 
\end{proposition}
The proposition follows from the fact that the vectors formed by the first $J$ autocovariances are, after rescaling, the columns of a Vandermonde-type matrix generated by $\alpha_1^0,\ldots,\alpha_J^0$. Distinct autoregressive coefficients therefore imply linear independence of the corresponding spectral profiles. Combining Proposition~\ref{prop:L2_consistency_AR1} with Proposition~\ref{prop:ARMA_L2_consistency_condition} and Theorem~\ref{thm:sufficient_cond_L2_consistency} establishes consistency of the $L^2$-estimator for the AR$(1)$ model under the stated conditions. The same argument can be extended to autoregressive processes of arbitrary finite order; details are provided in the Supplementary Material.

\section{Model selection}
\label{sub:SpecEst_ModelSeln} 

In practice, the order of the underlying linear process is typically unknown. We therefore consider a finite collection of candidate models and propose a bootstrap procedure for assessing their fit and selecting among them. The basic principle is that, if a fitted candidate model provides an adequate description of the data, samples generated from that model should reproduce relevant spectral features of the observed sample. Let $\mathcal{M}_1,\ldots,\mathcal{M}_M$ denote the candidate models. For each $\mathcal{M}_m$, we first estimate its joint spectral distribution using the procedure of Section~\ref{sub:SpecEst_Algm}. The resulting fitted model is then used to generate bootstrap samples. Model fit is assessed by comparing a chosen matrix-valued statistic computed from each bootstrap sample with the corresponding statistic computed from the observed data. More precisely, let $\mbf H_s$ denote a statistic computed from the observed sample, such as the lag-$\tau$ symmetrized sample autocovariance matrix, and let $\mbf H_{m,b}^*$ denote its counterpart computed from the $b$th bootstrap sample generated under candidate model $\mathcal{M}_m$. Since the methodology developed in this paper focuses on spectral distributions rather than individual matrix entries, we compare these matrices through their ordered eigenvalues. For a Hermitian matrix $\mbf H$, let 
\[
\operatorname{eig}(\mbf H) = \bigl(\lambda_1(\mbf H),\ldots,\lambda_p(\mbf H)\bigr)^T, 
\] 
where the eigenvalues are arranged in decreasing order.

For each candidate model $\mathcal{M}_m$, proceed as follows: 
\begin{enumerate} 
\item Fit $\mathcal{M}_m$ to the observed data and obtain the estimated joint spectral distribution of its coefficient matrices and innovation covariance matrix. 
\item Construct a bootstrap data-generating process from the fitted model. For $b=1,\ldots,B$, generate a bootstrap sample 
\[ 
X_{t}^{*(m,b)} = \sum_{\ell=0}^{\infty} \mbf A_{\ell,m}^* Z_{t-\ell}^{*(m,b)}, 
\] 
where the coefficient matrices and innovation covariance are constructed from the fitted spectral distribution under $\mathcal{M}_m$, and the bootstrap innovations are Gaussian with the fitted covariance structure. 
\item From the $b$th bootstrap sample, compute $\mbf H_{m,b}^*$ and define the bootstrap discrepancy 
\[ 
L_{m,b} = \big\| \operatorname{eig}(\mbf H_{m,b}^*) - \operatorname{eig}(\mbf H_s) \big\|_2^2. 
\]
\item Define the model-selection score 
\[ 
\bar L_m = \frac{1}{B} \sum_{b=1}^B L_{m,b}. 
\] 
\end{enumerate} 
The selected model is 
\[ \hat m = \arg\min_{1\leq m\leq M} \bar L_m. 
\]
Thus, the procedure favors candidate models whose fitted bootstrap distributions reproduce the empirical eigenvalue distribution of the chosen sample statistic. The statistic $\mbf H_s$ can be selected to emphasize features of particular interest. For example, using symmetrized sample autocovariance matrices at positive lags focuses the comparison on temporal dependence, whereas the lag-zero matrix also reflects the contemporaneous covariance structure. Several lags can likewise be combined in a single criterion if a broader assessment of the dependence structure is desired. The proposed criterion should be viewed as a bootstrap measure of relative fit among the candidate models rather than as a formal hypothesis test. Its performance depends on the choice of the diagnostic statistic $\mbf H_s$ and on whether the collection of candidate models contains an adequate approximation to the data-generating process.

\section{Simulation study}
\label{sec:simulation} 

\subsection{Estimation of the spectral distribution} 

We investigate the finite-sample performance of the proposed estimator through a series of simulation experiments. The study is designed to assess the effects of the dimension $p$, the sample size $n$, the aspect ratio $p/n$, the underlying linear process, and the number of weight functions used in the estimation procedure. Unless stated otherwise, the simulations use the following implementation settings. The weight functions $g\in\mathcal G$ are shifted upward by $0.05$. This avoids numerical instability in regions where $g(\theta)$ is close to zero, since $M_g(z,\theta)=0$ there and the corresponding integrand in $\mathcal S_g(z)$ reduces to $-1/z$. Four fixed-point iterations are used in computing $K_g(z, \bs\lambda)$. The set $\mathcal Z\subset\mathbb C^+$ is a rectangular grid whose real parts consist of 25 equally spaced points from $-2$ to $2$ and whose imaginary parts consist of five equally spaced points from
$0.1$ to $0.5$.  The innovations are iid centered Gaussian random vectors with covariance matrix $\Sigma$.

For estimation, we primarily use the $L^1$-discrepancy 
\[ 
\mathcal L_1 = \sum_{g\in\mathcal G} \sum_{z\in\mathcal Z} \Big| \hat{\mathcal S}_g^{(n)}(z) - \mathcal S_g(z) \Big|. 
\] 
We also compare its performance with the $L^2$-discrepancy used in the consistency analysis of Section~\ref{sec:consistency}. The simulation study consists of two stages. The first focuses on AR$(1)$ processes and examines, in particular, the effect of the number of weight functions. The second considers the more complex ARMA$(1,1)$ and AR$(2)$ models and investigates performance over a wider range of $(p,n)$ combinations.

\medskip

\noindent {\it AR(1) models.} We first consider the model 
\[ 
X_t=\mbf A_1X_{t-1}+Z_t, 
\] 
where the spectral distributions of $\mbf A_1$ and $\Sigma$ are finite mixtures of point masses. For example, $F^{A_1} = 0.5\,\delta_{-0.5} + 0.5\,\delta_{0.8}$ means that $\mbf A_1 = \mbf U\mbf D_{A_1}\mbf U^*$, where one half of the diagonal entries of $\mbf D_{A_1}$ are $-0.5$ and the other half are $0.8$. Similarly, $F^{A_1}=\delta_{0.5}$ corresponds to $\mbf A_1=0.5\mbf I_p$.

We consider the following five specifications: \begin{align*} 
\text{Case 1.1:}\qquad\qquad& F^{A_1}=\delta_{0.5}, & F^\Sigma &= 0.5\,\delta_1+0.5\,\delta_2, \\ 
\text{Case 1.2:}\qquad\qquad& F^{A_1}=\delta_{0.5}, & F^\Sigma &= 0.75\,\delta_1+0.25\,\delta_2, \\ 
\text{Case 1.3:}\qquad\qquad& F^{A_1}
= 0.5\,\delta_{-0.5}+0.5\,\delta_{0.8}, & F^\Sigma &= \delta_1, \\ 
\text{Case 1.4:}\qquad\qquad& F^{A_1} = 0.25\,\delta_{-0.5}+0.75\,\delta_{0.8}, & F^\Sigma &= \delta_1, \\ 
\text{Case 1.5:}\qquad\qquad& F^{A_1} = 0.5\,\delta_{-0.5}+0.5\,\delta_{0.8}, & F^\Sigma &= 0.5\,\delta_1+0.5\,\delta_2. 
\end{align*}
For each model, we consider $(p,n)=(400,1600)$ and $(600,2400)$ and three choices of $\mathcal G$, consisting of 4, 8, and 12 B-spline weight functions. For each configuration, 20 independent samples are generated. The candidate grid is prespecified and contains the true support points. We assess estimation accuracy using the $L^2$-distance between the true and estimated marginal spectral distribution functions, 
\[ 
d_{L^2}(F^0,\hat F) = \left\{ \int \bigl(F^0(x)-\hat F(x)\bigr)^2\,dx \right\}^{1/2}. 
\] 
The mean, median, and standard deviation of this distance across simulation replications are reported in Tables~\ref{table:case1.L1_AR} and \ref{table:case1.L1_Sig}.

The results show a substantial improvement when the number of B-spline weight functions is increased from 4 to 8. Increasing the number from 8 to 12 generally yields a smaller additional improvement relative to the associated increase in computational cost. We also repeated the experiments using the $L^2$-discrepancy 
\[ 
\sum_{g\in\mathcal G} \sum_{z\in\mathcal Z} \Big| \hat{\mathcal S}_g^{(n)}(z) - \mathcal S_g(z) \Big|^2 
\] 
in the estimation criterion. In these experiments, the $L^1$-criterion generally produced more accurate estimates, and we therefore use it for the remaining simulations. Results for the $L^2$-criterion are reported in Section~B.2.2 of the Supplementary Material.

\medskip

\begin{figure}[H] \label{fig:ARMA_AR}
\begin{center}
\resizebox{\textwidth}{!}{%
\begin{tabular}{cc}
$(p,n)=(400,1600)$ & $(p,n)=(400,800)$ \\ \vspace{-20pt}
\includegraphics[width=2.9in, height=3in, bb = 50 50 650 650]{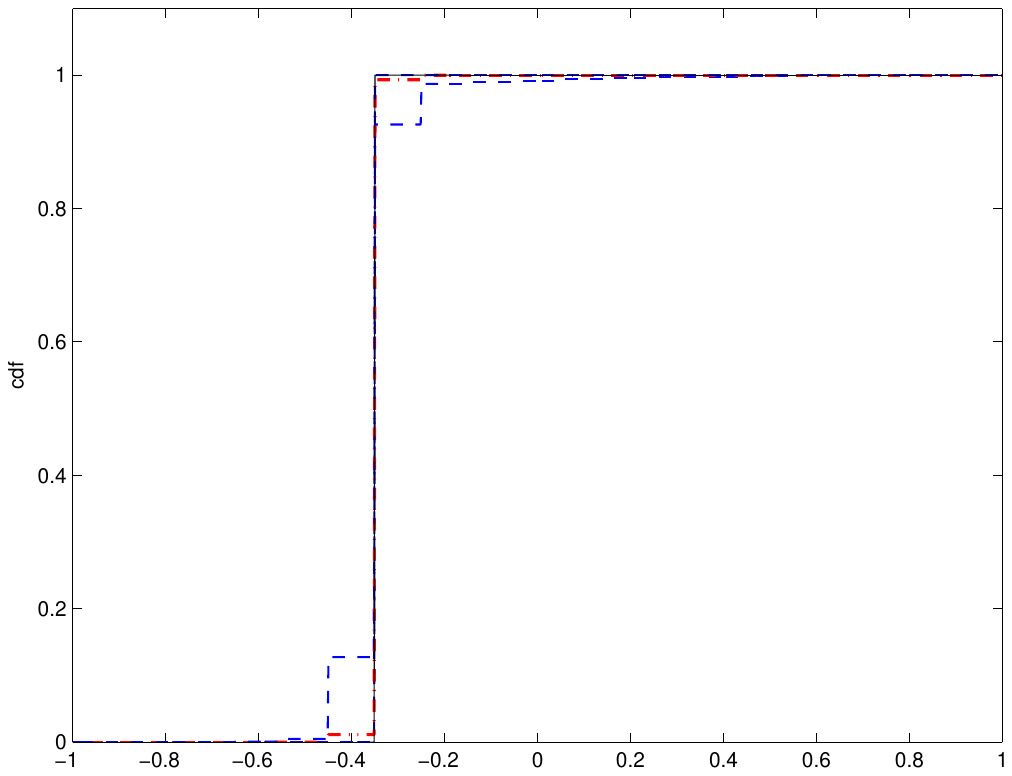} 
& \includegraphics[width=2.9in, height=3in, bb = 50 50 650 650]{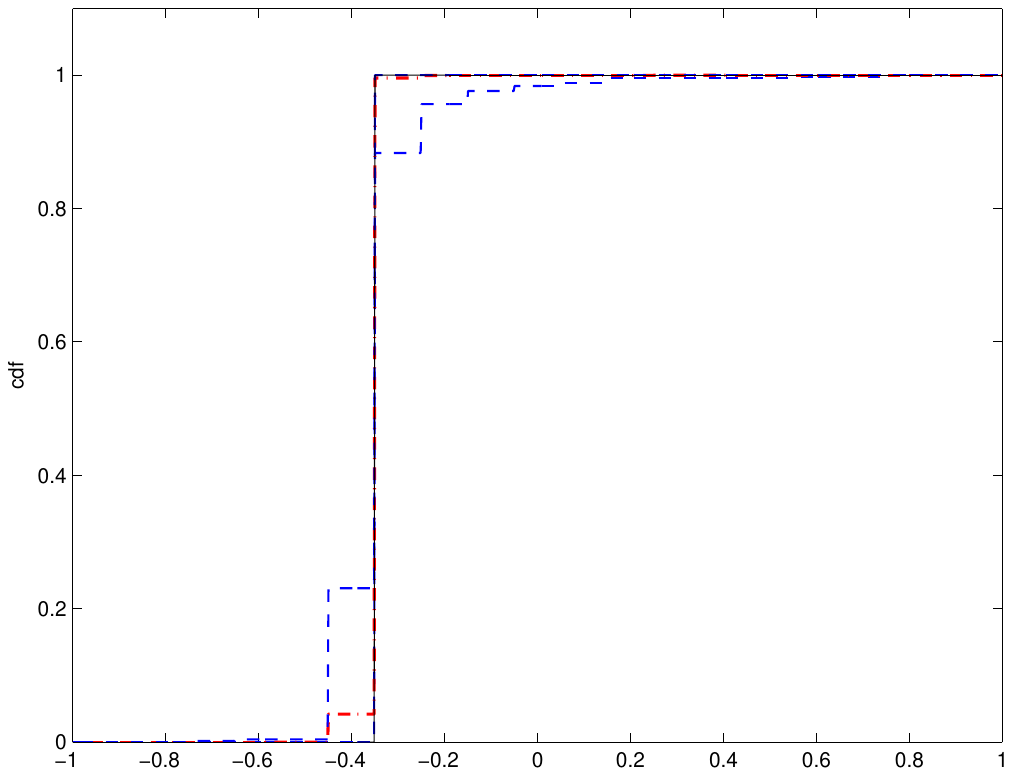}\\ \vspace{-10pt}
$(p,n)=(200,800)$ & $(p,n)=(200,400)$ \vspace{-5pt}\\ \vspace{-10pt}
\includegraphics[width=2.9in, height=3in, bb = 50 50 650 650]{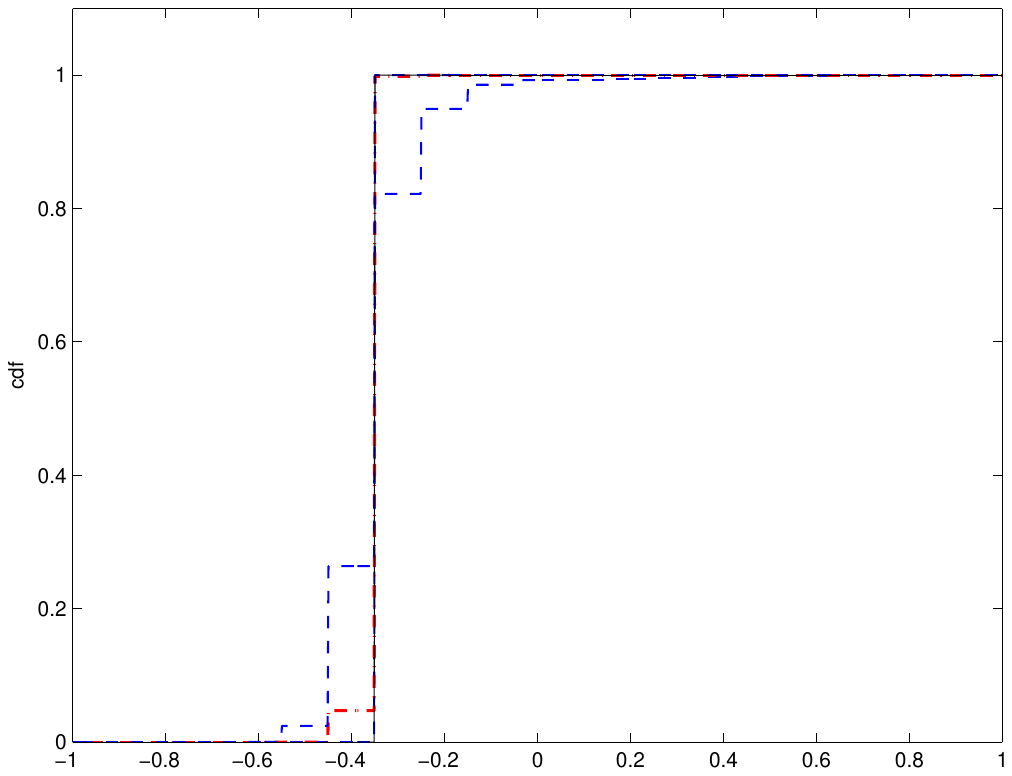}
& \includegraphics[width=2.9in, height=3in, bb = 50 50 650 650]{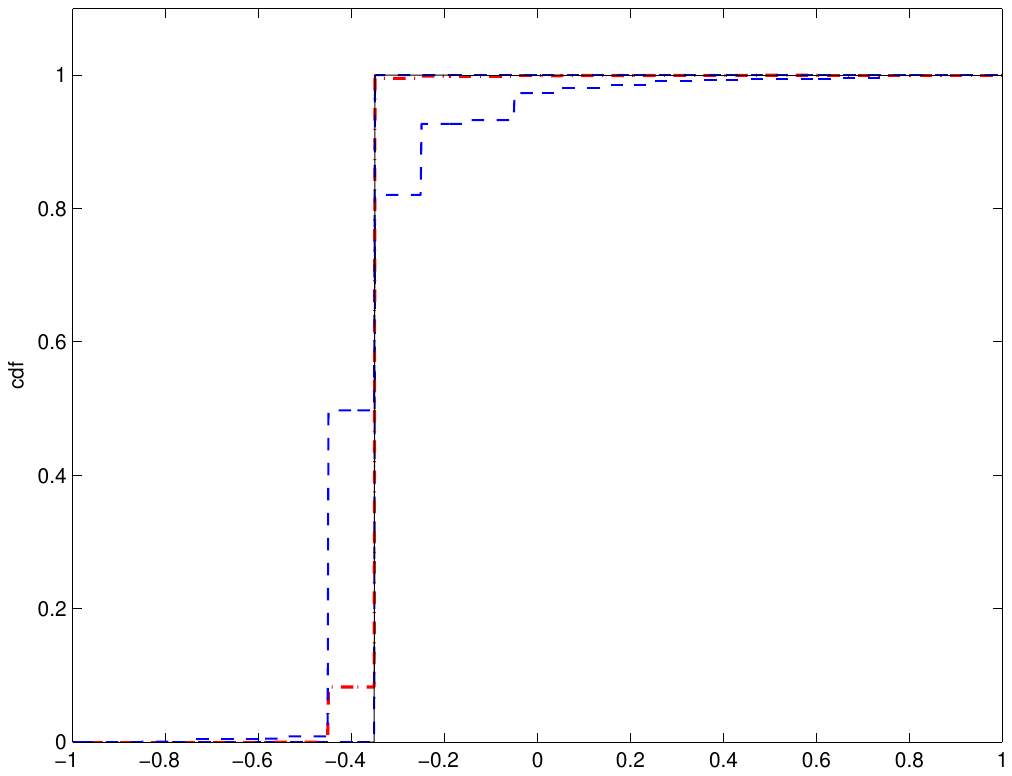} \vspace{-20pt}\\ \vspace{-20pt}
$(p,n)=(200,100)$ & \vspace{-10pt}\\ \vspace{-20pt}
\hspace{10pt}\includegraphics[width=2.8in, height=2.8in, bb = 50 50 650 650]{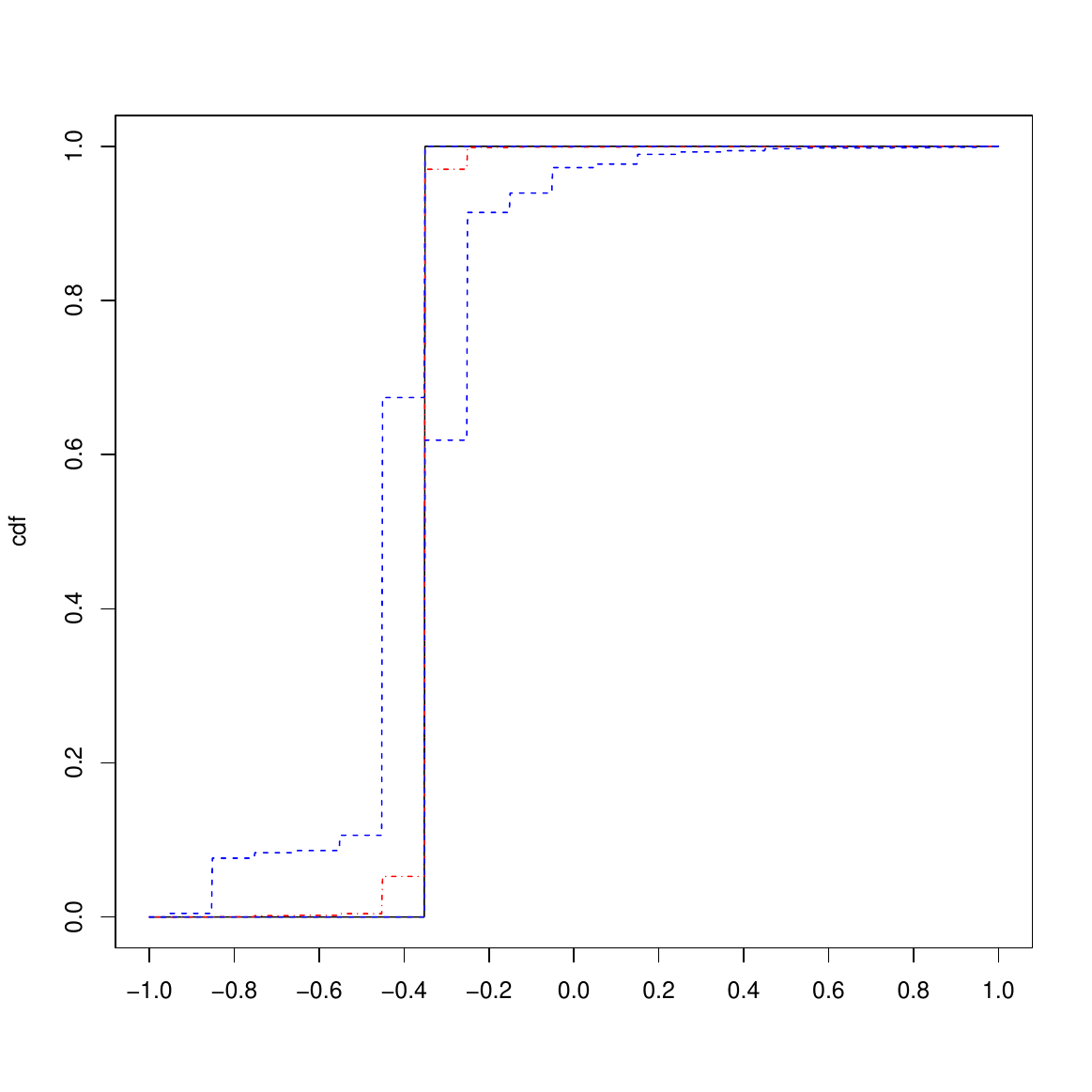}
& \\
\end{tabular}
}
\vspace{15pt}
         \caption{Plot of Median and $90\%$ confidence band for spectral cdf of \textbf{AR coefficient matrix} corresponding to the \textbf{case 2.1}. \textbf{Dash-Dot Red} curve: median, \textbf{Dashed Blue} curve: $90\%$ confidence band, \textbf{Black Solid} curve: true spectral cdf}
\end{center}
\end{figure}

\begin{figure}[H] \label{fig:ARMA_MA}
\begin{center}
\resizebox{\textwidth}{!}{%
\begin{tabular}{cc}
$(p,n)=(400,1600)$ & $(p,n)=(400,800)$ \\ \vspace{-20pt}
\includegraphics[width=2.9in, height=3in, bb = 50 50 650 650]{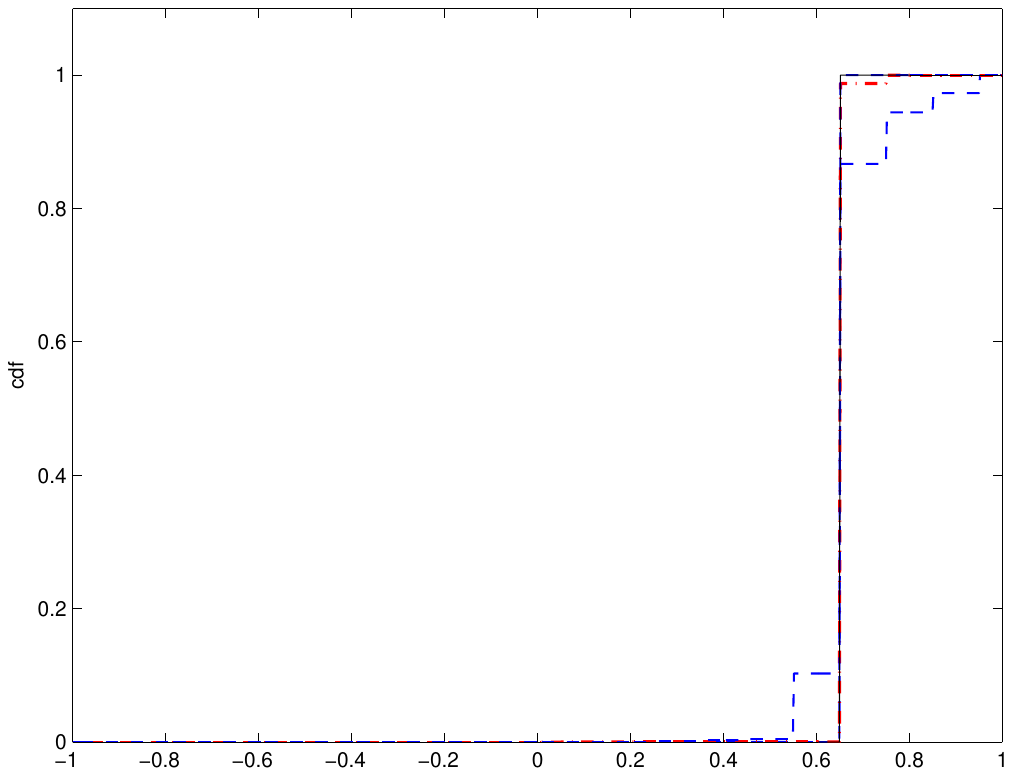}
& \includegraphics[width=2.9in, height=3in, bb = 50 50 650 650]{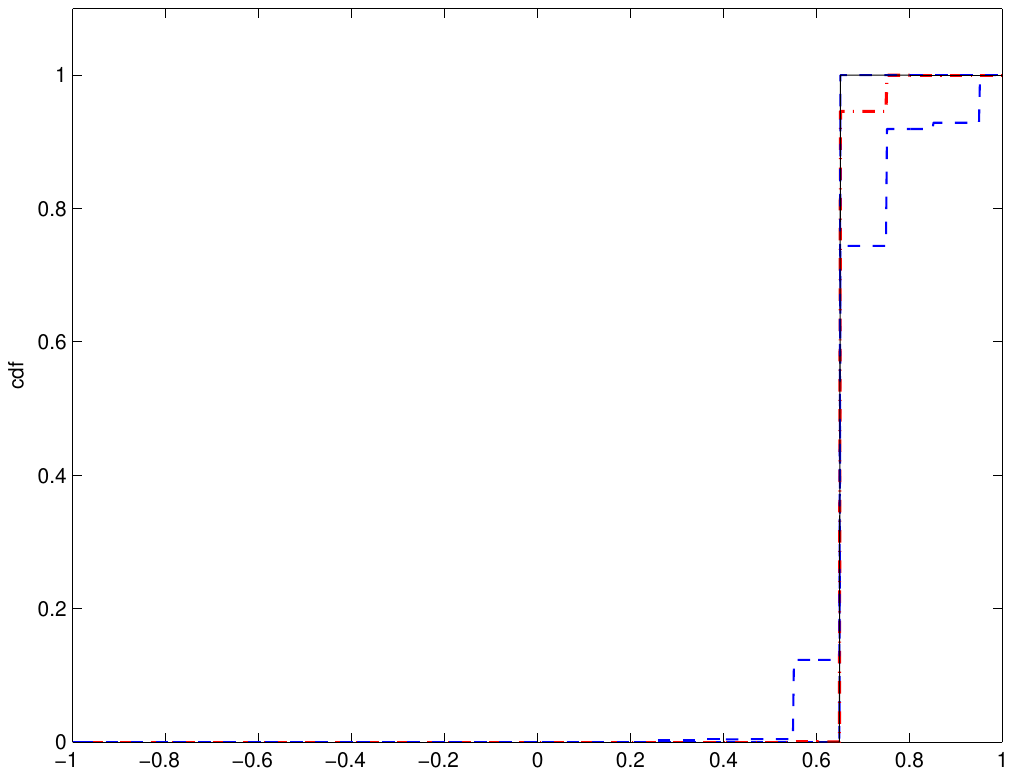}\vspace{-20pt}\\ \vspace{-10pt}
$(p,n)=(200,800)$ & $(p,n)=(200,400)$ \vspace{-5pt}\\ \vspace{-10pt}
\includegraphics[width=2.9in, height=3in, bb = 50 50 650 650]{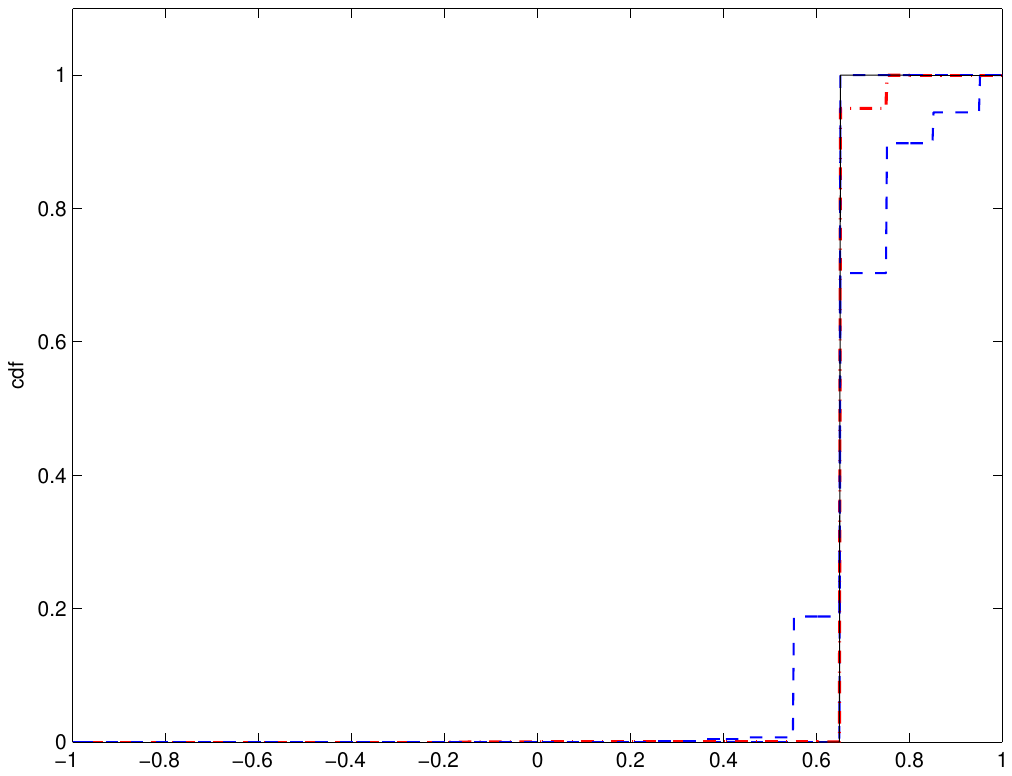}
& \includegraphics[width=2.9in, height=3in, bb = 50 50 650 650]{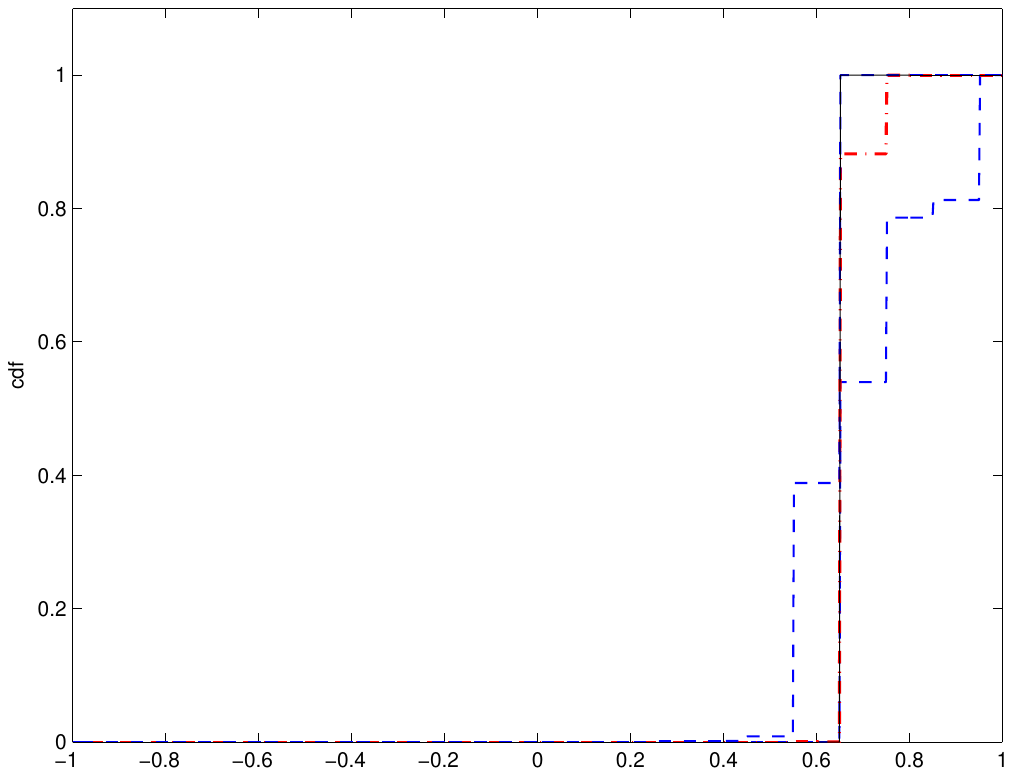}\vspace{-20pt}\\ \vspace{-20pt}
$(p,n)=(200,100)$ & \vspace{-0pt}\\ \vspace{-0pt}
\hspace{10pt}\includegraphics[width=2.9in, height=2.9in, bb = 50 50 650 650]{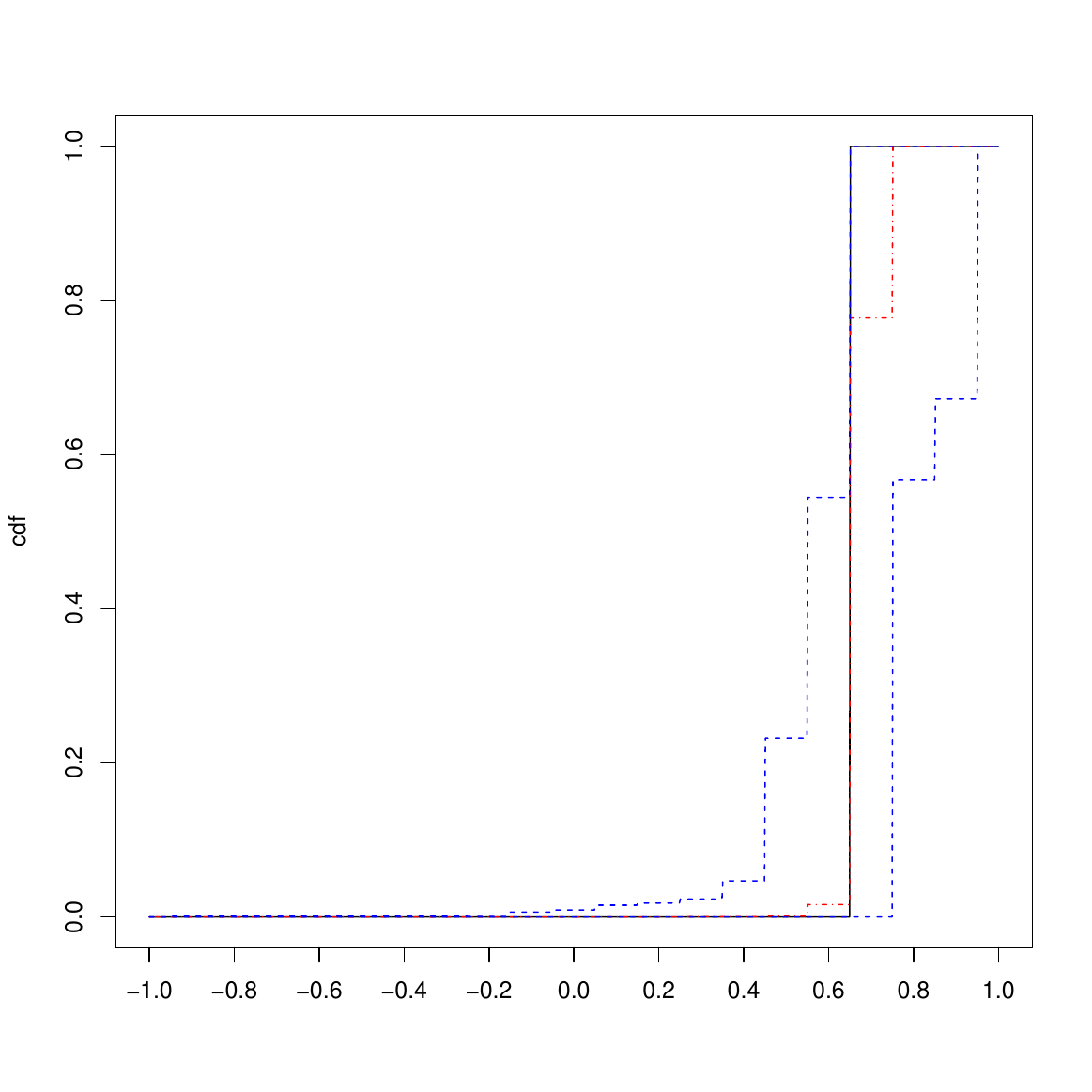}
& \\
\end{tabular}
}
\vspace{0pt}
         \caption{Plot of Median and $90\%$ confidence band for spectral cdf of \textbf{MA coefficient matrix} corresponding to the \textbf{case 2.1}. \textbf{Dash-Dot Red} curve: median, \textbf{Dashed Blue} curve: $90\%$ confidence band, \textbf{Black Solid} curve: true spectral cdf}
\end{center}
\end{figure}

\begin{figure}[H] \label{fig:AR2_A1}
\begin{center}
\resizebox{\textwidth}{!}{%
\begin{tabular}{cc}
$(p,n)=(400,1600)$ & $(p,n)=(400,800)$ \\ \vspace{-20pt}
\includegraphics[width=2.9in, height=3in, bb = 50 50 650 650]{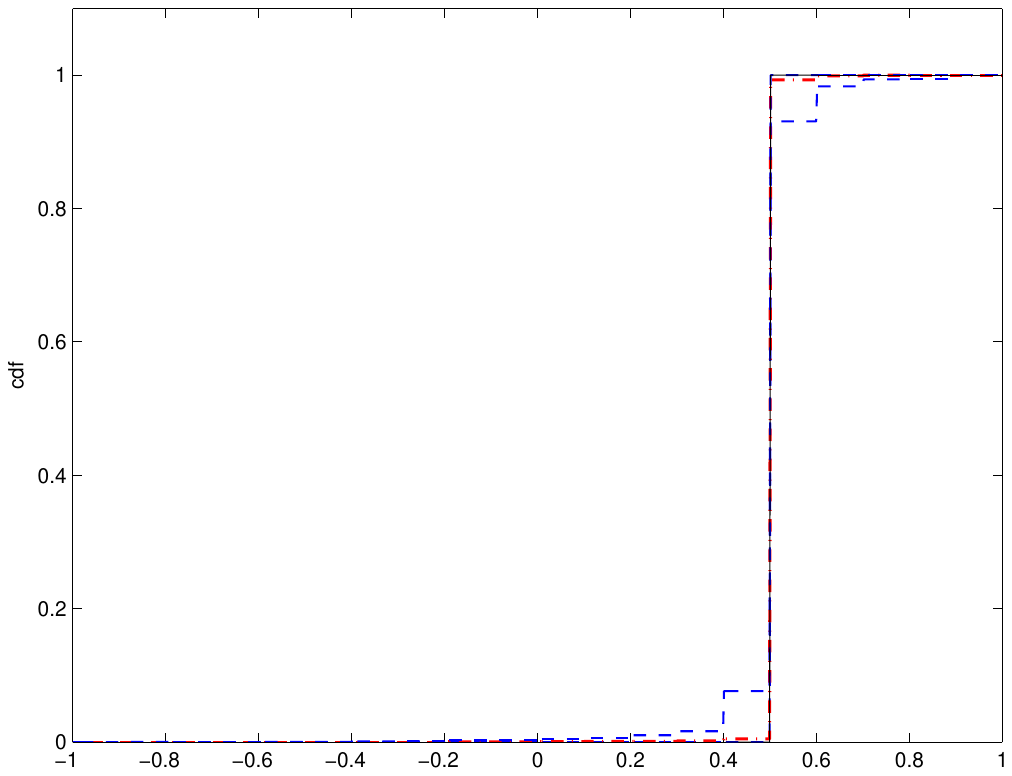}
& \includegraphics[width=2.9in, height=3in, bb = 50 50 650 650]{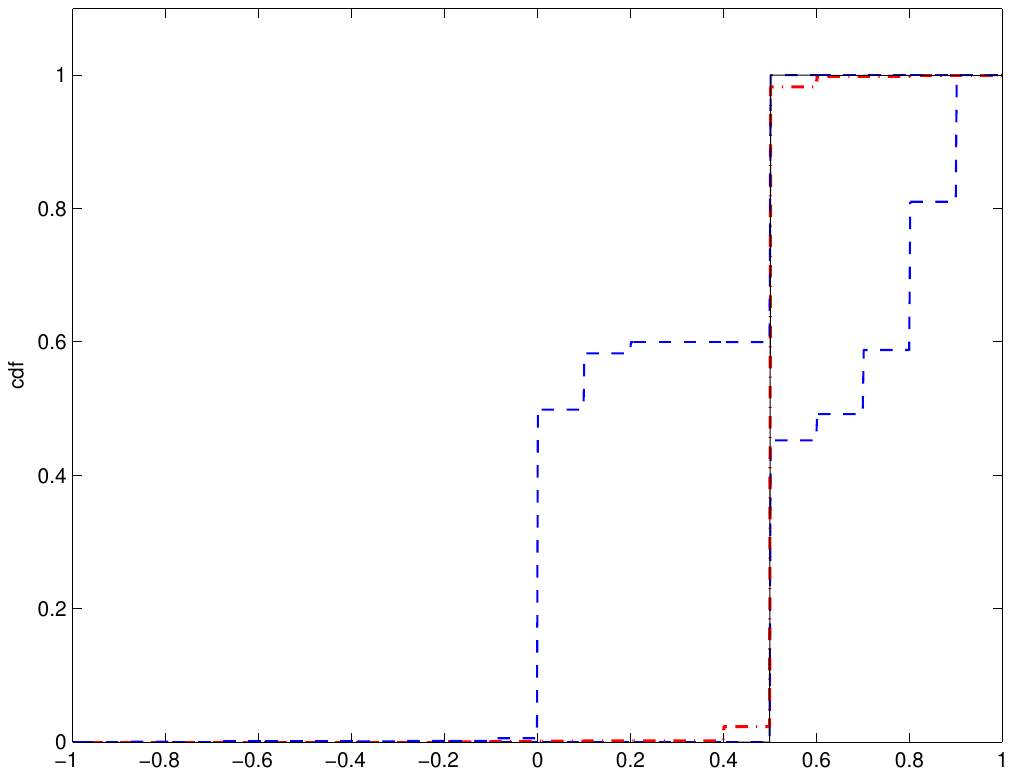}\vspace{-20pt}\\ \vspace{-10pt}
$(p,n)=(200,800)$ & $(p,n)=(200,400)$ \vspace{-5pt}\\ \vspace{-10pt}
\includegraphics[width=2.9in, height=3in, bb = 50 50 650 650]{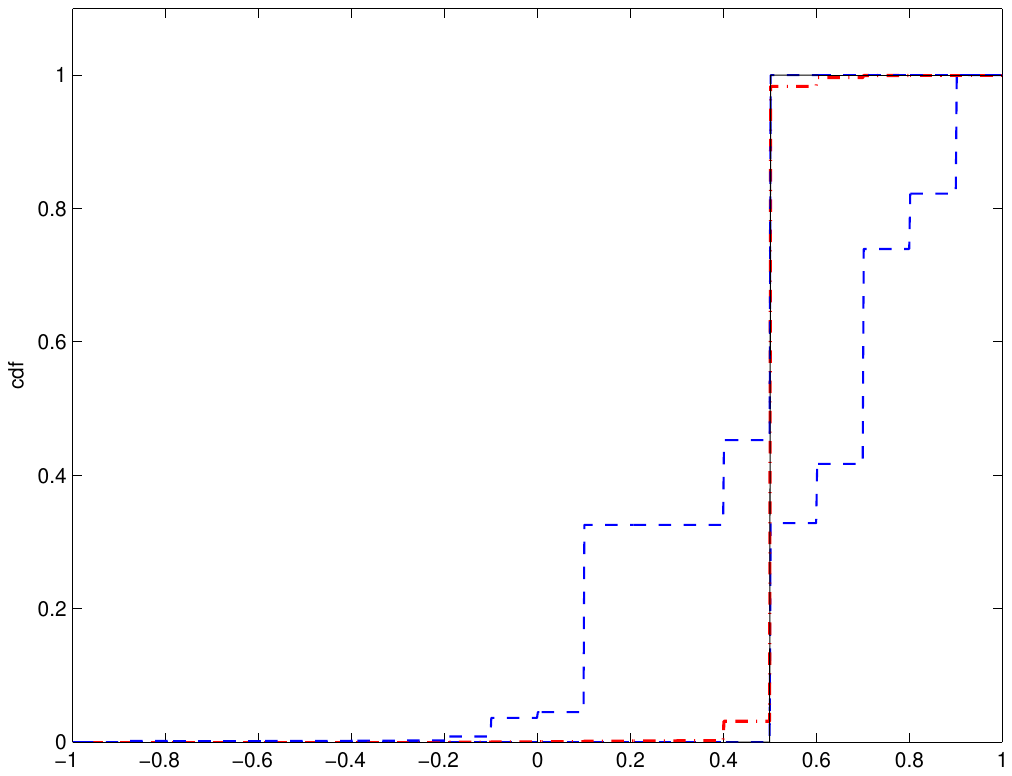}
& \includegraphics[width=2.9in, height=3in, bb = 50 50 650 650]{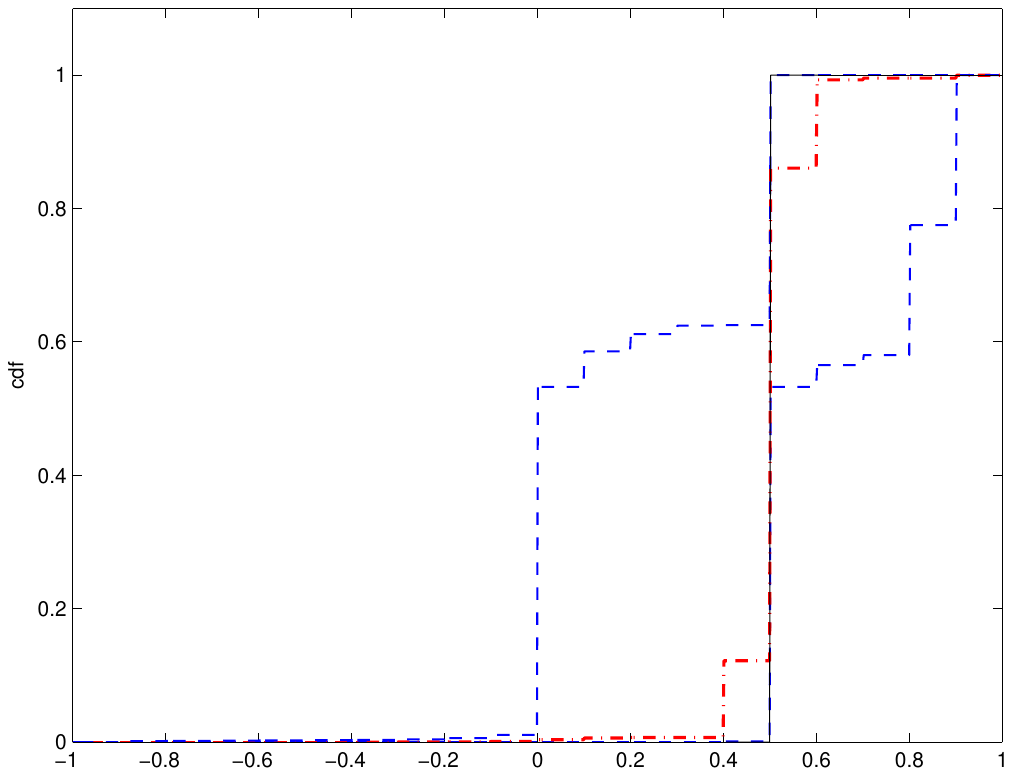} \vspace{-20pt}\\ \vspace{-20pt}
$(p,n)=(200,100)$ & \vspace{-10pt}\\ \vspace{-30pt}
\hspace{10pt}\includegraphics[width=2.8in, height=2.8in, bb = 50 50 650 650]{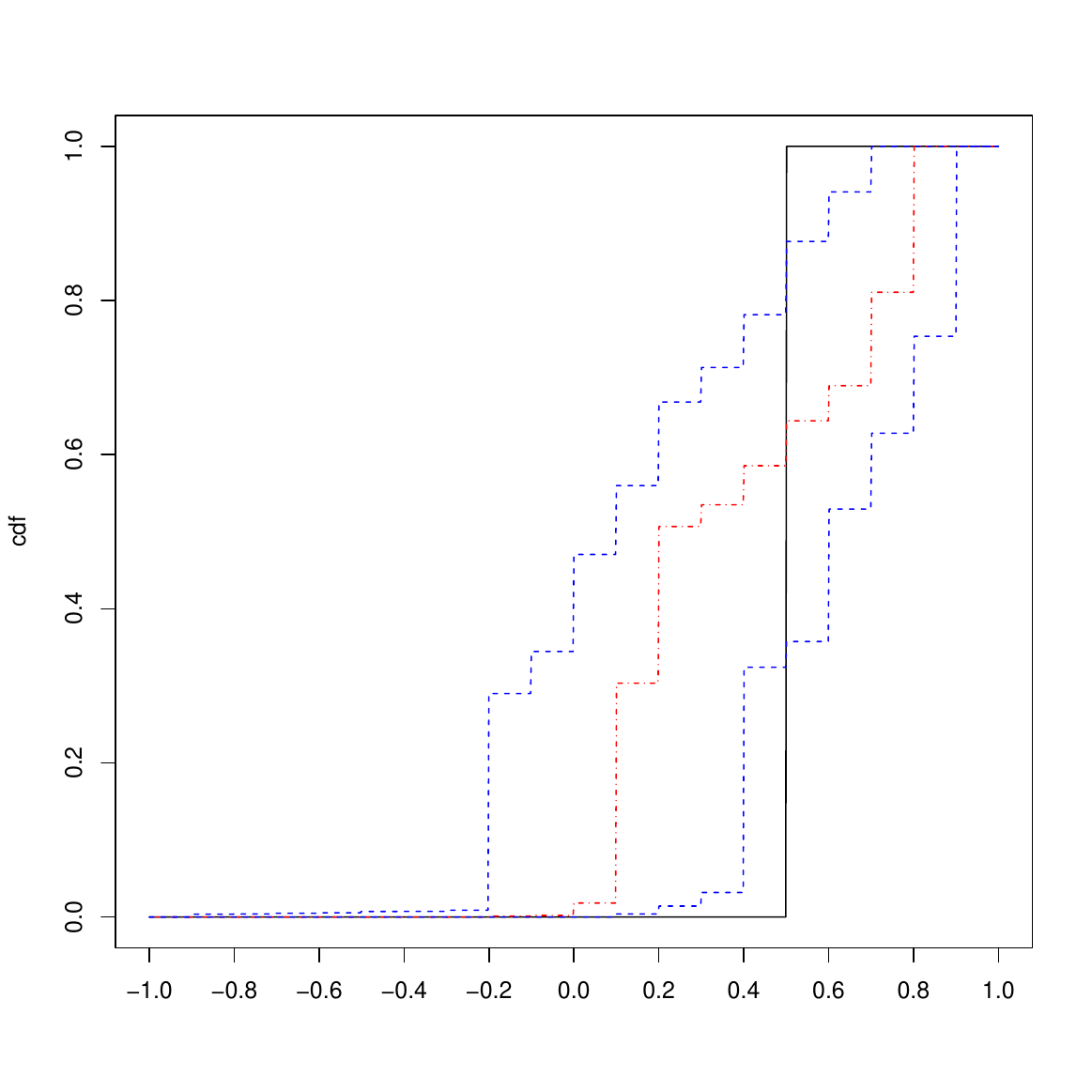}
& \\
\end{tabular}
}
\vspace{25pt}
         \caption{Plot of Median and $90\%$ confidence band for spectral cdf of $\boldsymbol{A_1}$ corresponding to the \textbf{case 2.2}. \textbf{Dash-Dot Red} curve: median, \textbf{Dashed Blue} curve: $90\%$ confidence band, \textbf{Black Solid} curve: true spectral cdf}
\end{center}
\end{figure}

\begin{figure}[H] \label{fig:AR2_Sigma}
\begin{center}
\resizebox{\textwidth}{!}{%
\begin{tabular}{cc}
$(p,n)=(400,1600)$ & $(p,n)=(400,800)$ \\ \vspace{-20pt}
\includegraphics[width=2.9in, height=3in, bb = 50 50 650 650]{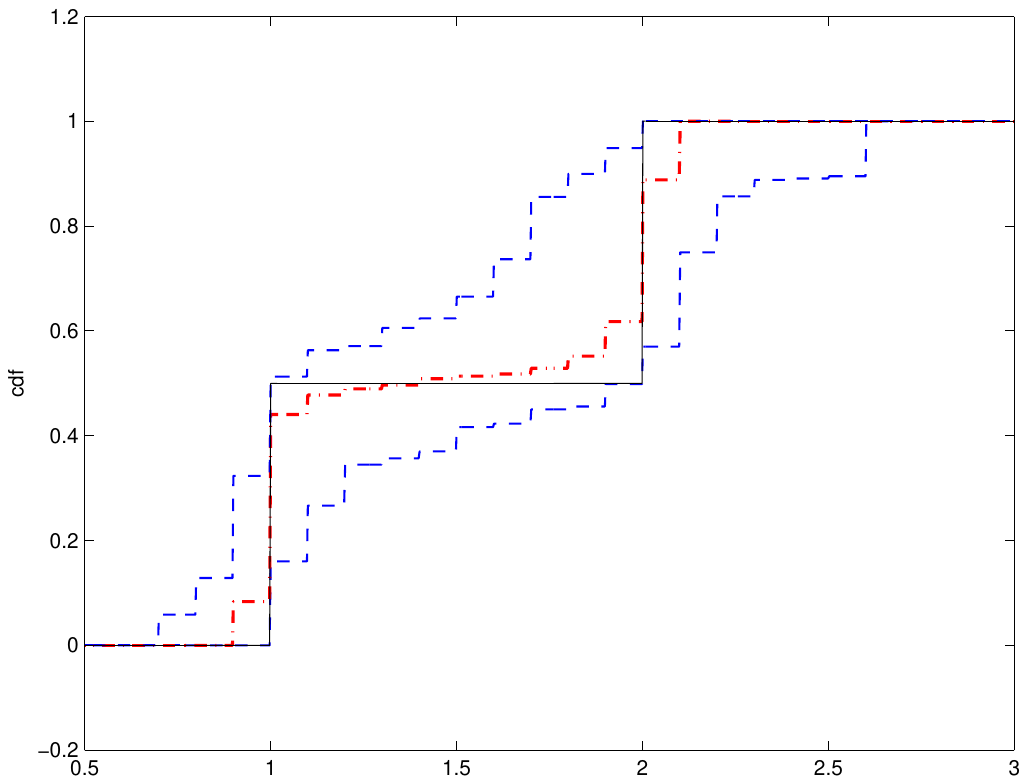}
& \includegraphics[width=2.9in, height=3in, bb = 50 50 650 650]{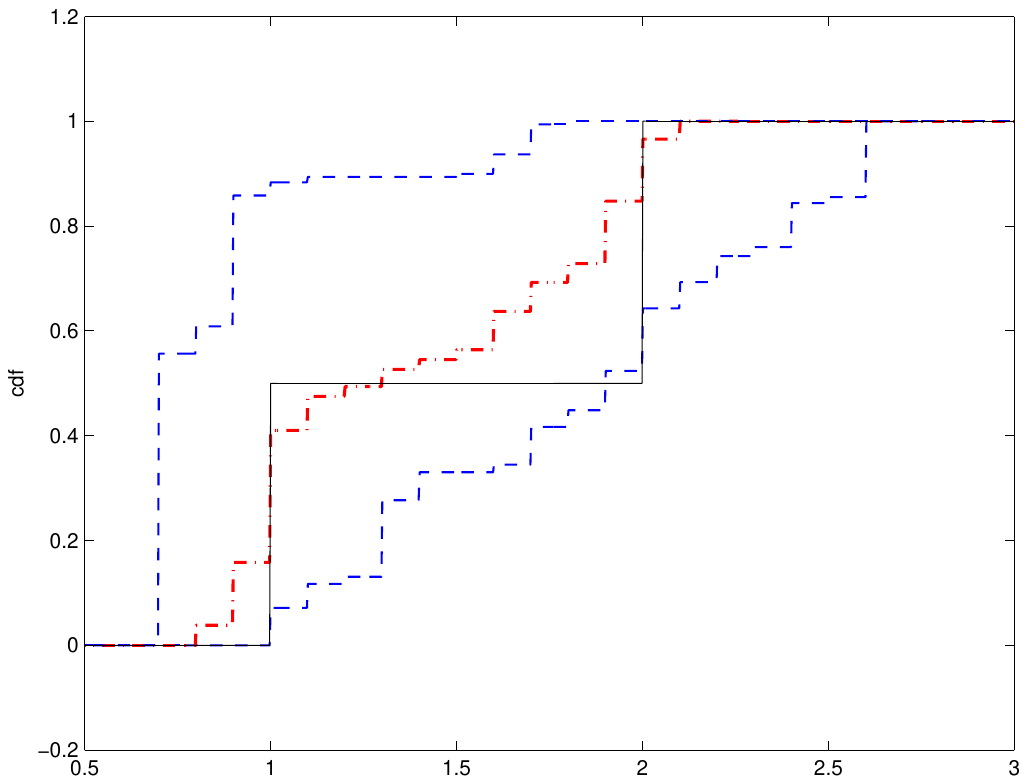}\vspace{-20pt}\\ \vspace{-10pt}
$(p,n)=(200,800)$ & $(p,n)=(200,400)$ \vspace{-5pt}\\ \vspace{-10pt}
\includegraphics[width=2.9in, height=3in, bb = 50 50 650 650]{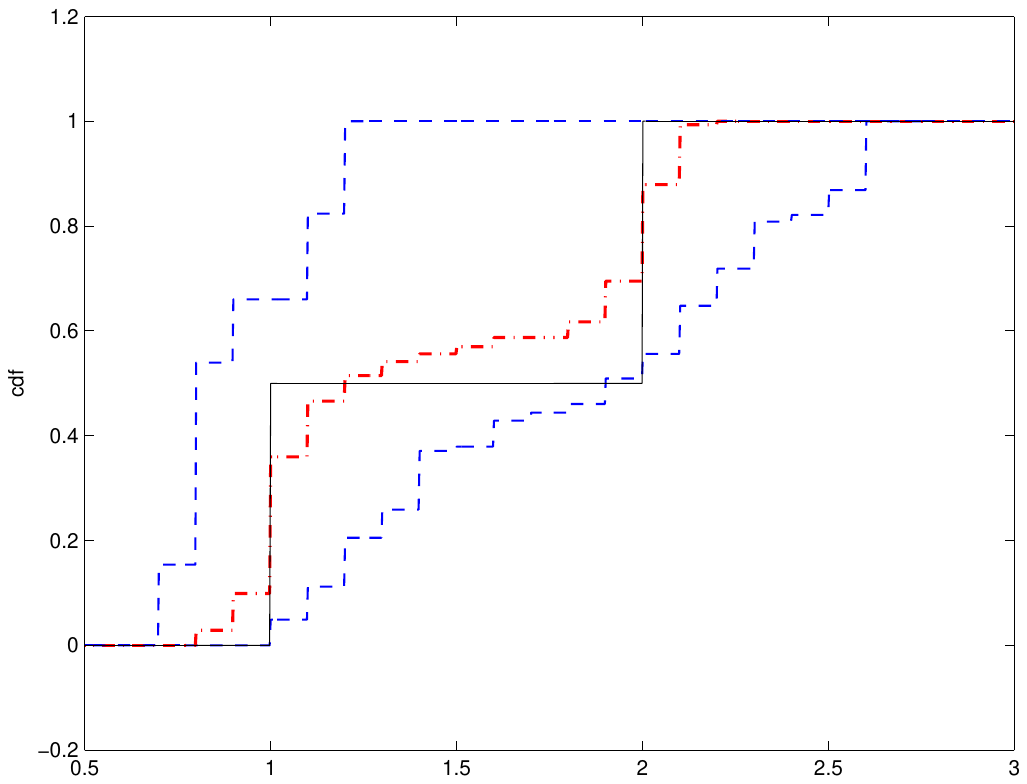}
& \includegraphics[width=2.9in, height=3in, bb = 50 50 650 650]{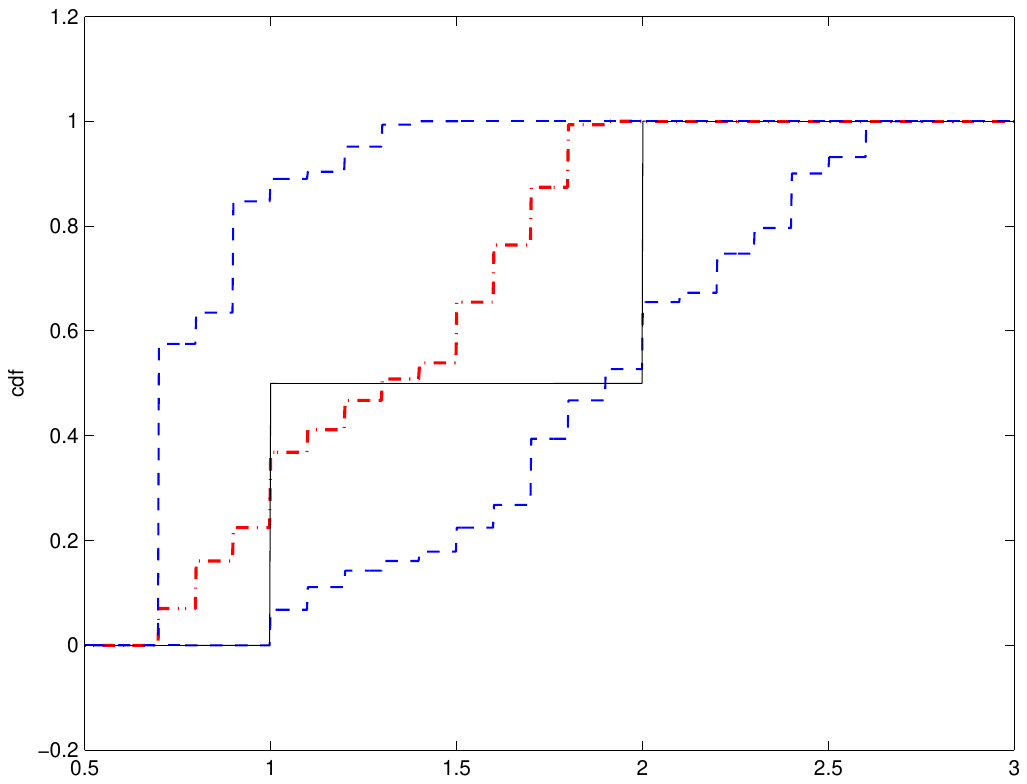}\vspace{-20pt}\\ \vspace{-20pt}
$(p,n)=(200,100)$ & \vspace{-10pt}\\ \vspace{-30pt}
\hspace{10pt}\includegraphics[width=2.8in, height=2.8in, bb = 50 50 650 650]{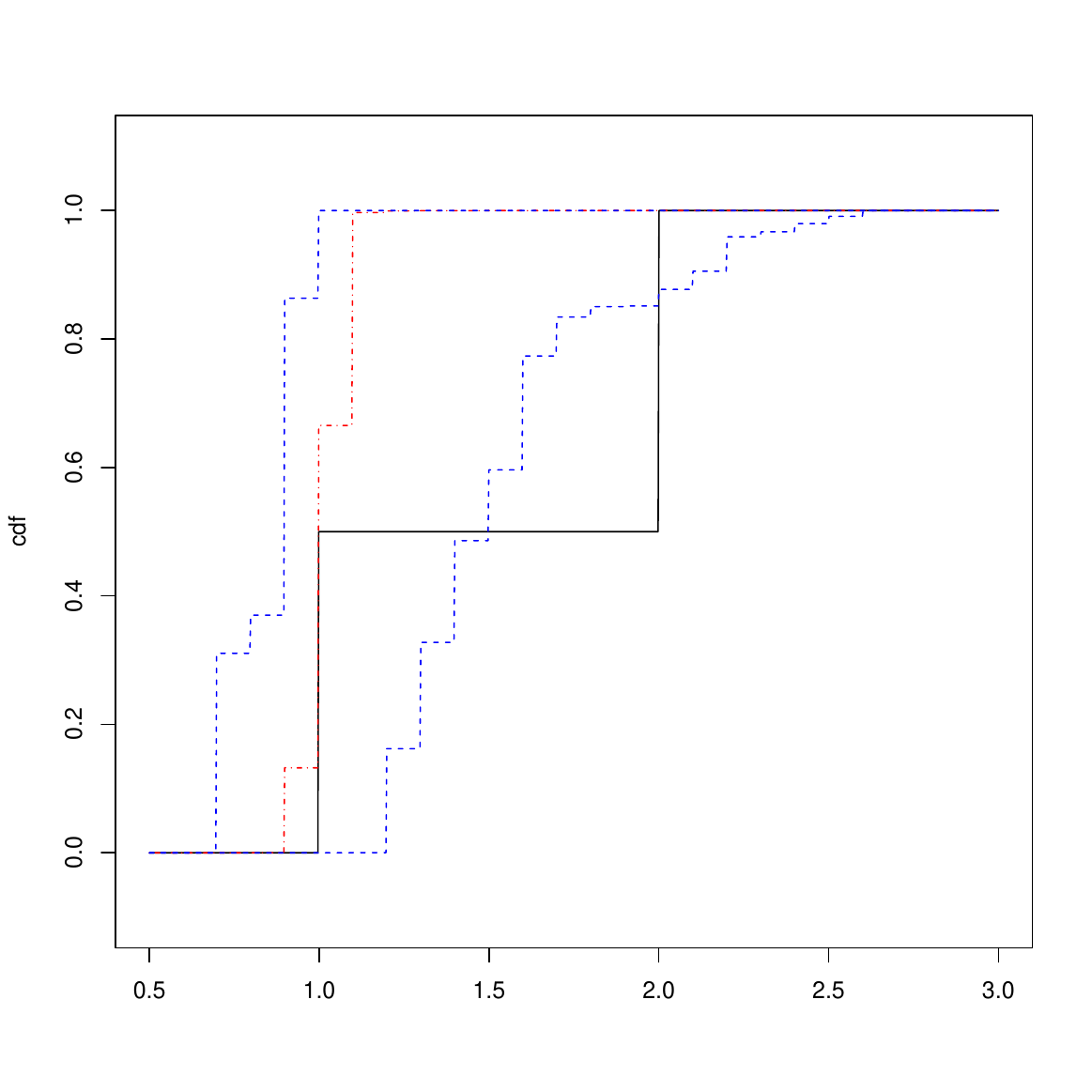}
& \\
\end{tabular}
}
\vspace{25pt}
         \caption{Plot of Median and $90\%$ confidence band for spectral cdf of $\boldsymbol{\Sigma}$ corresponding to the \textbf{case 2.2}. \textbf{Dash-Dot Red} curve: median, \textbf{Dashed Blue} curve: $90\%$ confidence band, \textbf{Black Solid} curve: true spectral cdf}
\end{center}
\end{figure}

\noindent {\it ARMA$(1,1)$ and AR$(2)$ models.} In the second stage, we fix $\mathcal G$ to contain eight B-spline weight functions and consider Case 1.5 together with the following two models: \begin{align*} 
\text{Case 2.1:}\qquad & X_t-\mbf A_1X_{t-1} = Z_t+\mbf M_1Z_{t-1}, \\[-2mm] & F^{A_1}=\delta_{-0.35}, \qquad F^{M_1}=\delta_{0.65}, \qquad F^\Sigma=0.5\,\delta_1+0.5\,\delta_2; \\[2mm] 
\text{Case 2.2:}\qquad & X_t = \mbf A_1X_{t-1}
+ \mbf A_2X_{t-2} + Z_t, \\[-2mm] & F^{A_1}=\delta_{0.5}, \qquad F^{A_2}=\delta_{-0.8}, \qquad F^\Sigma=0.5\,\delta_1+0.5\,\delta_2. 
\end{align*}
We consider 
\[ 
(p,n) = (400,1600),\ (400,800),\ (200,800),\ (200,400),\ (200,100). 
\] 
For computational tractability, we impose the product structure for the joint spectral distribution described in Section~\ref{sub:SpecEst_Algm}. Thus, the eigenvalue coordinates corresponding to the coefficient matrices and $\Sigma$ are treated as independent in the candidate joint distribution. For each configuration, 100 independent samples are generated. Figures~\ref{fig:ARMA_AR}.1--\ref{fig:ARMA_MA}.3 report the median and pointwise $90\%$ simulation bands for the estimated marginal spectral distribution functions in Case 2.1, while Figures~\ref{fig:AR2_A1}.4--\ref{fig:AR2_Sigma}.6 give the corresponding results for Case 2.2. As expected, estimation becomes more precise as the amount of information relative to dimension increases; in particular, the simulation bands generally become narrower as $p/n$ decreases.

\subsection{Model selection}
\label{subsec:model_selection_sim}

We next investigate the ability of the bootstrap model-selection procedure in Section~\ref{sub:SpecEst_ModelSeln} to identify the order of an autoregressive process. Data are generated from the AR$(2)$ model
\[
X_t=\mathbf{A}_1X_{t-1}+\mathbf{A}_2X_{t-2}+Z_t,
\]
where $\mathbf{A}_1=0.5\,\mathbf{I}_p$, $\mathbf{A}_2=-0.8\,\mathbf{I}_p$, and $Z_1,Z_2,\ldots$ are independent Gaussian random vectors with covariance matrix $\Sigma$ having spectral distribution $F^\Sigma=0.5\,\delta_1+0.5\,\delta_2$.

We generate $R=100$ independent data sets from this model. For each data set, we fit both an AR$(1)$ and an AR$(2)$ model using the proposed
estimation procedure. From each fitted model, we then generate $B=100$ bootstrap samples and compute the discrepancy $L_b$, $b=1,\ldots,B$,
defined in Section~\ref{sub:SpecEst_ModelSeln}. Model fit is assessed using symmetrized sample autocovariance matrices at lags $\tau=0,\ldots,5$. For each lag and each simulation replication, the
candidate model with the smaller mean bootstrap discrepancy
\[
\bar L = \frac{1}{B}\sum_{b=1}^B L_b
\]
is preferred.

The results are reported in Table~B.5 of the Supplementary Material. For each of the six lags considered, the AR$(2)$ model is selected in all $R=100$ simulation replications. Thus, in this setting, the proposed bootstrap criterion consistently distinguishes the correctly specified
AR$(2)$ model from the underfitted AR$(1)$ model.

\section{Data analysis}
\label{sec:Data_Analysis} 

We illustrate the proposed methodology using daily closing prices of $p=486$ companies included in the S\&P 500, observed from May 8, 2012 to October 17, 2016. Prices were adjusted for stock splits, and the analysis is based on the corresponding log returns 
\[ 
R_{i,t} = \log X_{i,t}-\log X_{i,t-1}, 
\qquad i=1,\ldots,p,\quad t=1,\ldots,n, 
\] 
where $n=1111$. Let \[ \mbf R = [\bs R_1,\ldots,\bs R_n] \] denote the resulting $p\times n$ return matrix. Our analysis proceeds in two stages. We first remove the dominant common factor from the return series and examine the remaining cross-sectional dependence. We then apply the proposed methodology to determine whether the residual series exhibits temporal dependence and, if so, which of several candidate linear-process models provides the best fit.

\medskip

\begin{figure}[H] \label{fig:scree}
\begin{center}
\includegraphics[width=3.0in, height=3in, bb = 50 50 600 600]{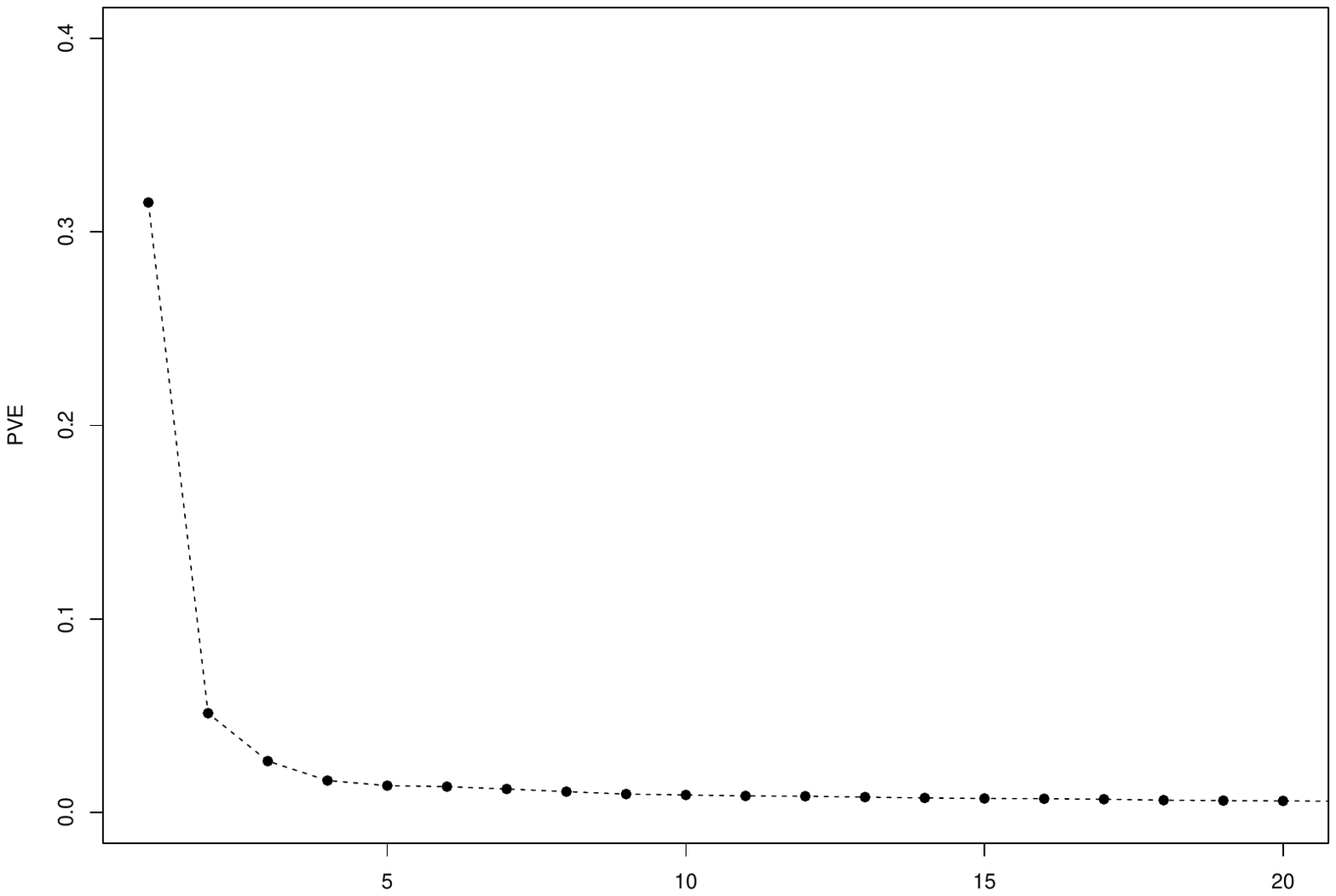} \\

         \caption{Plot of the proportion of variation explained (PVE) by leading factors of the log return series.}
\end{center}
\end{figure}

\begin{figure} \label{fig:pairwise_corr_density}
\begin{center}
\includegraphics[width=3.0in, height=3in, bb = 50 50 600 600]{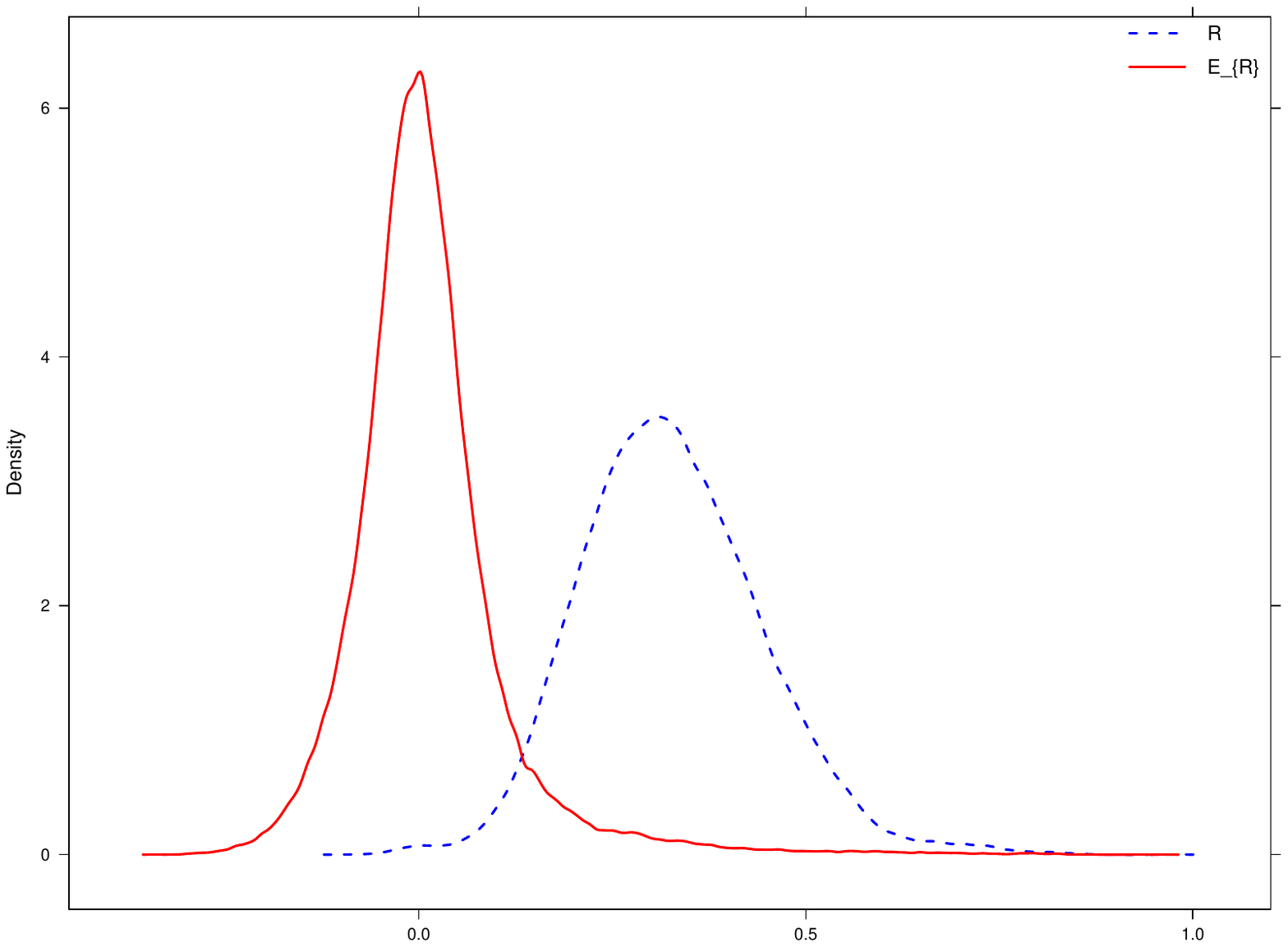}\\

         \caption{Plot of the densities of the pairwise correlation between coordinates of the series $\mathbf{R}$, in blue, and $\mathbf{E}_R^{(1)}$, in red.}
\end{center}
\end{figure}

\noindent{\it Removal of the dominant factor.} Consider the factor representation 
\[ 
\bs R_t-E(\bs R_t) = \mbf L\bs F_t+\bs E_t, 
\qquad t=1,\ldots,n, 
\] 
where $\bs F_t$ is a vector of common factors, $\mbf L$ is the matrix of factor loadings, and $\bs E_t$ is the idiosyncratic component. To identify the dominant empirical factors, let $\bar{\bs R}$ denote the sample mean and consider the singular value decomposition 
\[ 
\mbf R-\bar{\bs R}\mbf 1_n^T 
= \mbf U\bs\Lambda\mbf V^T, 
\] 
where $\mbf U=[\bs u_1|\cdots|\bs u_p]$ and $\mbf V=[\bs v_1|\cdots|\bs v_n]$. We define the residual matrix after removing the first $s$ principal components by 
\[ 
\mbf E_R^{(s)} 
= \mbf R-\bar{\bs R}\mbf 1_n^T 
- \sum_{j=1}^s \sqrt{\lambda_j}\, \bs u_j\bs v_j^T. 
\]
Figure~\ref{fig:scree} shows the proportion of variation explained by the leading principal components. The first component alone accounts for approximately $31.5\%$ of the total variation in the return series, indicating a strong common market component. Figure~ \ref{fig:pairwise_corr_density} compares the distributions of pairwise correlations before and after removal of this component. The correlations of the original returns are predominantly positive, whereas those of the residual series $\mbf E_R^{(1)}$ are concentrated much more closely around zero. The removal of the dominant factor substantially reduces contemporaneous cross-sectional dependence, but it does not imply that the residual series is temporally independent. We therefore use the proposed methodology to investigate the dynamic structure remaining in $\mbf E_R^{(1)}$.

\medskip

\noindent {\it Dynamic modeling of the residual series.} Let $\bs E_t^{(1)}$ denote the $t$th column of $\mbf E_R^{(1)}$. We compare the following models for $(\bs E_t^{(1)})$: 
\begin{itemize} 
\item an iid model, denoted by $\mathrm{Ind}$; 
\item an AR$(1)$ model; 
\item an AR$(2)$ model with product-form joint spectral distribution $F^{A,\Sigma} = F^{A_1}F^{A_2}F^\Sigma$; 
\item an ARMA$(1,1)$ model with product-form joint spectral distribution $F^{A,\Sigma} = F^{AR}F^{MA}F^\Sigma$, denoted by $\mathrm{ARMA}_{\mathrm{Ind}}(1,1)$. 
\end{itemize} 
The product representations impose independence among the eigenvalue coordinates of the corresponding coefficient matrices and innovation covariance, as discussed in Section~\ref{sub:SpecEst_Algm}. The iid model provides a benchmark for assessing whether the residual series contains temporal dependence that cannot be explained by contemporaneous covariance alone. Its covariance spectrum is estimated from $\mbf E_R^{(1)}$ by applying the proposed procedure with the autoregressive coefficient set equal to zero.

We compare the fitted models using the bootstrap procedure of Section~\ref{sub:SpecEst_ModelSeln}. The diagnostic statistics are the symmetrized sample autocovariance matrices at lags $\tau=0,\ldots,5$. For each fitted model, $B=500$ bootstrap samples are generated, and the discrepancy between the ordered eigenvalues of each bootstrap autocovariance matrix and those of the corresponding matrix computed from the observed residual series is recorded. Table~\ref{table:model_Selection} reports the percentage of bootstrap samples yielding each observed ordering of the four model discrepancies. Here $\mathcal M_1\prec\mathcal M_2$ means that $\mathcal M_1$ produces a smaller bootstrap discrepancy than $\mathcal M_2$. For example, at lag zero the ordering \[ \mathrm{ARMA}_{\mathrm{Ind}}(1,1) \prec \mathrm{AR}(1) \prec \mathrm{Ind} \prec \mathrm{AR}(2) \] occurs in $83\%$ of the bootstrap samples.

The rankings vary with the lag used in the comparison, indicating that different models reproduce different aspects of the dependence structure. At lag zero, the ARMA$(1,1)$ model performs best in every bootstrap sample, whereas at lags one and two the iid model most often has the smallest discrepancy. At lags three through five, the ARMA$(1,1)$ model again performs particularly well. Taken together, these results suggest that the residual series exhibits dependence features that are not uniformly captured by a single low-order autoregressive model. They also illustrate the importance of examining several lags when assessing the fit of competing dynamic models.

\begin{table}[ht]
\centering
\caption{Bootstrap frequencies (in percent) of the observed rankings of
the candidate models, based on the eigenvalue discrepancy for
symmetrized sample autocovariance matrices at lags $\tau=0,\ldots,5$.
The notation $\mathcal M_1\prec\mathcal M_2$ indicates that
$\mathcal M_1$ has smaller bootstrap discrepancy than $\mathcal M_2$.}
\label{table:model_Selection}
\begin{tabular}{lrrrrrr}
\hline
& \multicolumn{6}{c}{Lag $\tau$} \\
\cline{2-7}
Model ranking & 0 & 1 & 2 & 3 & 4 & 5 \\
\hline
$\mathrm{ARMA}_{\mathrm{Ind}}(1,1)
 \prec \mathrm{AR}(1)
 \prec \mathrm{Ind}
 \prec \mathrm{AR}(2)$
& 83.0 & 0 & 0 & 94.0 & 5.0 & 0 \\

$\mathrm{ARMA}_{\mathrm{Ind}}(1,1)
 \prec \mathrm{Ind}
 \prec \mathrm{AR}(1)
 \prec \mathrm{AR}(2)$
& 17.0 & 0 & 0 & 5.8 & 5.0 & 10.6 \\

$\mathrm{Ind}
 \prec \mathrm{ARMA}_{\mathrm{Ind}}(1,1)
 \prec \mathrm{AR}(1)
 \prec \mathrm{AR}(2)$
& 0 & 47.2 & 4.8 & 0 & 0 & 0.2 \\

$\mathrm{Ind}
 \prec \mathrm{ARMA}_{\mathrm{Ind}}(1,1)
 \prec \mathrm{AR}(2)
 \prec \mathrm{AR}(1)$
& 0 & 52.8 & 95.2 & 0.2 & 0.2 & 0.6 \\

$\mathrm{ARMA}_{\mathrm{Ind}}(1,1)
 \prec \mathrm{Ind}
 \prec \mathrm{AR}(2)
 \prec \mathrm{AR}(1)$
& 0 & 0 & 0 & 0 & 94.8 & 88.6 \\
\hline
\end{tabular}
\end{table}

\section{Conclusion}
\label{sec:Conclusion}

We have developed a random-matrix-based framework for estimating the joint spectral distribution of the coefficient matrices and innovation covariance of a class of high-dimensional stationary linear processes. Under the assumption of simultaneous diagonalizability, we derived the limiting spectral distributions of weighted sample periodograms and used their Stieltjes transforms to construct a minimum-distance estimator of the underlying joint spectrum. We established consistency of the estimator when the support of the joint spectral distribution is known, and also proposed an estimator of the spectral density matrix and a bootstrap procedure for comparing candidate linear-process models. Simulation studies demonstrate the finite-sample performance of the proposed methods, while the analysis of S\&P 500 returns illustrates their use for investigating dependence remaining after removal of a dominant common factor.

Several extensions of the proposed framework are of interest. First, simultaneous diagonalizability is an important structural restriction. A natural direction for future work is to determine whether the limiting spectral results and the resulting estimation procedure remain valid under weaker conditions, such as approximate simultaneous diagonalizability or other forms of structured noncommutativity. Second, the consistency theory developed here assumes that the joint spectral distribution is supported on a known finite grid and is therefore essentially parametric. Although the estimation procedure can be applied approximately to more general distributions by using a sufficiently fine grid, establishing consistency when the support is unknown or the population distribution is continuous would provide a genuinely nonparametric extension of the theory Finally, estimation of the joint spectral distribution does not by itself recover the coefficient matrices. More direct estimation of these matrices, and of the spectral density matrix, will require a better understanding of the common eigenbasis and, in particular, further theoretical results on the eigenvectors of high-dimensional sample periodograms.

\begin{appendix}

\section{Proofs}
\label{sec:appendix}

\subsection{Outline of the proof of Theorem~\ref{theorem1}}
\label{subsec:Proof_Theorem1}

We provide an outline of the proof of Theorem~\ref{theorem1}; detailed technical arguments can be found in the Ph.D.\ thesis of the first author \cite{NamdariThesis2018}. The main step is the derivation of approximating equations for the Stieltjes transform of $F_g^{(n)}$. To illustrate the argument without unnecessary technical complications, we first consider an MA$(1)$ process with Gaussian innovations. We then describe the main steps needed to establish existence and uniqueness of the limiting equations and to extend the result to linear processes with non-Gaussian innovations.

Suppose, for simplicity, that $\Sigma=\mbf I_p$ and that $(X_t\colon t\in\mbb Z)$ is a complex-valued Gaussian MA$(1)$ process satisfying the assumptions of Theorem~\ref{theorem1}. The key idea is to replace the lag operator by a circulant approximation. This allows the time-domain representation of the process to be transformed into the Fourier domain, where the resulting matrix has a particularly simple structure.

Let $\mbf X_n=[X_1:\cdots:X_n]$ and $\mbf Z_n=[Z_1:\cdots:Z_n]$. For the MA$(1)$ process, $\mbf X_n=\mbf Z_n+\mbf A_1\mbf Z_n\mbf L$, where $\mbf L=[0:e_1:\cdots:e_{n-1}]$ is the lag operator. We replace $\mbf L$ by the circulant matrix $\tilde{\mbf L}=[e_n:e_1:\cdots:e_{n-1}]$, and define $\tilde{\tilde{\mbf X}}_n=\mbf Z_n+\mbf A_1\mbf Z_n\tilde{\mbf L}$. Since $\tilde{\mbf L}$ is circulant, it is diagonalized by the discrete Fourier basis: $\tilde{\mbf L}=\mbf U_{\tilde L}\bs\Lambda_{\tilde L}\mbf U_{\tilde L}^*$, where
\[
\bs\Lambda_{\tilde L}
=
\operatorname{diag}(\eta_1,\ldots,\eta_n),
\qquad
\eta_t=e^{i\theta_t},
\qquad
\theta_t=\frac{2\pi t}{n},
\]
and
\[
\mbf U_{\tilde L}
=
[u_1:\cdots:u_n],
\qquad
u_t
=
\frac{1}{\sqrt n}
\left(
\eta_t,\eta_t^2,\ldots,\eta_t^n
\right)^T.
\]
Rotating the circulant approximation into the Fourier basis, define $\breve{\mbf X}_n=\tilde{\tilde{\mbf X}}_n\mbf U_{\tilde L}$. This matrix differs from the discrete Fourier transform $\tilde{\mbf X}_n$ of the original observations only through the
circulant approximation of the lag operator. The Gaussian assumption, together with the unitary Fourier transformation, gives a representation
of $\breve{\mbf X}_n$ in terms of independent Gaussian vectors that is particularly convenient for the resolvent analysis below.

Define
\[
\breve{\mbf S}
=
\frac{1}{n}
\mbf W_g
\breve{\mbf X}_n^*
\breve{\mbf X}_n
\mbf W_g
=
\frac{1}{n}\mbf Y\mbf Y^*
=
\frac{1}{n}
\sum_{j=1}^p
Y_jY_j^*,
\]
where $\mbf Y=[Y_1:\cdots:Y_p]$, $Y_j=\mbf W_g\bs\Psi^*(\bs\lambda_j)\breve Z_j$, with $\bs\Psi(\bs\lambda_j)= \operatorname{diag}
\left\{\psi(\bs\lambda_j,\theta_1),
\ldots, \psi(\bs\lambda_j,\theta_n)
\right\}$, and $\breve Z_1,\ldots,\breve Z_p$ are independent Gaussian random vectors.

For $z\in\mbb C^+$, let $\breve{\mbf R}(z)=(\breve{\mbf S}-z\mbf I)^{-1}$
denote the resolvent of $\breve{\mbf S}$. For
$j=1,\ldots,p$, define the rank-one-deleted matrix and its resolvent by
\[
\breve{\mbf S}_j
=
\breve{\mbf S}
-
\frac{1}{n}Y_jY_j^*,
\qquad
\breve{\mbf R}_j(z)
=
(\breve{\mbf S}_j-z\mbf I)^{-1}.
\]
The circulant approximation is a finite-rank perturbation of the original lag representation. Consequently, the empirical spectral distributions of $\tilde{\mbf S}_g^{(n)}$ and $\breve{\mbf S}$ have the same limiting behavior. It therefore suffices to establish the almost sure convergence of the empirical spectral distribution of $\breve{\mbf S}$. Equivalently, we study the
pointwise almost sure convergence of its Stieltjes transform
\[
\breve{\mathcal S}_n(z)
=
\frac{1}{n}
\operatorname{tr}\{\breve{\mbf R}(z)\},
\qquad z\in\mbb C^+.
\]
The remainder of the proof derives deterministic approximating equations for $\breve{\mathcal S}_n(z)$ and the associated Stieltjes kernel, and  then establishes convergence to the unique solution of the limiting system in Theorem~\ref{theorem1}. For notational simplicity, we suppress the subscript $g$ on $\mbf W_g$ in the calculations below.

\subsubsection{Derivation of the approximating equations} 

We next derive finite-sample approximating equations whose limiting
counterparts are the equations characterizing $\mathcal{S}_g$ and $K_g$
in Theorem~\ref{theorem1}. The argument is based on constructing a
deterministic equivalent for the resolvent
\[
\breve{\mbf R}(z)
=
(\breve{\mbf S}-z\mbf I)^{-1}.
\]
More precisely, we seek a matrix $\mbf H(z)$ such that, for every
sequence of $n\times n$ Hermitian matrices $\mbf C_n$ with uniformly
bounded operator norm,
\begin{equation}\label{eq:det_equiv}
\frac{1}{n}
\operatorname{tr}
\left[
\left\{
(\mbf I+\mbf H(z))^{-1}
+
z\breve{\mbf R}(z)
\right\}
\mbf C_n
\right]
\longrightarrow 0.
\end{equation}
For the heuristic derivation below, we use the notation $\approx$ to
indicate that the difference between two quantities is asymptotically
negligible. Observe that
\begin{align*}
(\mbf I+\mbf H(z))^{-1}
+
z\breve{\mbf R}(z)
&=
\breve{\mbf R}(z)
[\breve{\mbf R}(z)]^{-1}
(\mbf I+\mbf H(z))^{-1}
\\
&\quad+
z\breve{\mbf R}(z)
(\mbf I+\mbf H(z))
(\mbf I+\mbf H(z))^{-1}
\\
&=
\breve{\mbf R}(z)
\{\breve{\mbf S}+z\mbf H(z)\}
(\mbf I+\mbf H(z))^{-1}.
\end{align*}
Thus, it is enough to choose $\mbf H(z)$ so that
\[
\breve{\mbf R}(z)\breve{\mbf S}
+
z\breve{\mbf R}(z)\mbf H(z)
\approx 0.
\]
Recall that
\[
\breve{\mbf S}
=
\frac{1}{n}\sum_{k=1}^pY_kY_k^*,
\qquad
\breve{\mbf S}_k
=
\breve{\mbf S}-\frac{1}{n}Y_kY_k^*,
\]
and let
\[
\breve{\mbf R}_k(z)
=
(\breve{\mbf S}_k-z\mbf I)^{-1}.
\]
Using the Sherman--Morrison formula,
\begin{align*}
\breve{\mbf R}(z)\breve{\mbf S}
&=
\frac{1}{n}
\sum_{k=1}^p
\left(
\breve{\mbf S}_k-z\mbf I
+\frac{1}{n}Y_kY_k^*
\right)^{-1}
Y_kY_k^*
\\
&=
\frac{1}{n}
\sum_{k=1}^p
\left\{
\breve{\mbf R}_k(z)
-
\frac{
n^{-1}\breve{\mbf R}_k(z)Y_kY_k^*
\breve{\mbf R}_k(z)
}{
1+n^{-1}Y_k^*\breve{\mbf R}_k(z)Y_k
}
\right\}
Y_kY_k^*
\\
&=
\frac{1}{n}
\sum_{k=1}^p
\frac{
\breve{\mbf R}_k(z)Y_kY_k^*
}{
1+n^{-1}Y_k^*\breve{\mbf R}_k(z)Y_k
}.
\end{align*}
Since $Y_k=\mbf W\bs\Psi(\bs\lambda_k)\breve Z_k$, this becomes
\[
\breve{\mbf R}(z)\breve{\mbf S}
=
\frac{1}{n}
\sum_{k=1}^p
\frac{
\breve{\mbf R}_k(z)
\mbf W\bs\Psi(\bs\lambda_k)
\breve Z_k\breve Z_k^*
\bs\Psi^*(\bs\lambda_k)\mbf W
}{
1+
n^{-1}\breve Z_k^*
\bs\Psi^*(\bs\lambda_k)\mbf W
\breve{\mbf R}_k(z)\mbf W
\bs\Psi(\bs\lambda_k)\breve Z_k
}.
\]
Standard concentration arguments for quadratic forms, together with the
rank-one relation between $\breve{\mbf R}_k(z)$ and
$\breve{\mbf R}(z)$, suggest the approximation
\begin{align*}
\breve{\mbf R}(z)\breve{\mbf S}
\approx
\breve{\mbf R}(z)
\frac{1}{n}
\sum_{k=1}^p
\frac{
\bs\Psi(\bs\lambda_k)\mbf W^2
\bs\Psi^*(\bs\lambda_k)
}{
1+
n^{-1}
\operatorname{tr}
\left[
\breve{\mbf R}(z)
\bs\Psi(\bs\lambda_k)\mbf W^2
\bs\Psi^*(\bs\lambda_k)
\right]
}.
\end{align*}
This motivates the choice
\begin{equation}\label{eq:H_z}
\mbf H(z)
=
-\frac{1}{zn}
\sum_{k=1}^p
\frac{
\bs\Psi(\bs\lambda_k)\mbf W^2
\bs\Psi^*(\bs\lambda_k)
}{
1+
n^{-1}
\operatorname{tr}
\left[
\breve{\mbf R}(z)
\bs\Psi(\bs\lambda_k)\mbf W^2
\bs\Psi^*(\bs\lambda_k)
\right]
}.
\end{equation}

\medskip
\noindent{\it Approximating equation for the Stieltjes transform.} Define
\[
\breve{\mathcal S}_n(z)
=
\frac{1}{n}\operatorname{tr}\{\breve{\mbf R}(z)\}
\]
and
\[
\breve K_n(z,\bs\lambda)
=
\frac{1}{n}
\operatorname{tr}
\left[
\breve{\mbf R}(z)
\bs\Psi(\bs\lambda)\mbf W^2
\bs\Psi^*(\bs\lambda)
\right].
\]
Taking $\mbf C_n=\mbf I$ in the deterministic-equivalent relation gives
\[
\breve{\mathcal S}_n(z)
\approx
-\frac{1}{zn}
\operatorname{tr}
\{(\mbf I+\mbf H(z))^{-1}\}.
\]
Since the matrices appearing in $\mbf H(z)$ are diagonal in the Fourier
basis, this yields
\begin{equation}\label{eq:St_hat}
\breve{\mathcal S}_n(z)
\approx
-\frac{1}{n}
\sum_{t=1}^n
\left[
z
-
c_n
\left\{
\frac{1}{p}
\sum_{k=1}^p
\frac{
g(\theta_t)h(\bs\lambda_k,\theta_t)
}{
1+\breve K_n(z,\bs\lambda_k)
}
\right\}
\right]^{-1},
\end{equation}
where $c_n={p}/{n}$. It is useful to introduce the finite-sample quantity
\[
\breve M_n(z,\theta_t)
=
\frac{1}{p}
\sum_{k=1}^p
\frac{
g(\theta_t)h(\bs\lambda_k,\theta_t)
}{
1+\breve K_n(z,\bs\lambda_k)
}.
\]
Then (\ref{eq:St_hat}) takes the more transparent form
\[
\breve{\mathcal S}_n(z)
\approx
\frac{1}{n}
\sum_{t=1}^n
\frac{1}{
c_n\breve M_n(z,\theta_t)-z
}.
\]
This is the finite-sample analogue of equation~(\ref{eq:s}) in
Theorem~\ref{theorem1}.

\medskip

\noindent{\it Approximating equation for the kernel.} Next, fix $\bs\lambda$ and take $\mbf C_n=\bs\Psi(\bs\lambda)\mbf W^2\bs\Psi^*(\bs\lambda)$ in (\ref{eq:det_equiv}). We obtain
\begin{align*}
\breve K_n(z,\bs\lambda)
&=
\frac{1}{n}
\operatorname{tr}
\left[
\breve{\mbf R}(z)
\bs\Psi(\bs\lambda)\mbf W^2
\bs\Psi^*(\bs\lambda)
\right]
\\
&\approx
-\frac{1}{zn}
\operatorname{tr}
\left[
(\mbf I+\mbf H(z))^{-1}
\bs\Psi(\bs\lambda)\mbf W^2
\bs\Psi^*(\bs\lambda)
\right]
\\
&=
-\frac{1}{n}
\sum_{t=1}^n
\frac{
g(\theta_t)h(\bs\lambda,\theta_t)
}{
z-c_n\breve M_n(z,\theta_t)
}.
\end{align*}
Equivalently,
\begin{equation}\label{eq:K_hat_new}
\breve K_n(z,\bs\lambda)
\approx
\frac{1}{n}
\sum_{t=1}^n
\frac{
g(\theta_t)h(\bs\lambda,\theta_t)
}{
c_n\breve M_n(z,\theta_t)-z
}.
\end{equation}
Equation~(\ref{eq:K_hat_new}) is the finite-sample analogue of
equation~(\ref{eq:K}) in Theorem~\ref{theorem1}.

Thus, the approximating system consists of
\[
\breve M_n(z,\theta_t)
=
\frac{1}{p}
\sum_{k=1}^p
\frac{
g(\theta_t)h(\bs\lambda_k,\theta_t)
}{
1+\breve K_n(z,\bs\lambda_k)
},
\]
together with
\[
\breve{\mathcal S}_n(z)
\approx
\frac{1}{n}
\sum_{t=1}^n
\frac{1}{
c_n\breve M_n(z,\theta_t)-z
}
\]
and
\[
\breve K_n(z,\bs\lambda)
\approx
\frac{1}{n}
\sum_{t=1}^n
\frac{
g(\theta_t)h(\bs\lambda,\theta_t)
}{
c_n\breve M_n(z,\theta_t)-z
}.
\]
As $p,n\to\infty$ with $p/n\to c$, the empirical distribution of
$\bs\lambda_1,\ldots,\bs\lambda_p$ converges to $F^{A,\Sigma}$ and the
Riemann sums over the Fourier frequencies converge to their corresponding
integrals. The preceding system therefore leads naturally to the limiting
equations (\ref{eq:s}) and (\ref{eq:K}) of Theorem~\ref{theorem1}.

\medskip

\noindent{\it Completion of the convergence argument.} The heuristic derivation above is made rigorous in two main steps.
\begin{enumerate}
\item[(P.1)]
For every fixed $z\in\mbb C^+$ and every sequence of Hermitian matrices
$\mbf C_n$ with uniformly bounded operator norm, establish
\[
\frac{1}{zn}
\operatorname{tr}
\left[
(\mbf I+\mbf H(z))^{-1}\mbf C_n
\right]
+
\frac{1}{n}
\operatorname{tr}
\left[
\breve{\mbf R}(z)\mbf C_n
\right]
\stackrel{a.s.}{\longrightarrow}0.
\]
This validates the deterministic-equivalent approximation used above.

\item[(P.2)]
Establish existence, uniqueness, and continuity of the solution of the
limiting system and use these properties to identify the limit of the
finite-sample approximating equations. More precisely, one proceeds as
follows:
\begin{enumerate}
\item[(S.1)]
Construct an event $\Omega_0$ with $\mathbb P(\Omega_0)=1$ on which the
required convergence statements hold.

\item[(S.2)]
For every $\omega\in\Omega_0$, extract a subsequence $(n_\ell)$ for which $\breve K_{z,n_\ell}(\bs\lambda,\omega)\rightarrow K_g(\bs\lambda,z)$ uniformly in $\bs\lambda$ and pointwise in $z\in\mbb C^+$, and show that the limit is the unique Stieltjes-kernel solution of
equation~(\ref{eq:K}).

\item[(S.3)]
Use uniqueness of the limiting solution to extend the convergence from the subsequence to the full sequence: $\breve K_n(z,\bs\lambda,\omega) \rightarrow K_g(z,\bs\lambda)$.

\item[(S.4)]
Finally, substitute the convergence of $\breve K_n$ into the approximating equation for $\breve{\mathcal S}_n$ to obtain $\breve{\mathcal S}_n(z,\omega)\rightarrow\mathcal S_g(z)$ for $z\in\mbb C^+$.
\end{enumerate}
\end{enumerate}

\subsubsection{Extension to non-Gaussian innovations}

We next indicate how the Gaussian assumption can be relaxed. The argument proceeds by truncating the innovations and subsequently centering and rescaling them. One first shows that these modifications do not affect the limiting spectral distribution. It therefore remains to compare the Gaussian model with a model having non-Gaussian innovations with matching first and second moments. To this end, consider
\[
\bar{\mathbf S}_n
=
\frac{1}{n}
\sum_{j=1}^p
\bar Y_j\bar Y_j^*,
\]
where $\bar Y_j=\mathbf W\boldsymbol{\Psi}^*(\boldsymbol{\lambda}_j)\xi_j$, $\xi_j=(\xi_{j1},\ldots,\xi_{jn})^T\in\mathbb C^n$, and the variables $\xi_{jk}$ are iid, have mean zero and variance one, and possess a finite fourth moment, but are not necessarily Gaussian. The goal is to show that $\bar{\mathbf S}_n$ and its Gaussian counterpart $\breve{\mathbf S}_n$ have the same limiting spectral distribution. Let
\[
s_n(z)
=
\frac{1}{n}
\operatorname{tr}
\bigl\{(\breve{\mathbf S}_n-z\mathbf I)^{-1}\bigr\},
\qquad
\bar s_n(z)
=
\frac{1}{n}
\operatorname{tr}
\bigl\{(\bar{\mathbf S}_n-z\mathbf I)^{-1}\bigr\}
\]
denote the corresponding Stieltjes transforms. After establishing tightness of the two sequences of empirical spectral distributions, it is enough to show that their Stieltjes transforms have the same limit. This is accomplished in two steps. For every fixed $z\in\mathbb C^+$,
\begin{align}
\mathbb E\!\left[
\bar s_n(z)-s_n(z)
\right]
&\longrightarrow 0,
\label{eq:nonGaussian_mean}
\\
\bar s_n(z)-\mathbb E[\bar s_n(z)]
&\longrightarrow 0
\qquad\text{in probability}.
\label{eq:nonGaussian_concentration}
\end{align}
The first step is proved using the generalized Lindeberg principle of \cite{c06}. The Gaussian innovation vectors are replaced successively by
non-Gaussian vectors having the same first two moments, and the cumulative effect of these replacements on the Stieltjes transform is shown to be asymptotically negligible. The second step follows from a concentration argument based on McDiarmid's inequality. Together, these results show that the Gaussian and non-Gaussian models have the same limiting Stieltjes transform, and hence the same limiting spectral distribution.

\subsection{Proof of Theorem~\ref{thm:sufficient_cond_L2_consistency}}

For $g\in\mcal G$, $z\in\mcal Z$, and
$\bs\omega\in\Delta_J$, define
\begin{equation}\label{eq:f_g_z_omega}
f_g(z,\bs\omega)
=
\mcal S_g^0(z)
-
\frac{1}{2\pi}
\int_0^{2\pi}
\frac{d\theta}
{c\tilde M_g(z,\theta\mid\bs\omega)-z},
\end{equation}
where
\begin{equation}\label{eq:M_g_tilde_omega}
\tilde M_g(z,\theta\mid\bs\omega)
=
g(\theta)
\sum_{j=1}^J
\omega_j
\frac{
h(\bs\lambda_j^0,\theta)
}{
1+K_g^0(\bs\lambda_j^0,z)
}.
\end{equation}
At the true parameter,
\[
\tilde M_g(z,\theta\mid\bs\omega^0)
=
M_g^0(z,\theta),
\]
and hence
\[
f_g(z,\bs\omega^0)=0.
\]

Since $\mcal Z$ is contained in a compact subset of $\mbb C^+$, there
exists $a>0$ such that
\[
\Im(z)\geq a,
\qquad z\in\mcal Z.
\]
Together with the nonnegativity of $g$ and $h$ and the Stieltjes-kernel
properties of $K_g^0$, this implies that
$f_g(z,\bs\omega)$ is uniformly bounded over
$g\in\mcal G$, $z\in\mcal Z$, and $\bs\omega\in\Delta_J$.

Define the population $L^2$ criterion
\begin{equation}\label{eq:D_population}
D(\bs\omega)
=
\frac{1}{|\mcal G||\mcal Z|}
\sum_{g\in\mcal G}
\sum_{z\in\mcal Z}
|f_g(z,\bs\omega)|^2.
\end{equation}
In particular,
\[
D(\bs\omega^0)=0.
\]

The empirical criterion can be written as
\[
D_n(\bs\omega)
=
\frac{1}{|\mcal G||\mcal Z|}
\sum_{g\in\mcal G}
\sum_{z\in\mcal Z}
\left|
f_g(z,\bs\omega)
+
\varepsilon_n(g,z,\bs\omega)
\right|^2,
\]
where $\varepsilon_n(g,z,\bs\omega)$ collects the difference between
the empirical quantities entering the estimation criterion and their
population counterparts.

Note that $\hat{\mathcal S}_n^{(n)}\to\mathcal{S}_g(z)$ and $\hat K_g^{(n)}(z,\bs\lambda^0_j)\to K_g(z,\bs\lambda_j^0)$, both almost surely. This
together with $c_n\to c$ and the finiteness of $\mcal G$ and $\mcal Z$ imply
\begin{equation}\label{eq:residual_g_z_omega_limit}
\sup_{\bs\omega\in\Delta_J}
\max_{g\in\mcal G}
\max_{z\in\mcal Z}
|\varepsilon_n(g,z,\bs\omega)|
\stackrel{a.s.}{\longrightarrow}0.
\end{equation}
Since $f_g$ is uniformly bounded, it follows that
\begin{equation}\label{eq:uniform_D_convergence}
\sup_{\bs\omega\in\Delta_J}
|D_n(\bs\omega)-D(\bs\omega)|
\stackrel{a.s.}{\longrightarrow}0.
\end{equation}

We next examine the population criterion in a neighborhood of
$\bs\omega^0$. Recall that any admissible perturbation of the population
weights can be written as
\[
\bs\omega
=
\bs\omega^0+\mbf B\bs\eta,
\]
where $\bs\eta\in\mbb R^{J-1}$ and
$\bs\omega^0+\mbf B\bs\eta\in\Delta_J$.

For $g\in\mcal G$ and $z\in\mcal Z$, recall the vector
\[
\mbf a_g^0(z)
=
\frac{1}{2\pi}
\int_0^{2\pi}
\frac{
g(\theta)\mbf v_g^0(z,\theta)
}{
\{cM_g^0(z,\theta)-z\}^2
}\,d\theta,
\]
where
\[
\mbf v_g^0(z,\theta)
=
\left(
\frac{
h(\bs\lambda_j^0,\theta)
}{
1+K_g^0(\bs\lambda_j^0,z)
}
\right)_{j=1}^J.
\]
Differentiating $f_g$ at $\bs\omega^0$ along the tangent space of the
simplex gives
\[
\nabla_{\bs\eta}
f_g(z,\bs\omega^0)
=
-c\,\mbf B^T\mbf a_g^0(z),
\]
up to the convention used in the definition of
$\mcal M_{\mcal G,\mcal Z}$ for absorbing the factor $c$.
Consequently, a Taylor expansion around $\bs\omega^0$ yields
\begin{equation}\label{eq:D_Taylor}
D(\bs\omega^0+\mbf B\bs\eta)
=
c^2
\bs\eta^T
\mbf B^T
\mcal M_{\mcal G,\mcal Z}
(\bs\Lambda_J^0,\bs\omega^0)
\mbf B
\bs\eta
+
O(\|\bs\eta\|^3),
\end{equation}
again with the factor $c^2$ omitted if it is incorporated into the
definition of $\mcal M_{\mcal G,\mcal Z}$.

By assumption,
\[
\mbf B^T
\mcal M_{\mcal G,\mcal Z}
(\bs\Lambda_J^0,\bs\omega^0)
\mbf B
\]
is positive definite. Hence there exist constants $\rho>0$ and
$C>0$ such that
\begin{equation}\label{eq:local_identification}
D(\bs\omega)
\geq
C\|\bs\omega-\bs\omega^0\|^2
\end{equation}
whenever
$\bs\omega\in\Delta_J$ and
$\|\bs\omega-\bs\omega^0\|\leq\rho$.
In particular, for every $0<\epsilon<\rho$,
\[
\inf_{\substack{
\bs\omega\in\Delta_J\\
\epsilon\leq
\|\bs\omega-\bs\omega^0\|
\leq\rho}}
D(\bs\omega)
>0.
\]

Combining this separation with the uniform convergence in
(\ref{eq:uniform_D_convergence}) shows that, almost surely for all
sufficiently large $n$, the empirical criterion $D_n$ is strictly larger
on the annulus
\[
\left\{
\bs\omega\in\Delta_J:
\epsilon\leq
\|\bs\omega-\bs\omega^0\|
\leq\rho
\right\}
\]
than at $\bs\omega^0$. Thus, any local minimizer of $D_n$ lying in the
$\rho$-neighborhood of $\bs\omega^0$ satisfies
\[
\|\widehat{\bs\omega}-\bs\omega^0\|<\epsilon
\]
for all sufficiently large $n$. Since $\epsilon>0$ is arbitrary,
\[
\widehat{\bs\omega}
\stackrel{a.s.}{\longrightarrow}
\bs\omega^0.
\]
This proves local consistency. If, in addition, $\bs\omega^0$ is the
unique minimizer of $D$ over the entire simplex $\Delta_J$, then
(\ref{eq:uniform_D_convergence}) yields consistency of the global
minimum-distance estimator.

\subsection{Proof of Proposition
\ref{prop:L2_consistency_independent_observations}}

For independent observations, $X_t=\mbf A_0 Z_t$,
$\Sigma=\mbf A_0^2=\operatorname{Var}(X_t)$, and hence the true joint spectral distribution reduces to
\[
F^0\equiv F^\Sigma
=
\sum_{j=1}^J
\omega_j^0\delta_{\lambda_j^0},
\]
where $\lambda_1^0,\ldots,\lambda_J^0$ are distinct nonnegative real numbers. Since $\mcal G=\{g_0\}$ with $g_0(\theta)\equiv1$, we have $K_{g_0}^0(z,\lambda_j^0,)=\lambda_j^0 S^0(z)$, and
\[
M_{g_0}^0(z,\theta)
=
\sum_{j=1}^J
\omega_j^0
\frac{\lambda_j^0}
{1+\lambda_j^0S^0(z)},
\]
where $S^0=S_{g_0}^0$. In particular, $M_{g_0}^0$ does not depend on $\theta$. Recall that
\[
\mbf a_{g_0}^0(z)
=
\frac{1}{2\pi}
\int_0^{2\pi}
\frac{
\mbf v^0(z)
}{
\{cM_{g_0}^0(z)-z\}^2
}\,d\theta
=
\frac{\mbf v^0(z)}
{\{cM_{g_0}^0(z)-z\}^2},
\]
where
\[
\mbf v^0(z)
=
\left(
\frac{\lambda_j^0}
{1+\lambda_j^0S^0(z)}
\right)_{j=1}^J.
\]
Consequently,
\[
\mcal M_{\mcal G,\mcal Z}
(\bs\Lambda^0,\bs\omega^0)
=
\frac{1}{|\mcal Z|}
\sum_{z\in\mcal Z}
\frac{
\mbf v^0(z)\mbf v^0(z)^*
}{
|cM_{g_0}^0(z)-z|^4
}.
\]
Since $\mcal Z$ is contained in a compact subset of $\mbb C^+$,
the factors $|cM_{g_0}^0(z)-z|^{-4}$ are bounded above and bounded
away from zero. Hence it suffices to show that
\[
\tilde{\mcal M}_{\mcal Z}
=
\frac{1}{|\mcal Z|}
\sum_{z\in\mcal Z}
\mbf v^0(z)\mbf v^0(z)^*
\]
is positive definite on the tangent space of the simplex. Using
\[
\frac{\lambda_j^0}
{1+\lambda_j^0S^0(z)}
=
\frac{1}{S^0(z)}
-
\frac{1}{\{S^0(z)\}^2}
\frac{1}{1/S^0(z)+\lambda_j^0},
\]
we can write
\[
\mbf v^0(z)
=
\frac{1}{S^0(z)}\mbf 1_J
-
\frac{1}{\{S^0(z)\}^2}
\tilde{\mbf v}^0(z),
\]
where
\[
\tilde{\mbf v}^0(z)
=
\left(
\frac{1}
{1/S^0(z)+\lambda_j^0}
\right)_{j=1}^J.
\]
Since $\mbf B^T\mbf 1_J=0$, it is therefore enough to establish positive definiteness of the Gram matrix generated by $\tilde{\mbf v}^0(z)$. Choose $J$ points $z_1,\ldots,z_J\in\mcal Z$ such that
\[
u_k:=\frac{1}{S^0(z_k)},
\qquad k=1,\ldots,J,
\]
are distinct. Consider the $J\times J$ matrix
\[
\mbf C
=
\left(
\frac{1}{u_k+\lambda_j^0}
\right)_{1\leq j,k\leq J}.
\]
This is a Cauchy matrix. Since the $\lambda_j^0$'s are distinct, the $u_k$'s are distinct, and $u_k+\lambda_j^0\neq0$, its determinant is 
\[
\det(\mbf C)
=
\frac{
\displaystyle
\prod_{1\leq j<\ell\leq J}
(\lambda_\ell^0-\lambda_j^0)
\prod_{1\leq k<r\leq J}
(u_r-u_k)
}{
\displaystyle
\prod_{j=1}^J\prod_{k=1}^J
(u_k+\lambda_j^0)
},
\]
up to an irrelevant sign. Hence, $\det(\mbf C)\neq0$. It follows that $\tilde{\mbf v}^0(z_1),\ldots,\tilde{\mbf v}^0(z_J)$ span $\mbb C^J$, and consequently $\mbf B^T\mcal M_{\mcal G,\mcal Z}(\bs\Lambda_J^0,\bs\omega^0)\mbf B$ is positive definite. This establishes the required local
identifiability condition.

\subsection{Proof of Proposition
\ref{prop:ARMA_L2_consistency_condition}}

Suppose that the spectral-profile Gram matrix
\[
\mbf G(\bs\Lambda^0)
=
\frac{1}{2\pi}
\int_0^{2\pi}
\mbf h^0(\theta)
\mbf h^0(\theta)^T\,d\theta
\]
is positive definite, where
\[
\mbf h^0(\theta)
=
\left(
h(\bs\lambda_1^0,\theta),
\ldots,
h(\bs\lambda_J^0,\theta)
\right)^T.
\]
We show that one can construct a finite collection of localized weight functions $\mcal G$ for which $\mbf B^T\mcal M_{\mcal G,\mcal Z}(\bs\Lambda^0,\bs\omega^0)\mbf B$ is positive definite. For $\delta>0$, let $N_\delta=\frac{2\pi}{\delta}$ be an integer and let $\theta_{k,\delta}$, $k=1,\ldots,N_\delta$, be equally spaced points in $[0,2\pi)$. Define
\[
\mbf R_\delta(\bs\Lambda^0)
=
\frac{1}{N_\delta}
\sum_{k=1}^{N_\delta}
\mbf h^0(\theta_{k,\delta})
\mbf h^0(\theta_{k,\delta})^T.
\]
Since this is a Riemann-sum approximation to $\mbf G(\bs\Lambda^0)$, we have $\mbf R_\delta(\bs\Lambda^0) \rightarrow \mbf G(\bs\Lambda_J^0)$ as $\delta\to0$. Hence there exists $\delta_1>0$ such that
$\mbf R_\delta(\bs\Lambda_J^0)$ is positive definite whenever
$0<\delta<\delta_1$.

Let $\underline g$ be a nonnegative, symmetric, continuously differentiable function with compact support, satisfying
\[
\int_{\mbb R}\underline g(u)\,du=1.
\]
For $k=1,\ldots,N_\delta$, define the localized weight functions
\[
g_{k,\delta}(\theta)
=
\underline g
\left(
\frac{\theta-\theta_{k,\delta}}{\delta}
\right),
\]
with the natural periodic interpretation on $[0,2\pi]$, and set $\mcal G_\delta=\{g_{k,\delta}:k=1,\ldots,N_\delta\}$. Let
\[
h_\star(\bs\Lambda^0)
=
\sup_{\theta\in[0,2\pi]}
\max_{1\leq j\leq J}
h(\bs\lambda_j^0,\theta).
\]
Since $\mcal Z$ is contained in a compact subset of $\mbb C^+$, there exist constants $0<\underline a<\overline a<\infty$ such that $\Im(z)\geq\underline a$, $|z|\leq\overline a$, and $z\in\mcal Z$. The Stieltjes-kernel properties imply that the denominators appearing below are uniformly bounded away from zero. In particular, $|cM_{g_{k,\delta}}^0(z,\theta)-z|\geq\underline a$. Moreover, because the support of $g_{k,\delta}$ has length of order
$\delta$,
\[
\max_{\substack{
1\leq j\leq J\\
1\leq k\leq N_\delta\\
z\in\mcal Z}}
\left|
K_{g_{k,\delta}}^0
(z,\bs\lambda_j^0)
\right|
=
O(\delta).
\]
It follows that $\mbf v_{g_{k,\delta}}^0(z,\theta)=\mbf h^0(\theta)+O(\delta)$ uniformly on the support of $g_{k,\delta}$, where the error is uniform in $j$, $k$, and $z\in\mcal Z$.

Recall that
\[
\mbf a_g^0(z)
=
\frac{1}{2\pi}
\int_0^{2\pi}
\frac{
g(\theta)\mbf v_g^0(z,\theta)
}{
\{cM_g^0(z,\theta)-z\}^2
}\,d\theta.
\]
Using the change of variables
\[
u
=
\frac{\theta-\theta_{k,\delta}}{\delta}
\]
and the smoothness of the quantities in the integrand gives, uniformly
in $k$ and $z\in\mcal Z$,
\begin{equation}\label{eq:a_g_localized}
\mbf a_{g_{k,\delta}}^0(z)
=
\delta\,
q_{k,\delta}(z)
\mbf h^0(\theta_{k,\delta})
+
O(\delta^2),
\end{equation}
where
\[
q_{k,\delta}(z)
=
\frac{1}{
2\pi
\{cM_{g_{k,\delta}}^0
(z,\theta_{k,\delta})-z\}^2
}.
\]
The preceding bounds imply that $|q_{k,\delta}(z)|$ is uniformly bounded above and bounded away from
zero for sufficiently small $\delta$.

Consequently, for fixed $z\in\mcal Z$,
\begin{align}
&
\frac{1}{N_\delta}
\sum_{k=1}^{N_\delta}
\mbf a_{g_{k,\delta}}^0(z)
\mbf a_{g_{k,\delta}}^0(z)^*
\nonumber\\
&\qquad
=
\frac{\delta^2}{N_\delta}
\sum_{k=1}^{N_\delta}
|q_{k,\delta}(z)|^2
\mbf h^0(\theta_{k,\delta})
\mbf h^0(\theta_{k,\delta})^T
+
O(\delta^3).
\label{eq:M_G_Z_Lambda_ARMA_delta_approx}
\end{align}
Since the weights $|q_{k,\delta}(z)|^2$ are uniformly positive and
bounded, positive definiteness of
$\mbf R_\delta(\bs\Lambda^0)$ implies that the leading term on the
right-hand side of
(\ref{eq:M_G_Z_Lambda_ARMA_delta_approx}) is positive definite.
For sufficiently small $\delta$, the $O(\delta^3)$ remainder cannot
destroy this positive definiteness.

Averaging over the finite set $\mcal Z$ therefore yields
\[
\mcal M_{\mcal G_\delta,\mcal Z}
(\bs\Lambda^0,\bs\omega^0)
=
\frac{1}{|\mcal Z|}
\sum_{z\in\mcal Z}
\frac{1}{N_\delta}
\sum_{k=1}^{N_\delta}
\mbf a_{g_{k,\delta}}^0(z)
\mbf a_{g_{k,\delta}}^0(z)^*
>0
\]
for all sufficiently small $\delta$. In particular, $\mbf B^T\mcal M_{\mcal G_\delta,\mcal Z}(\bs\Lambda_J^0,\bs\omega^0)\mbf B>0$. Thus, the localized weight functions in $\mcal G_\delta$ provide the required local identifiability of the population weights.

\subsection{Proof of Proposition
\ref{prop:L2_consistency_AR1}}

For the univariate AR$(1)$ process corresponding to the support point $\bs\lambda_j^0=(\alpha_j^0,\sigma_j^0)$, we have $x_{j,t}=\alpha_j^0x_{j,t-1}+\sigma_j^0z_{j,t}$, and hence
\begin{equation}\label{eq:AR_autocov}
\gamma_\ell(\bs\lambda_j^0)
=
\frac{(\sigma_j^0)^2}
{1-(\alpha_j^0)^2}
(\alpha_j^0)^\ell,
\qquad
\ell=0,1,\ldots.
\end{equation}
Recall that $\mbf G(\bs\Lambda^0)=\bs\gamma_0^0(\bs\gamma_0^0)^T+2\sum_{\ell=1}^\infty\bs\gamma_\ell^0(\bs\gamma_\ell^0)^T$. For any $\mbf x\in\mbb R^J$, $\mbf x^T\mbf G(\bs\Lambda_J^0)\mbf x=(\mbf x^T\bs\gamma_0^0)^2+
2\sum_{\ell=1}^\infty(\mbf x^T\bs\gamma_\ell^0)^2$. It therefore suffices to show that $\bs\gamma_0^0,\ldots,\bs\gamma_{J-1}^0$ span $\mbb R^J$. Define
\[
\mbf D_0
=
\operatorname{diag}
\left(
\frac{(\sigma_1^0)^2}{1-(\alpha_1^0)^2},
\ldots,
\frac{(\sigma_J^0)^2}{1-(\alpha_J^0)^2}
\right).
\]
By assumption $\sigma_j^0>0$ and $|\alpha_j^0|<1$, so
$\mbf D_0$ is positive definite. From
(\ref{eq:AR_autocov}),
\[
[\bs\gamma_0^0:\cdots:\bs\gamma_{J-1}^0]
=
\mbf D_0
\begin{bmatrix}
1 & \alpha_1^0 & (\alpha_1^0)^2
  & \cdots & (\alpha_1^0)^{J-1}\\
1 & \alpha_2^0 & (\alpha_2^0)^2
  & \cdots & (\alpha_2^0)^{J-1}\\
\vdots & \vdots & \vdots
  & & \vdots\\
1 & \alpha_J^0 & (\alpha_J^0)^2
  & \cdots & (\alpha_J^0)^{J-1}
\end{bmatrix}.
\]
The matrix on the right is a Vandermonde matrix. Since $\alpha_1^0,\ldots,\alpha_J^0$ are distinct, its determinant is $\prod_{1\leq j<k\leq J}(\alpha_k^0-\alpha_j^0)\neq0$. Hence $[\bs\gamma_0^0:\cdots:\bs\gamma_{J-1}^0]$ is nonsingular. Consequently, for every nonzero $\mbf x\in\mbb R^J$, at least one of $\mbf x^T\bs\gamma_0^0,\ldots, \mbf x^T\bs\gamma_{J-1}^0$
is nonzero. Therefore,
\[
\mbf x^T
\mbf G(\bs\Lambda^0)
\mbf x
=
(\mbf x^T\bs\gamma_0^0)^2
+
2\sum_{\ell=1}^\infty
(\mbf x^T\bs\gamma_\ell^0)^2
>0
\]
and hence $\mbf G(\bs\Lambda^0)$ is positive definite, which proves the proposition.

\end{appendix}

\begin{acks}[Acknowledgments]
    Most of this work was done when the first author was at the University of California, Davis. The second author's research was partially funded by National Science Foundation grant DMS-2515821.
\end{acks}

\newpage

\section{Supplementary Material}

{\color{black}
\subsection{Sufficient condition for consistency for general AR processes }
Let $X_t = \mbf{A}_1 X_{t-1} +\dots + \mbf{A}_r X_{t-r} + \Sigma^{1/2}Z_t$, where $\Sigma$, $\mbf{A}_1,\dots,\mbf{A}_r$ are symmetric/Hermitian,  $\|\mbf{A}_k\| < 1, ~k=1,\dots,r$, and there exists
an orthogonal/unitary matrix $\mbf{U}$ such that, $\mbf{U}^* \mbf{A}_k \mbf{U} = \mbox{diag}(\alpha_{1,k},\ldots,\alpha_{p,k}), ~k=1,\dots,r$, and $\mbf{U}^* \Sigma \mbf{U}
= \mbox{diag}(\sigma_1^2,\ldots,\sigma_p^2)$, for $\alpha_{1,k},\ldots,\alpha_{p,k} \in \mbb{R}, ~k=1,\dots,r$
and $\sigma_1,\ldots,\sigma_p \geq 0$. 

Denote $\lambda = (\alpha_1,\dots,\alpha_r, \sigma)$ and suppose that the joint spectrum of $(\mbf{A}_1,\dots,\mbf{A}_r,\Sigma)$ is given by $F^{A_1,\dots,A_r,\Sigma} = \sum_{j=1}^J \omega_j \delta_{\lambda_j^0}$, where $\lambda_j^0 = (\alpha_{j,1}^0,\dots,\alpha_{j,r}^0,\sigma_j^0)$ with $\max_{1\leq j \leq J}|\alpha_{j,k}^0| < 1, ~k=1,\dots,r$, and $\min_{1\leq j \leq J}\sigma_j^0 > 0$. The univariate causal AR(r) process 
\[
\phi(\mathbb{B})x_{j,t} = \sigma_j^0z_{j,t}
\]
has an autocovariance function $\gamma_h(\lambda_j^0)$ of the form (see \cite{bd92} page 94)
\[
\gamma_h(\lambda_j^0) = \sum_{i=1}^{k}\sum_{\ell=0}^{r_i-1}\beta_{i,\ell}^{(j)}h^\ell\xi_i^{-k}, \qquad h\geq 0,
\]
where $\xi_i, i=1,\dots,k$ are the zeros (possibly complex) of $\phi(\mathbb{B})$, $r_i$ is the multiplicity of $\xi_i$, $\phi(\mathbb{B}) = \alpha_{j,1}^0\mathbb{B} + \dots + \alpha_{j,r}^0\mathbb{B}^r$, $\{z_{j,t}\}$ are white noise processes, $\mathbb{B}$ is the back-shift operator, and $\beta^{(j)}_{i,\ell}$ are uniquely determined through solving \[
\gamma_h(\lambda_j^0) - \alpha^0_{j,1}\gamma_{h-1}(\lambda_j^0) - \dots - \alpha^0_{j,r}\gamma_{h-r}(\lambda_j^0) = (\sigma_j^0)^2, \qquad 0\leq h \leq r. 
\] 
If we let $\gamma_h = (\gamma_h(\lambda_j^0))_{1\leq j \leq J}$ then 
\begin{align*}
\mbf G(\Lambda_J) =  \left(\gamma_0\gamma_0^T + 2\sum_{j=1}^{\infty}\gamma_j\gamma_j^T\right) 
\end{align*}
is nonsingular if and only if $ \boldsymbol{\Gamma} = [\gamma_0:\dots: \gamma_{J-1}]$ is nonsingular.
      
}

\subsection{Simulation Results: Tables}
\subsubsection{Absolute loss}

\begin{center}
\begin{table}[H] \label{table:case1.L1_AR}
\centering
\resizebox{\textwidth}{!}{%
\begin{tabular}{|>{\centering\arraybackslash}m{1.05cm}|>{\centering\arraybackslash}m{1.05cm}|>{\centering\arraybackslash}m{1.05cm}|>{\centering\arraybackslash}m{1.05cm}|>{\centering\arraybackslash}m{1.05cm}|>{\centering\arraybackslash}m{.89cm}|>{\centering\arraybackslash}m{.89cm}|>{\centering\arraybackslash}m{1.05cm}|>{\centering\arraybackslash}m{1cm}|>{\centering\arraybackslash}m{.89cm}|>{\centering\arraybackslash}m{.89cm}|>{\centering\arraybackslash}m{.89cm}|>{\centering\arraybackslash}m{.89cm}|}
	\hline
	\multicolumn{2}{|c|}{AR}&\multicolumn{2}{|c|}{Sigma}&&\multicolumn{2}{|c|}{p400n1600}&\multicolumn{3}{|c|}{p600n2400}\\
	\hline
	evals&weights&evals&weights&statistics&4g&8g&4g&8g&12g \\
	\hline
	0.5&1&(1,2)&(0.5,0.5)&mean median sd&0.0378    0.0377    0.0218
	&0.0122    0.0078    0.0119&0.0266    0.0203    0.0285& 0.0090    0.0067    0.0081& 0.0069    0.0054    0.0055\\
	\hline
	0.5&1&(1,2)&(0.75, 0.25)&mean median sd&0.0286    0.0274    0.0212&	 0.0102    0.0076    0.0085&	 0.0192    0.0159    0.0141&	 0.0084    0.0080    0.0060&	0.0059    0.0045    0.0043\\
	\hline
	(-0.5, 0.8)&(0.5, 0.5)&1&1&mean median sd&0.0155    0.0140    0.0092&	0.0095    0.0092    0.0054&	0.0096    0.0082    0.0070&	0.0061    0.0059    0.0034&	0.0053    0.0047    0.0037\\
	\hline
	(-0.5, 0.8)&(0.25,0.75)&1&1&mean median sd&0.0161    0.0152    0.0098&	0.0084    0.0074    0.0049&	0.0067    0.0052    0.0038&	0.0054    0.0053    0.0025&	 0.0057    0.0049    0.0033\\
	\hline
	(-0.5, 0.8)&(0.5,0.5)&(1,2)&(0.5,0.5)&mean median sd&0.1137    0.1199    0.0337&	0.0281    0.0205    0.0176&	 0.1096    0.1028    0.0230&	0.0255    0.0219    0.0143&	0.0140    0.0133    0.0090
	\\
	\hline
\end{tabular}
}
\caption{\textbf{Case 1}. Table contains mean, median, and standard deviation of  $L_1$ distance between true and estimated spectral cdf of AR(1) coefficient matrix, i.e. $d_{L_1}(F^A,\hat{F}^A)$ .}
\end{table}

\begin{table}[H] \label{table:case1.L1_Sig}
\centering
\resizebox{\textwidth}{!}{%
\begin{tabular}{|>{\centering\arraybackslash}m{1.05cm}|>{\centering\arraybackslash}m{1.05cm}|>{\centering\arraybackslash}m{1.05cm}|>{\centering\arraybackslash}m{1.05cm}|>{\centering\arraybackslash}m{1.05cm}|>{\centering\arraybackslash}m{.89cm}|>{\centering\arraybackslash}m{.89cm}|>{\centering\arraybackslash}m{1.05cm}|>{\centering\arraybackslash}m{1cm}|>{\centering\arraybackslash}m{.89cm}|>{\centering\arraybackslash}m{.89cm}|>{\centering\arraybackslash}m{.89cm}|>{\centering\arraybackslash}m{.89cm}|}
	\hline
	\multicolumn{2}{|c|}{AR}&\multicolumn{2}{|c|}{Sigma}&&\multicolumn{2}{|c|}{p400n1600}&\multicolumn{3}{|c|}{p600n2400}\\
	\hline
	evals&weights&evals&weights&statistics&4g&8g&4g&8g&12g \\
	\hline
	0.5&1&(1,2)&(0.5,0.5)&mean median sd&0.0435    0.0365    0.0284&	0.0329    0.0318    0.0185&	 0.0423    0.0308    0.0332&	0.0295    0.0225    0.0219&	0.0312    0.0262    0.0234\\
	\hline
	0.5&1&(1,2)&(0.75, 0.25)&mean median sd&0.0310    0.0231    0.0258&	 0.0295    0.0235    0.0200&	 0.0240    0.0185    0.0141&	 0.0222    0.0206    0.0085&	0.0212    0.0209    0.0094\\
	\hline
	(-0.5, 0.8)&(0.5, 0.5)&1&1&mean median sd&0.0344    0.0304    0.0240&	0.0215    0.0201    0.0151&	0.0191    0.0158    0.0120&	0.0137    0.0125    0.0062&	0.0149    0.0114    0.0124\\
	\hline
	(-0.5, 0.8)&(0.25,0.75)&1&1&mean median sd&0.0290    0.0248    0.0154&	0.0264    0.0249    0.0152&	0.0149    0.0132    0.0057&	0.0181    0.0143    0.0213&	0.0171    0.0147    0.0119\\
	\hline
	(-0.5, 0.8)&(0.5,0.5)&(1,2)&(0.5,0.5)&mean median sd&0.1640    0.1726    0.0284&	0.0697    0.0668    0.0270&	 0.1585    0.1499    0.0304&	0.0657    0.0642    0.0258&	0.0495    0.0454    0.0189\\
	\hline
\end{tabular}
}
\caption{\textbf{Case 1}. Table contains mean, median, and standard deviation of $L_1$ distance between true and estimated spectral cdf of $\Sigma$, i.e. $d_{L_1}(F^\Sigma,\hat{F}^\Sigma)$.}
\end{table}
\end{center}

\subsubsection{Quadratic loss}

\begin{center}
\begin{table}[H] \label{table:case1.L2_AR}
\centering
\resizebox{\textwidth}{!}{%
\begin{tabular}{|>{\centering\arraybackslash}m{1.05cm}|>{\centering\arraybackslash}m{1.05cm}|>{\centering\arraybackslash}m{1.05cm}|>{\centering\arraybackslash}m{1.05cm}|>{\centering\arraybackslash}m{1.05cm}|>{\centering\arraybackslash}m{.89cm}|>{\centering\arraybackslash}m{.89cm}|>{\centering\arraybackslash}m{1.05cm}|>{\centering\arraybackslash}m{1cm}|>{\centering\arraybackslash}m{.89cm}|>{\centering\arraybackslash}m{.89cm}|>{\centering\arraybackslash}m{.89cm}|>{\centering\arraybackslash}m{.89cm}|}
	\hline
	\multicolumn{2}{|c|}{AR}&\multicolumn{2}{|c|}{Sigma}&&\multicolumn{2}{|c|}{$p400n1600$}&\multicolumn{3}{|c|}{$p600n2400$}\\
	\hline
	values&weights&values&weights&statistics&$4g$&$8g$&$4g$&$8g$&$12g$ \\
	\hline
	0.5&1&(1,2)&(0.5,0.5)&mean median sd& 0.0709    0.1077    0.0539&	0.0108    0.0071    0.0145&	0.0729    0.1118    0.0512&	0.0082    0.0053    0.0096&	 0.0045    0.0032    0.0050\\
	\hline
	0.5&1&(1,2)&(0.75, 0.25)&mean median sd&0.0871    0.0991    0.0352&	 0.0083    0.0048    0.0090&	0.0810    0.0941    0.0387&	 0.0060    0.0044    0.0053&	0.0051    0.0032    0.0061\\
	\hline
	(-0.5, 0.8)&(0.5, 0.5)&1&1&mean median sd&0.0094    0.0067    0.0083&	 0.0055    0.0047    0.0042&	0.0062    0.0058    0.0056&	0.0044    0.0034    0.0044&	0.0035    0.0024    0.0028\\
	\hline
	(-0.5, 0.8)&(0.25,0.75)&1&1&mean median sd&0.0101    0.0099    0.0052&	0.0086    0.0069    0.0054&	0.0070    0.0059    0.0040&	0.0077    0.0087    0.0036&	0.0054    0.0045    0.0031\\
	\hline
	(-0.5, 0.8)&(0.5,0.5)&(1,2)&(0.5,0.5)&mean median sd&0.0240    0.0175    0.0195&	0.0166    0.0151    0.0105&	0.0164    0.0151    0.0098&	0.0114    0.0112    0.0062&	 0.0096    0.0091    0.0055\\
	\hline
\end{tabular}
}
\caption{\textbf{Case 1}. Table contains mean, median, and standard deviation of  $L_2$ distance between true and estimated spectral cdf of AR(1) coefficient matrix, i.e. $d_{L_2}(F^A,\hat{F}^A)$ .}
\end{table}
\end{center}

\begin{center}
\begin{table}[H] \label{table:case1.L2_AR-Sigma}
\centering
\resizebox{\textwidth}{!}{%
\begin{tabular}{|>{\centering\arraybackslash}m{1.05cm}|>{\centering\arraybackslash}m{1.05cm}|>{\centering\arraybackslash}m{1.05cm}|>{\centering\arraybackslash}m{1.05cm}|>{\centering\arraybackslash}m{1.05cm}|>{\centering\arraybackslash}m{.89cm}|>{\centering\arraybackslash}m{.89cm}|>{\centering\arraybackslash}m{1.05cm}|>{\centering\arraybackslash}m{1cm}|>{\centering\arraybackslash}m{.89cm}|>{\centering\arraybackslash}m{.89cm}|>{\centering\arraybackslash}m{.89cm}|>{\centering\arraybackslash}m{.89cm}|}
	\hline
	\multicolumn{2}{|c|}{AR}&\multicolumn{2}{|c|}{Sigma}&&\multicolumn{2}{|c|}{$p400n1600$}&\multicolumn{3}{|c|}{$p600n2400$}\\
	\hline
	values&weights&values&weights&statistics&$4g$&$8g$&$4g$&$8g$&$12g$ \\
	\hline
	0.5&1&(1,2)&(0.5,0.5)&mean median sd& 0.0948    0.0976    0.0127&	 0.0994    0.0981    0.0073&	0.0923    0.0932    0.0126&	0.1001    0.0994    0.0061&	0.1151    0.1144    0.0061\\
	\hline
	0.5&1&(1,2)&(0.75, 0.25)&mean median sd& 0.1154    0.1142    0.0115&	0.1010    0.1005    0.0050&	0.1114    0.1144    0.0103&	 0.1003    0.1013    0.0080&	0.1235    0.1253    0.0130\\
	\hline
	(-0.5, 0.8)&(0.5, 0.5)&1&1&mean median sd&0.0486    0.0529    0.0395&	0.0501    0.0285    0.0493&	0.0424    0.0514    0.0351&	0.0557    0.0741    0.0478&	0.0473    0.0100    0.0477\\
	\hline
	(-0.5, 0.8)&(0.25,0.75)&1&1&mean median sd&0.0676    0.0651    0.0104&	 0.0863    0.1227    0.0578&	 0.0641    0.0620    0.0077&	0.0998    0.1232    0.0518&	0.0757    0.1129    0.0554\\
	\hline
	(-0.5, 0.8)&(0.5,0.5)&(1,2)&(0.5,0.5)&mean median sd&0.0826    0.0802    0.0150&	0.1226    0.1201    0.0147&	0.0796    0.0798    0.0089&	0.1226    0.1216    0.0121&	 0.1149    0.1150    0.0069\\
	\hline
\end{tabular}
}
\caption{\textbf{Case 1}. Table contains mean, median, and standard deviation of $L_2$ distance between true and estimated spectral cdf of $\Sigma$, i.e. $d_{L_2}(F^\Sigma,\hat{F}^\Sigma)$.}\end{table}
\end{center}

\subsection{Simulation Results: Plots}
In what follows we present the plots of estimated marginal spectral cdf's of the coefficient matrices and the innovation matrix of the processes considered in Section 7 of the manuscript.
\newpage

\begin{figure}[H] \label{fig:Case5_AR}
\begin{center}
\resizebox{\textwidth}{!}{%
\begin{tabular}{cc}
$(p,n)=(400,1600)$ & $(p,n)=(400,800)$ \vspace{-32pt}\\ 
\includegraphics[width=2.9in, height=2.9in, bb = 50 50 650 650]{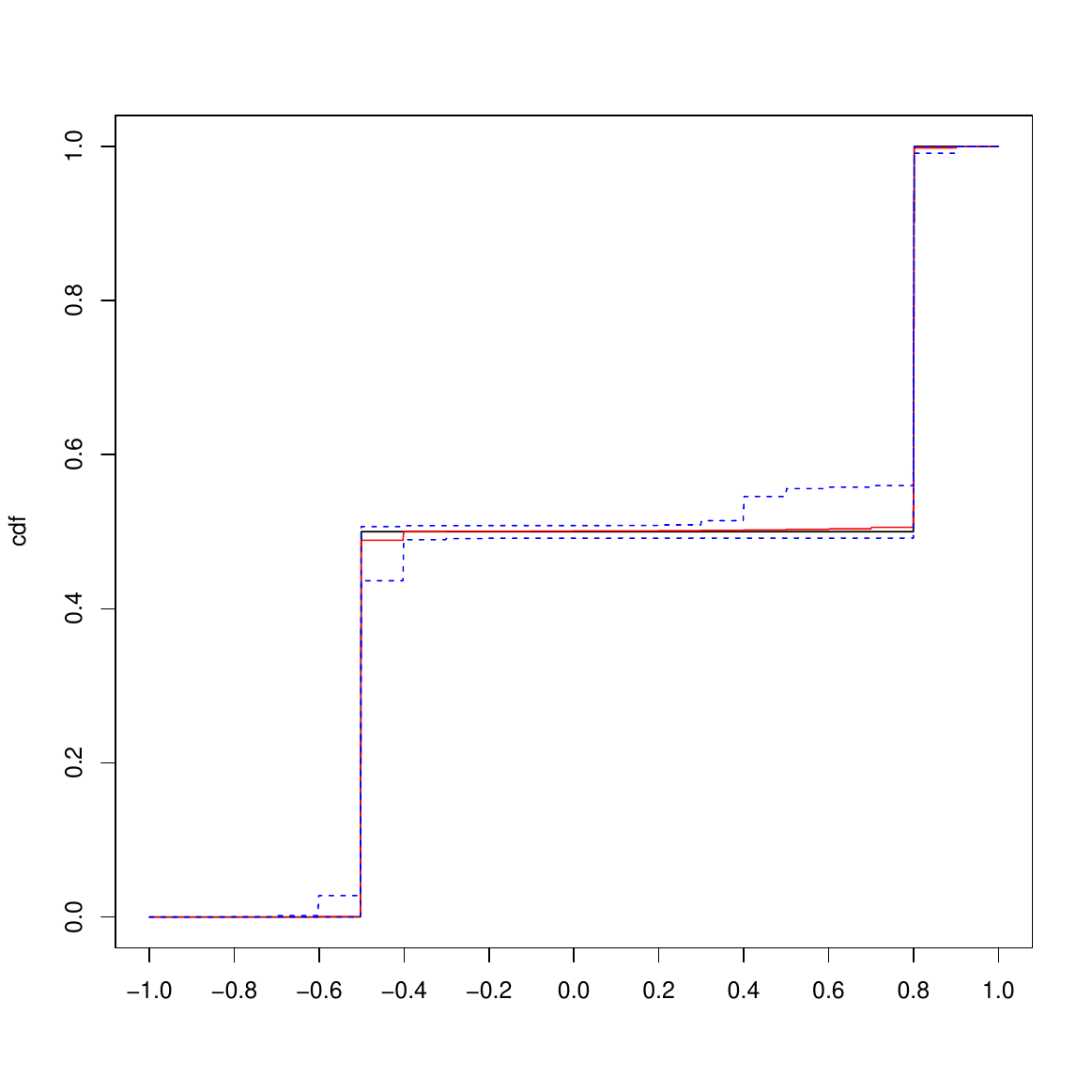} 
& \includegraphics[width=2.9in, height=2.9in, bb = 50 50 650 650]{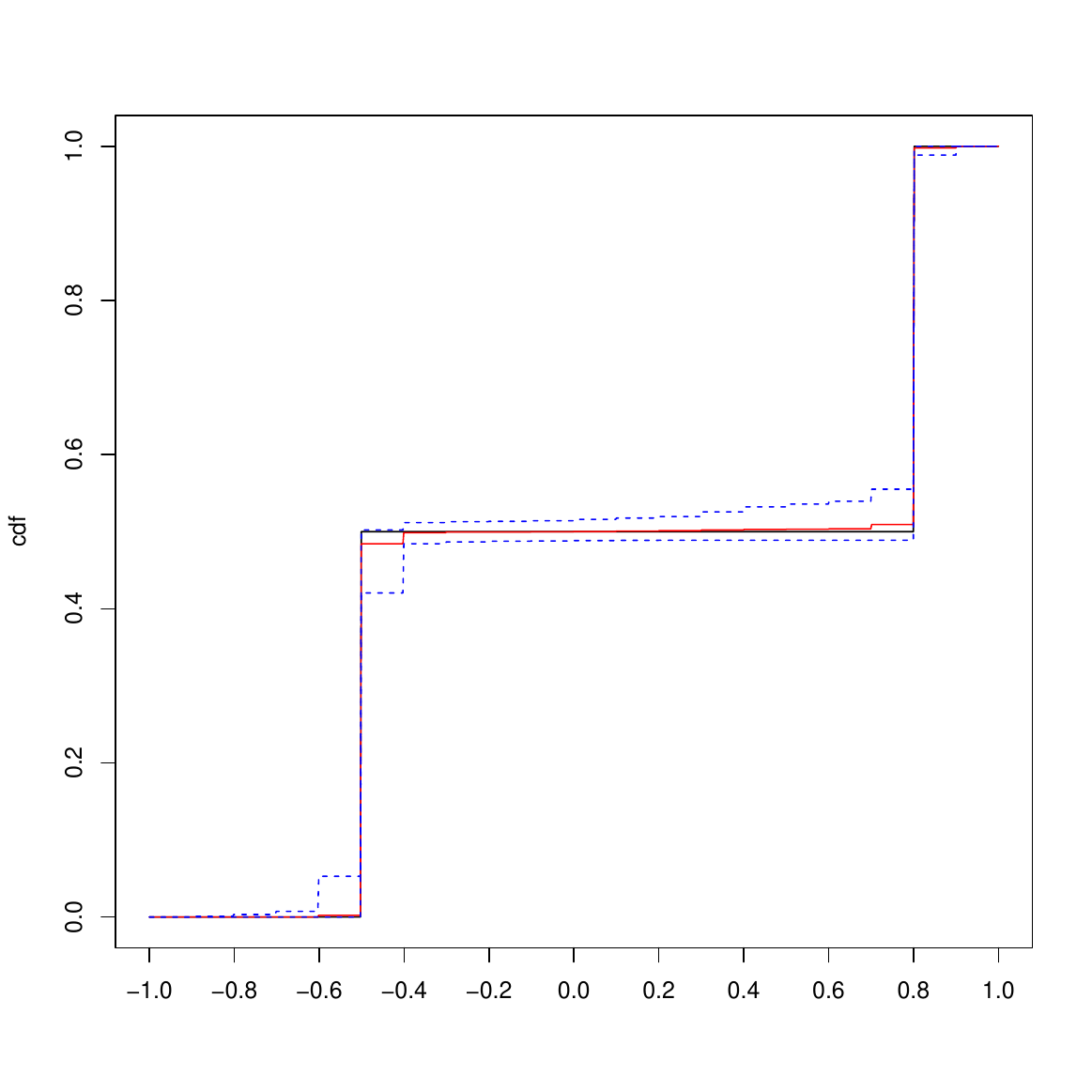}\vspace{7pt}\\  \vspace{-32pt}
$(p,n)=(200,800)$ & $(p,n)=(200,400)$\\ 
\hspace{10pt} \includegraphics[width=2.9in, height=2.9in, bb = 50 50 650 650]{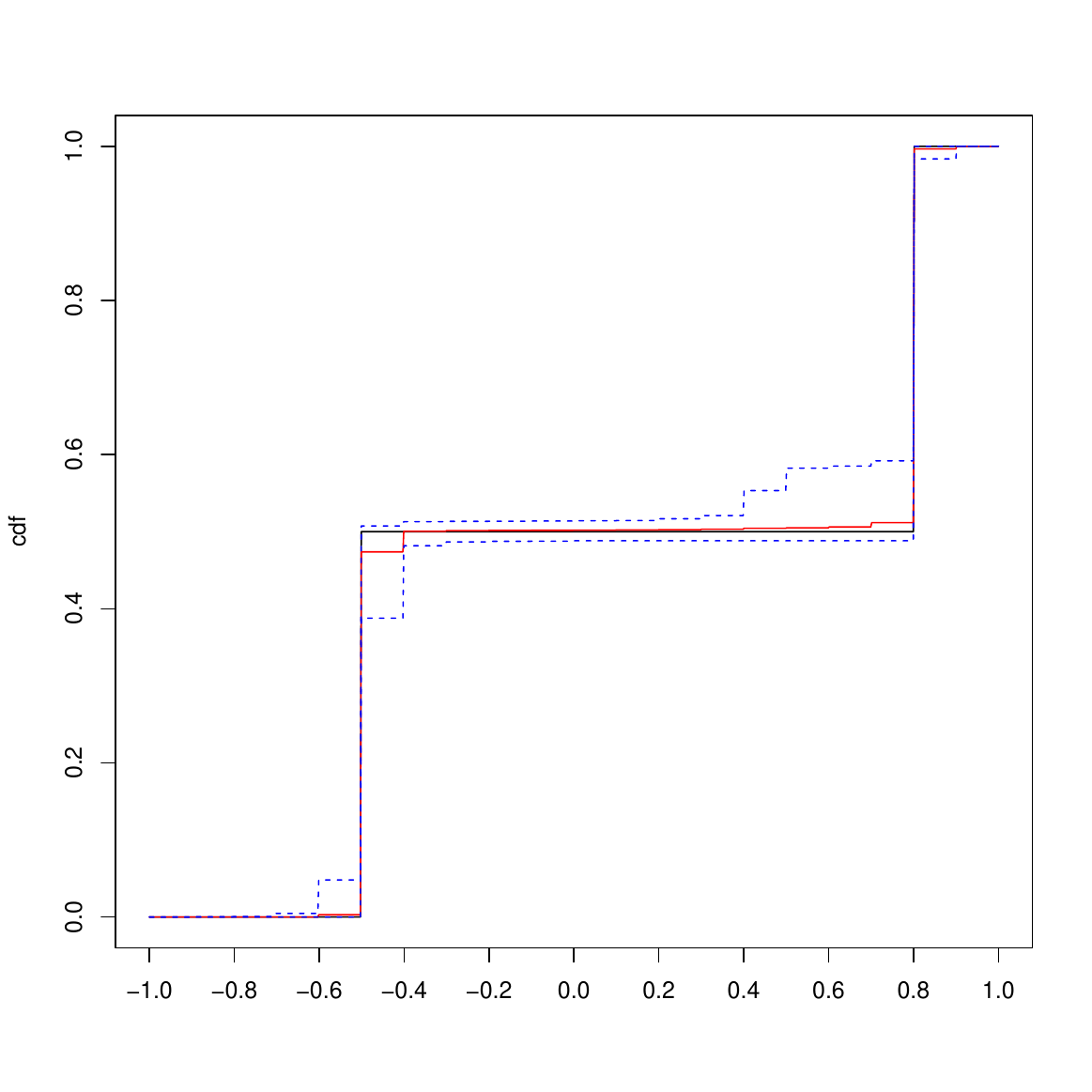}
& \includegraphics[width=2.9in, height=2.9in, bb = 50 50 650 650]{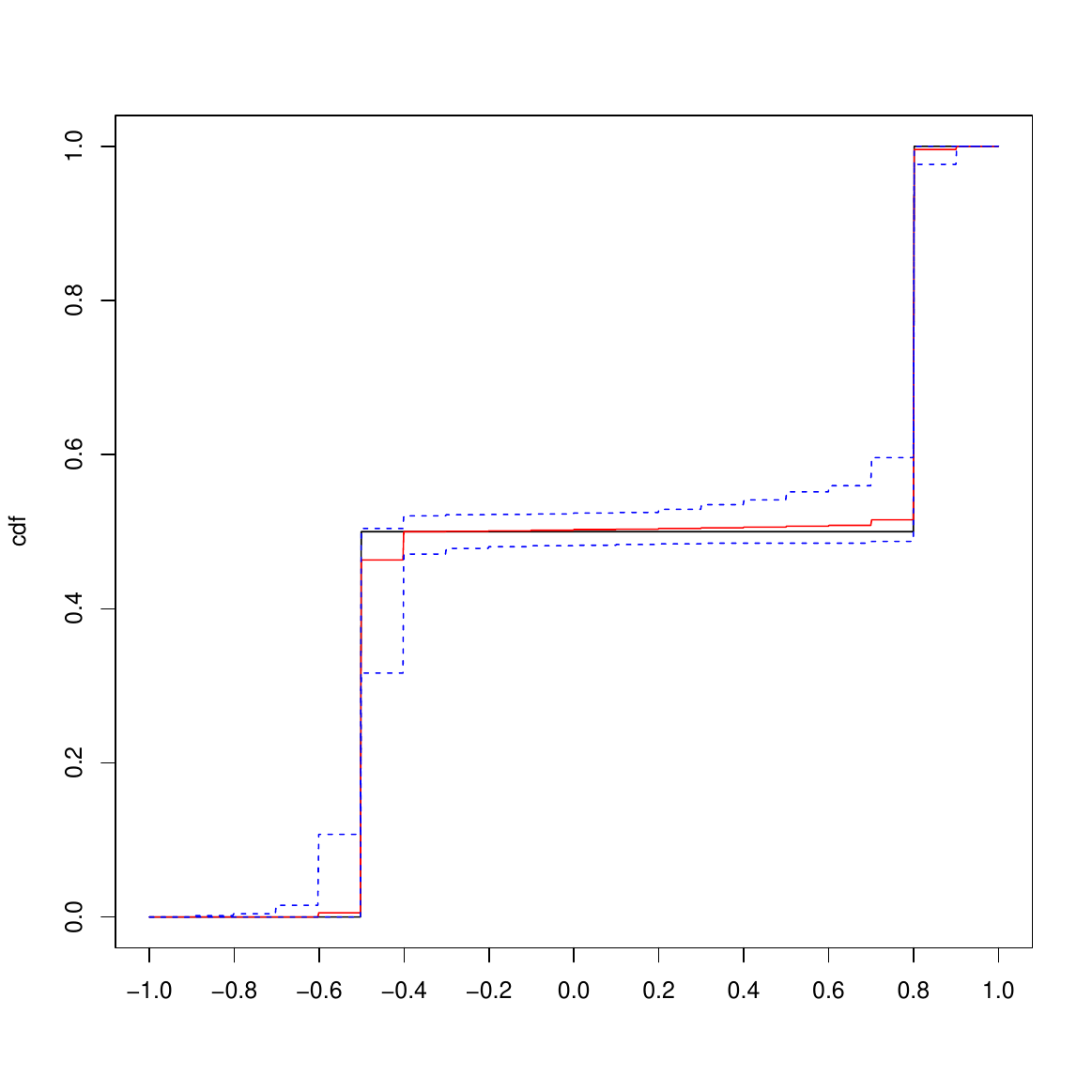} \vspace{7pt}\\ \vspace{-32pt}
$(p,n)=(200,100)$ & \\ 
\hspace{10pt}\includegraphics[width=2.9in, height=2.9in, bb = 50 50 650 650]{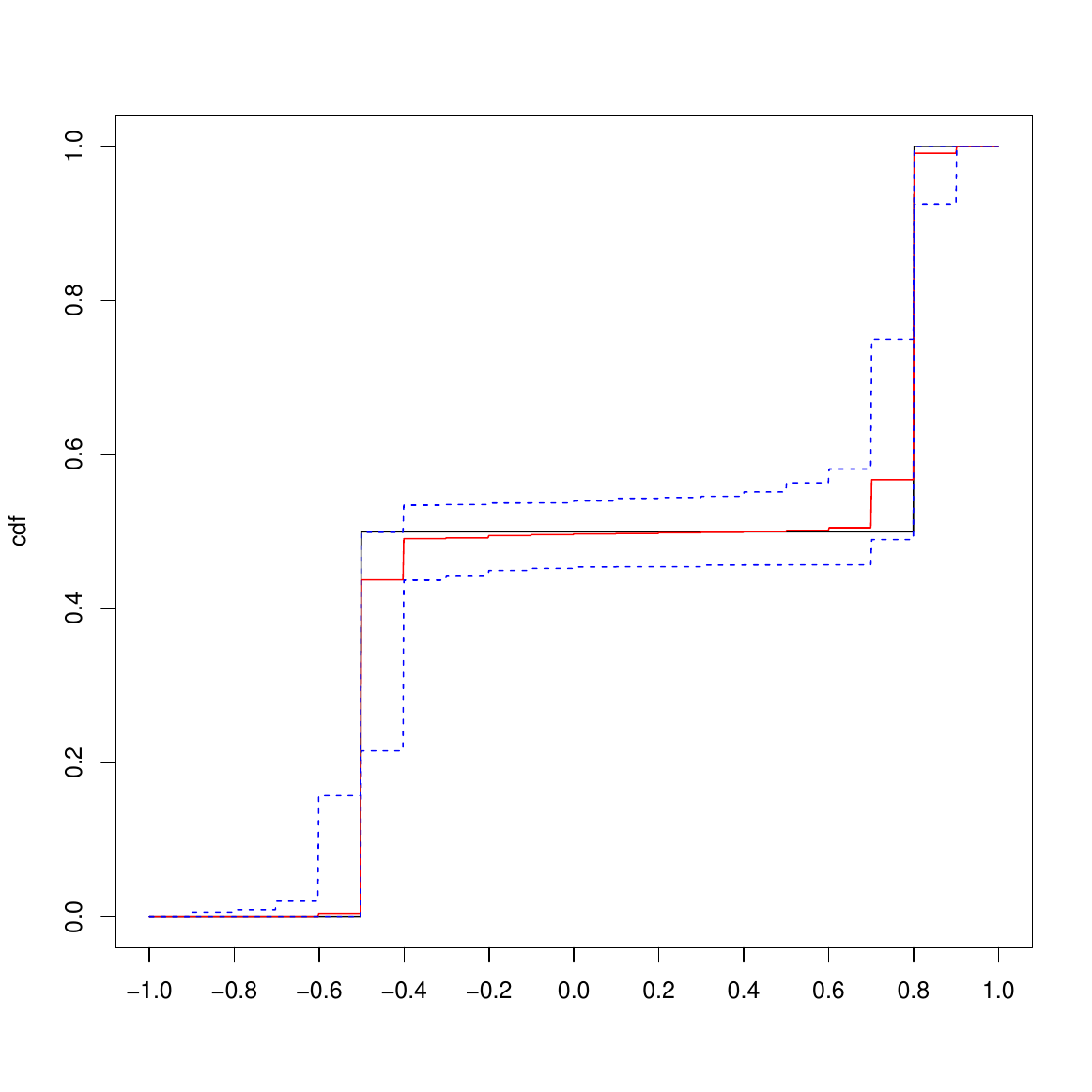}
& \\
\end{tabular}
}
\vspace{-1pt}
         \caption{Plot of Median and $90\%$ confidence band for spectral cdf of \textbf{AR coefficient matrix} corresponding to the \textbf{case 1.5}. \textbf{Solid Red} curve: median, \textbf{Dashed Blue} curve: $90\%$ confidence band, \textbf{Black Solid} curve: true spectral cdf}
\end{center}
\end{figure}

\begin{figure}[H] \label{fig:Case5_Sigma}
\begin{center}
\begin{tabular}{cc}
$(p,n)=(400,1600)$ & $(p,n)=(400,800)$ \vspace{-32pt}\\ 
\includegraphics[width=2.9in, height=2.9in, bb = 50 50 650 650]{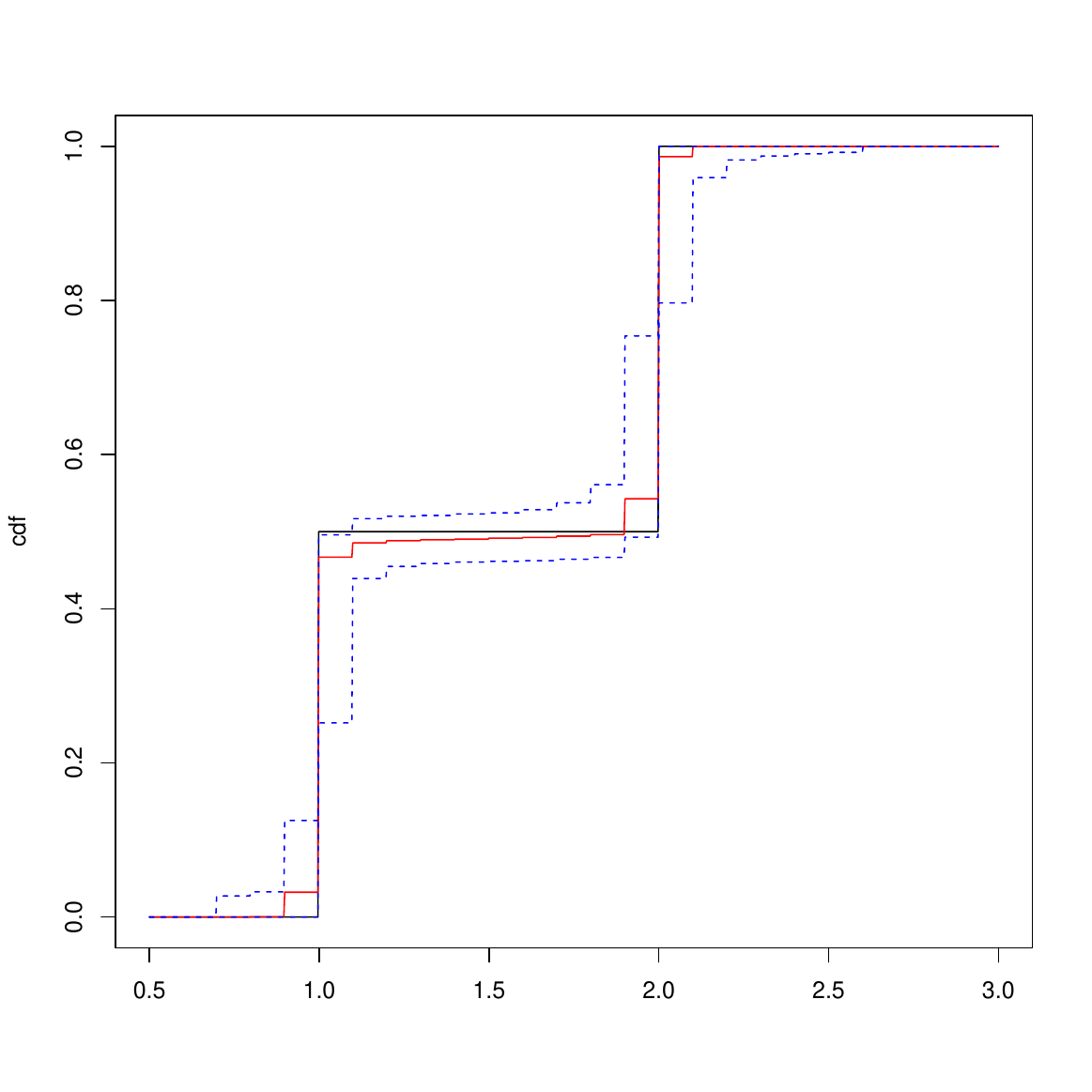}
& \includegraphics[width=2.9in, height=2.9in, bb = 50 50 650 650]{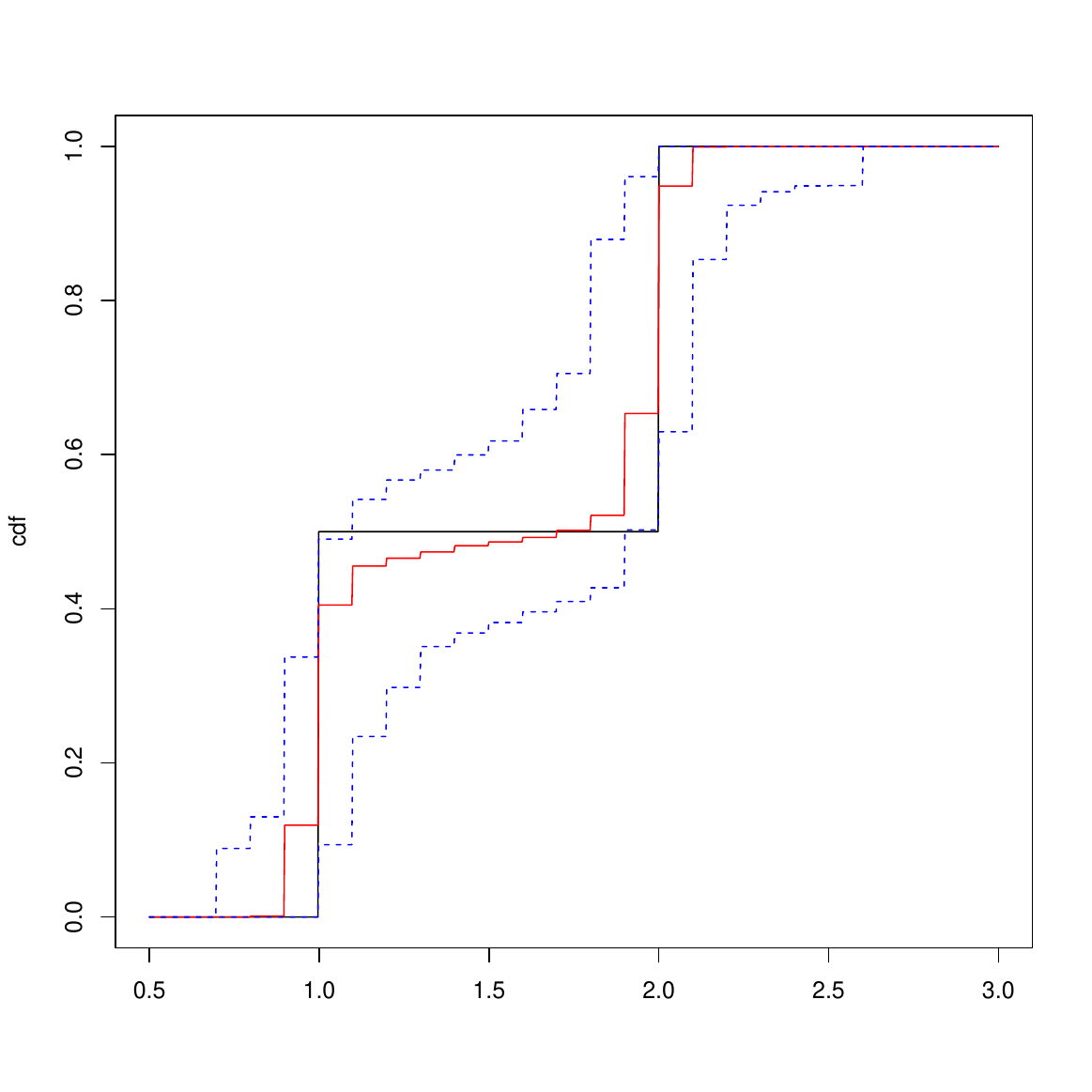}\vspace{7pt} \\ \vspace{-32pt}
$(p,n)=(200,800)$ & $(p,n)=(200,400)$\\ 
\hspace{10pt} \includegraphics[width=2.9in, height=2.9in, bb = 50 50 650 650]{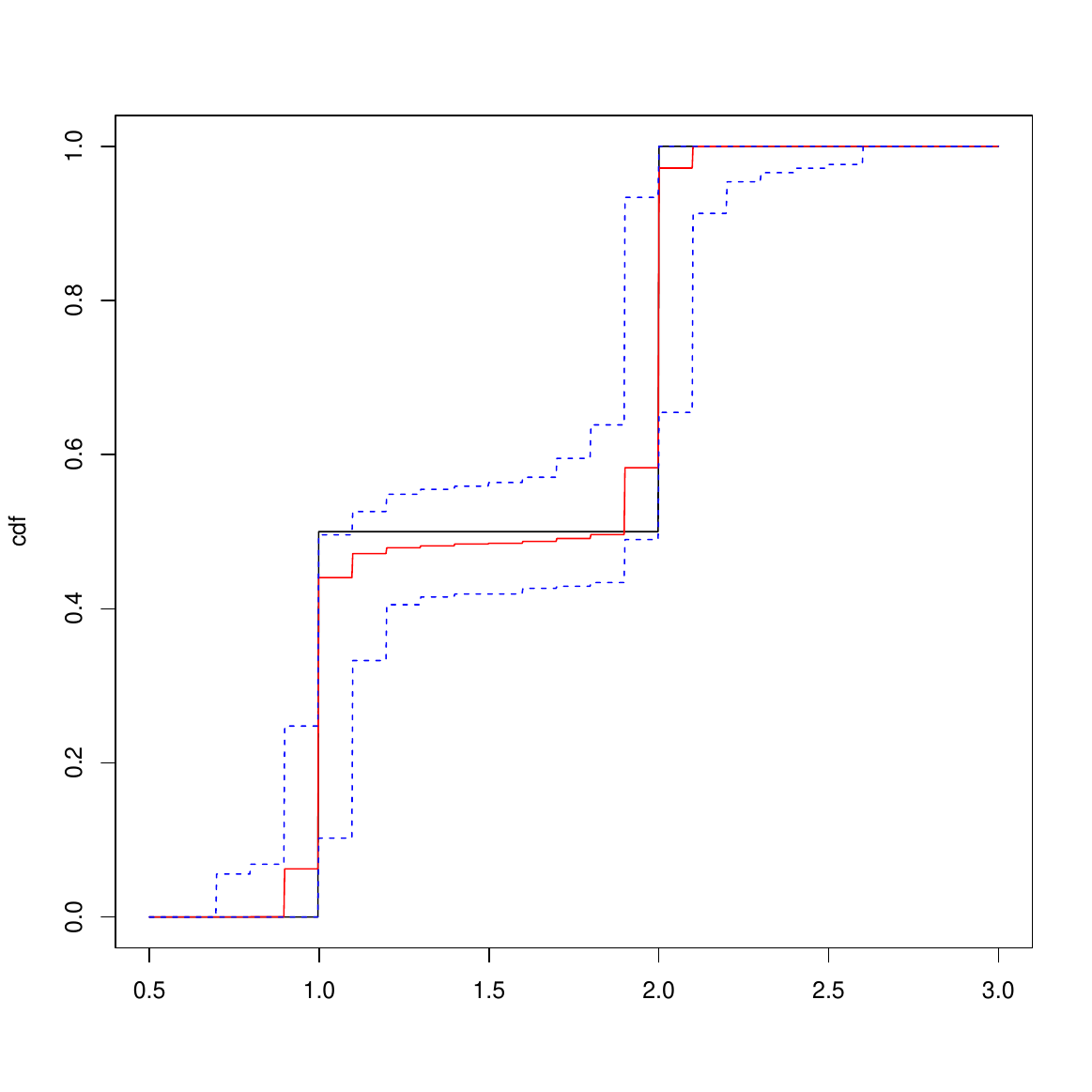} 
& \includegraphics[width=2.9in, height=2.9in, bb = 50 50 650 650]{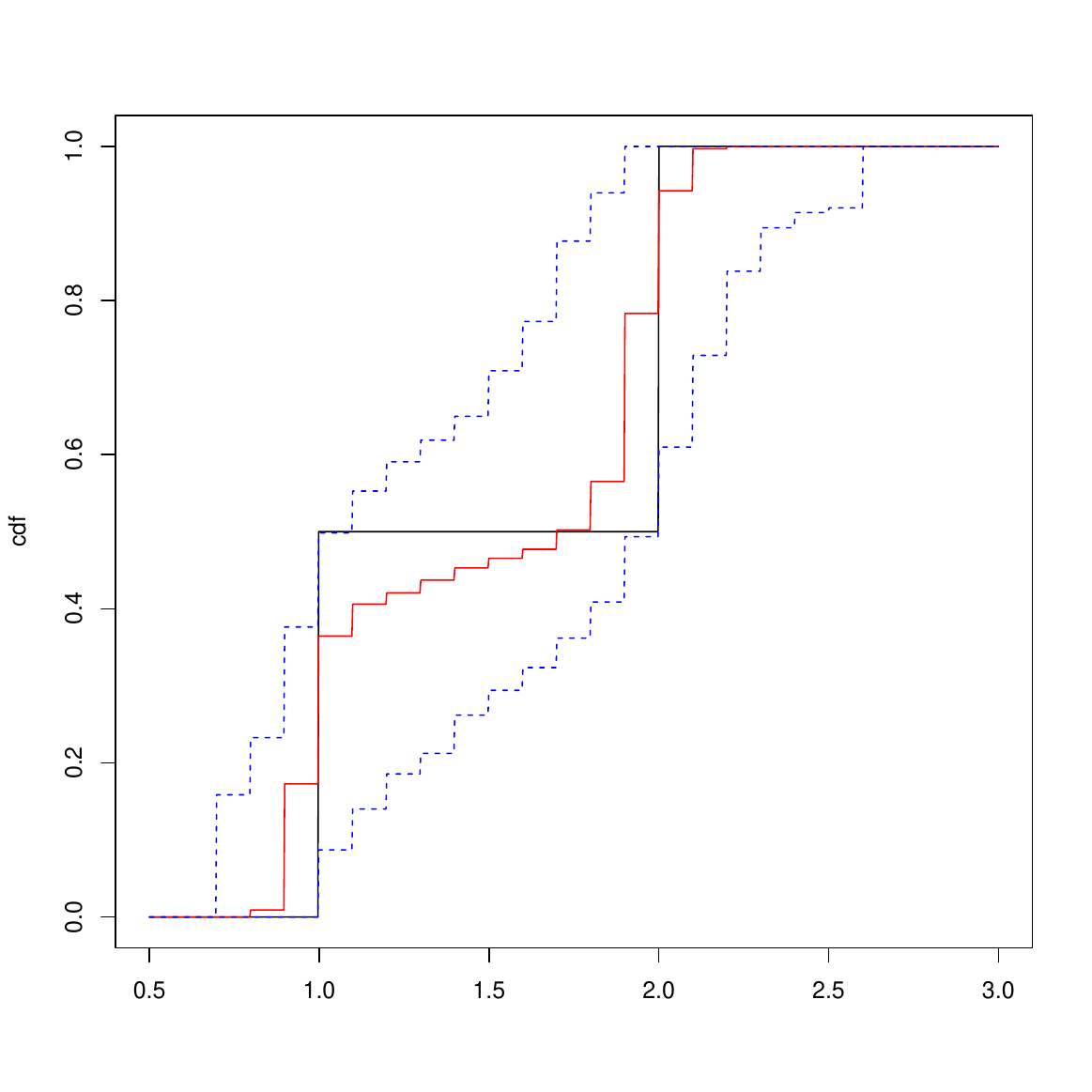}\vspace{7pt}\\ \vspace{-32pt}
$(p,n)=(200,100)$ & \\ \vspace{-10pt}
\hspace{10pt} \includegraphics[width=2.9in, height=2.9in, bb = 50 50 650 650]{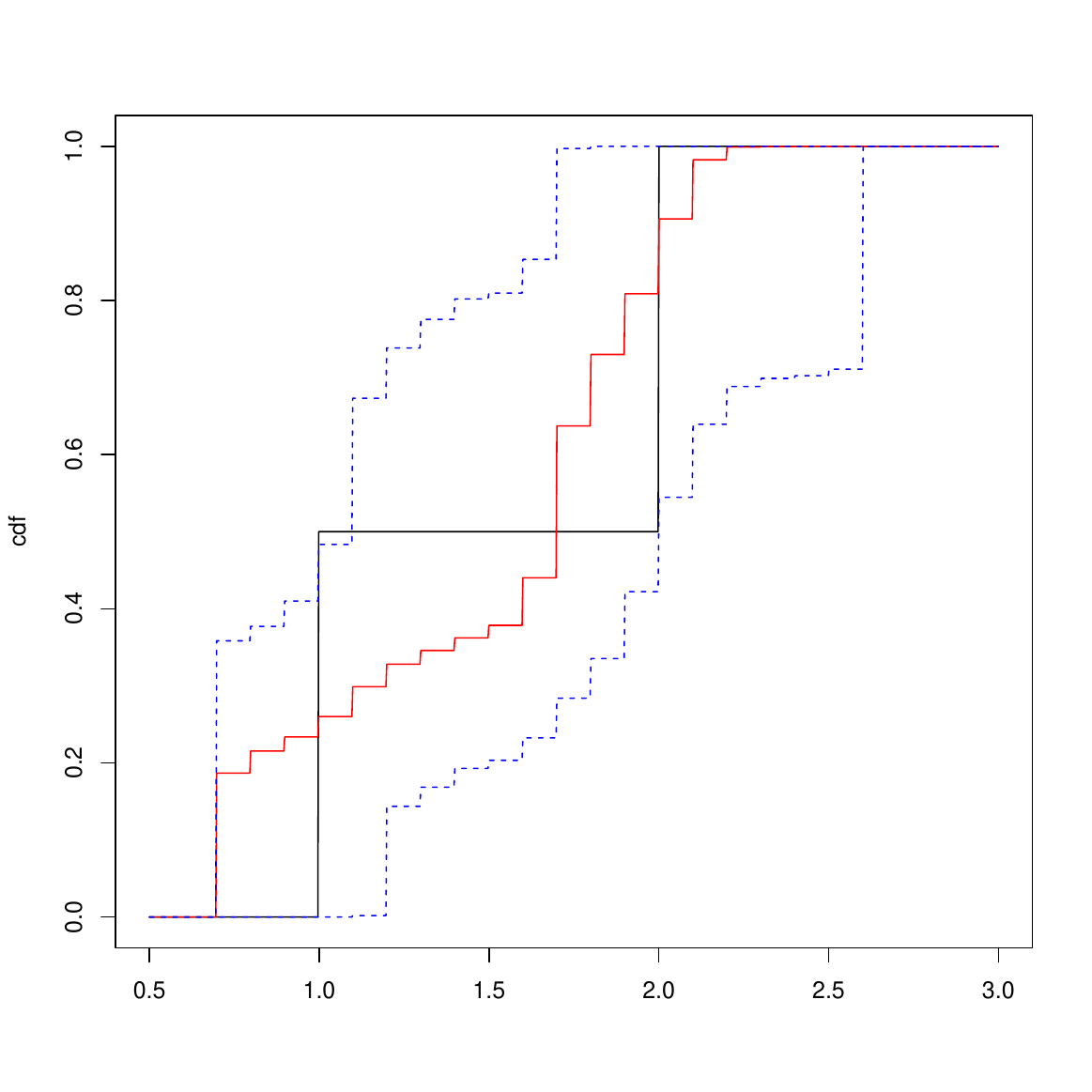}
& \\
\end{tabular}
\vspace{10pt}
         \caption{Plot of Median and $90\%$ confidence band for spectral cdf of $\boldsymbol{A_2}$ corresponding to the \textbf{case 1.5}. \textbf{Solid Red} curve: median, \textbf{Dashed Blue} curve: $90\%$ confidence band, \textbf{Black Solid} curve: true spectral cdf}
\end{center}
\end{figure}

\begin{figure}[H] \label{fig:ARMA_Sigma}
\begin{center}
\resizebox{\textwidth}{!}{%
\begin{tabular}{cc}
$(p,n)=(400,1600)$ & $(p,n)=(400,800)$ \\ \vspace{-20pt}
\includegraphics[width=2.9in, height=3in, bb = 50 50 650 650]{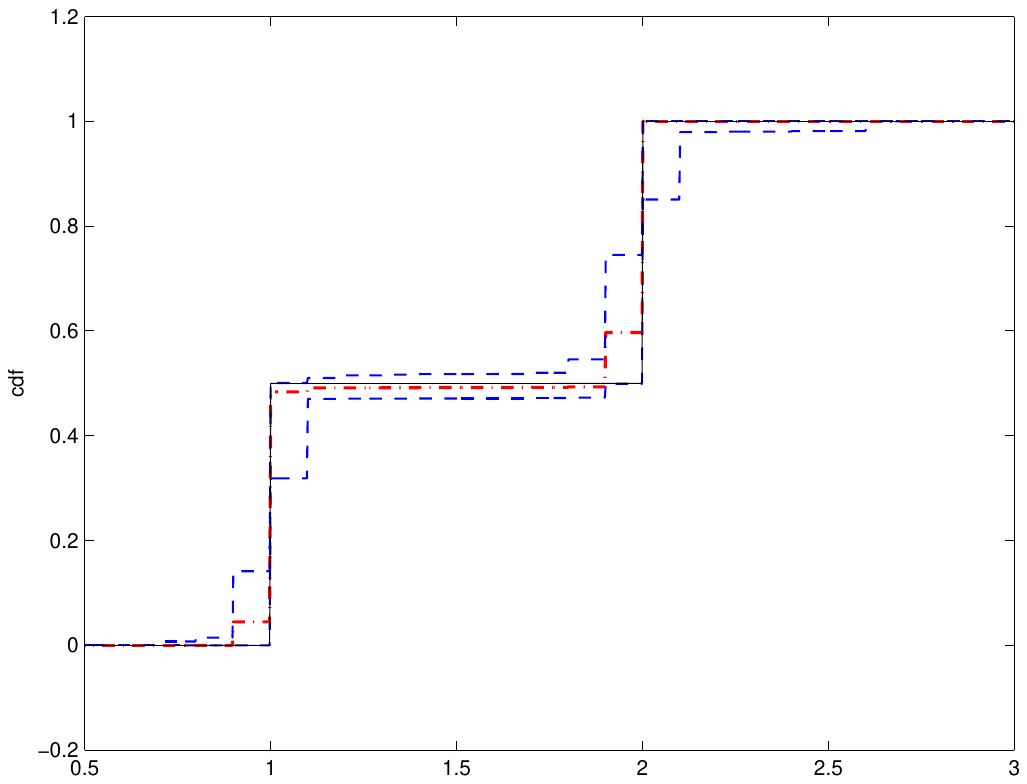}
& \includegraphics[width=2.9in, height=3in, bb = 50 50 650 650]{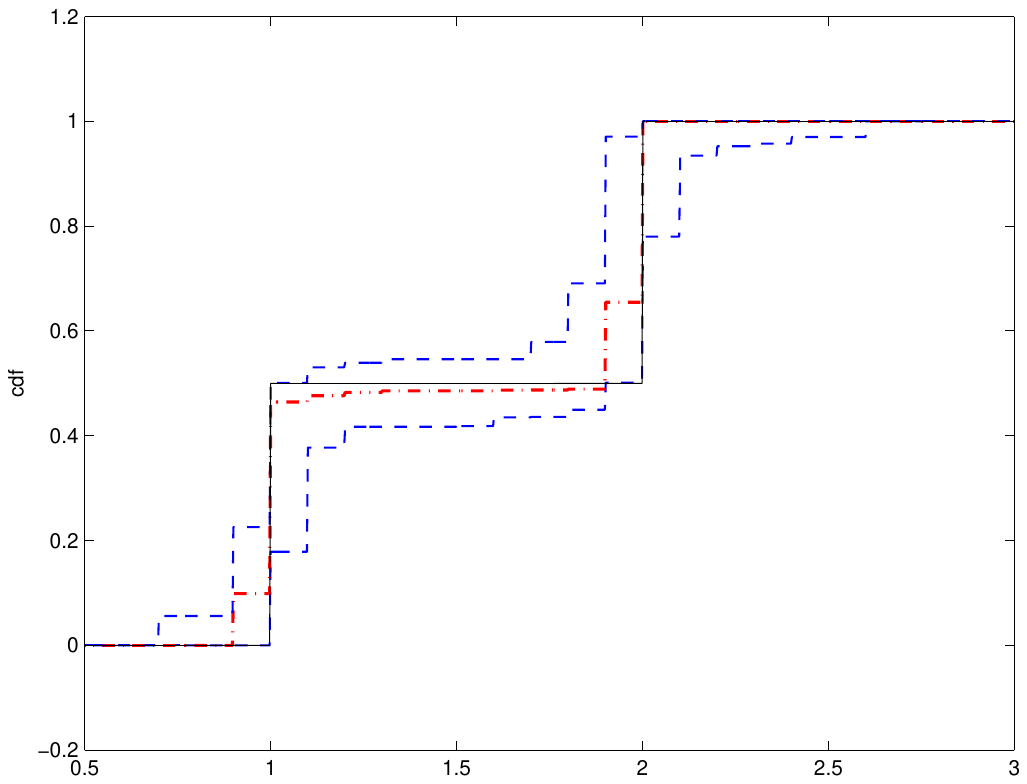}\vspace{-20pt}\\ \vspace{-10pt}
$(p,n)=(200,800)$ & $(p,n)=(200,400)$\vspace{-5pt}\\ \vspace{-10pt}
\includegraphics[width=2.9in, height=3in, bb = 50 50 650 650]{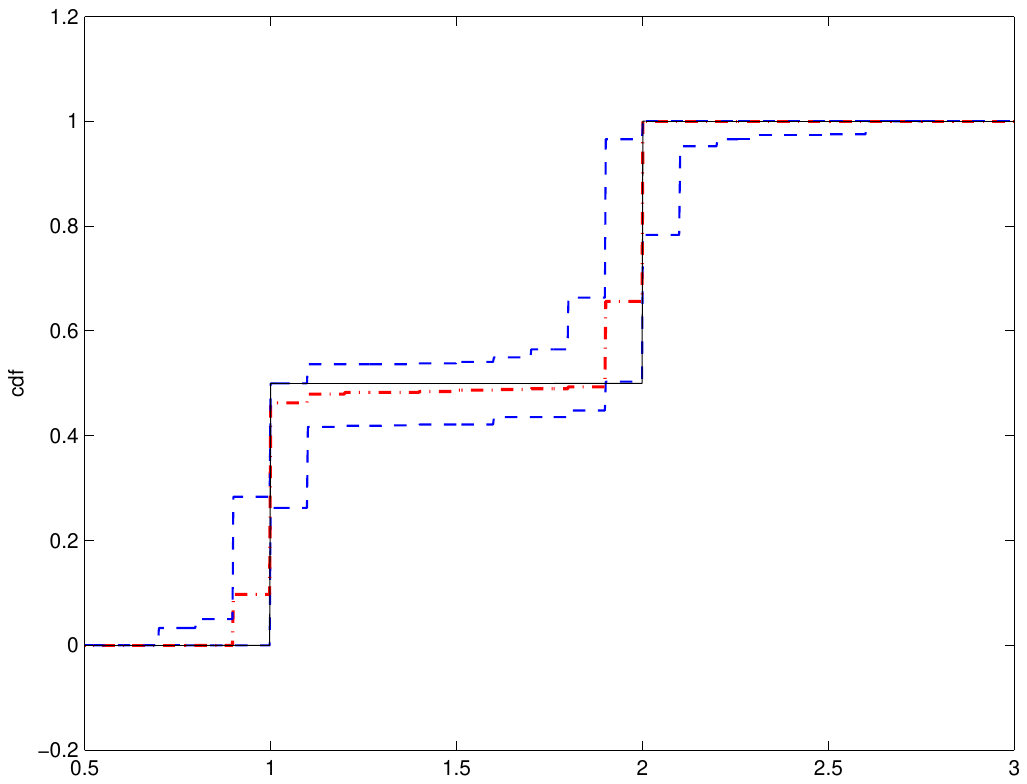}
& \includegraphics[width=2.9in, height=3in, bb = 50 50 650 650]{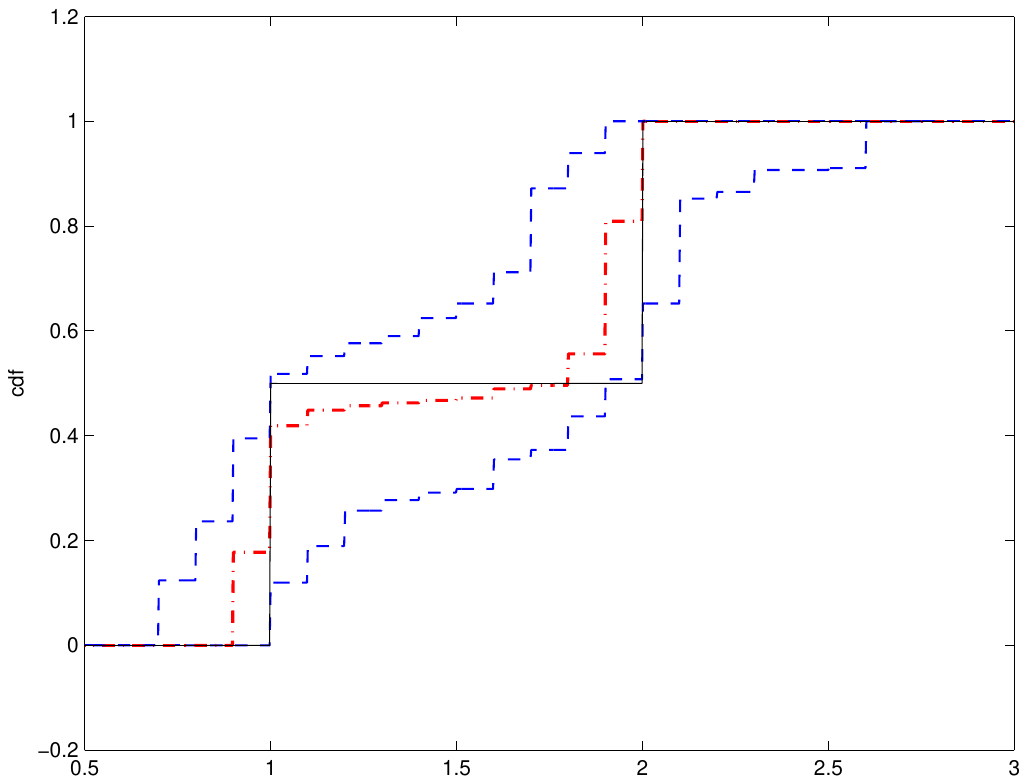}\vspace{-20pt}\\ \vspace{-20pt}
$(p,n)=(200,100)$ & \vspace{-10pt}\\ \vspace{-30pt}
\hspace{10pt}\includegraphics[width=2.8in, height=2.8in, bb = 50 50 650 650]{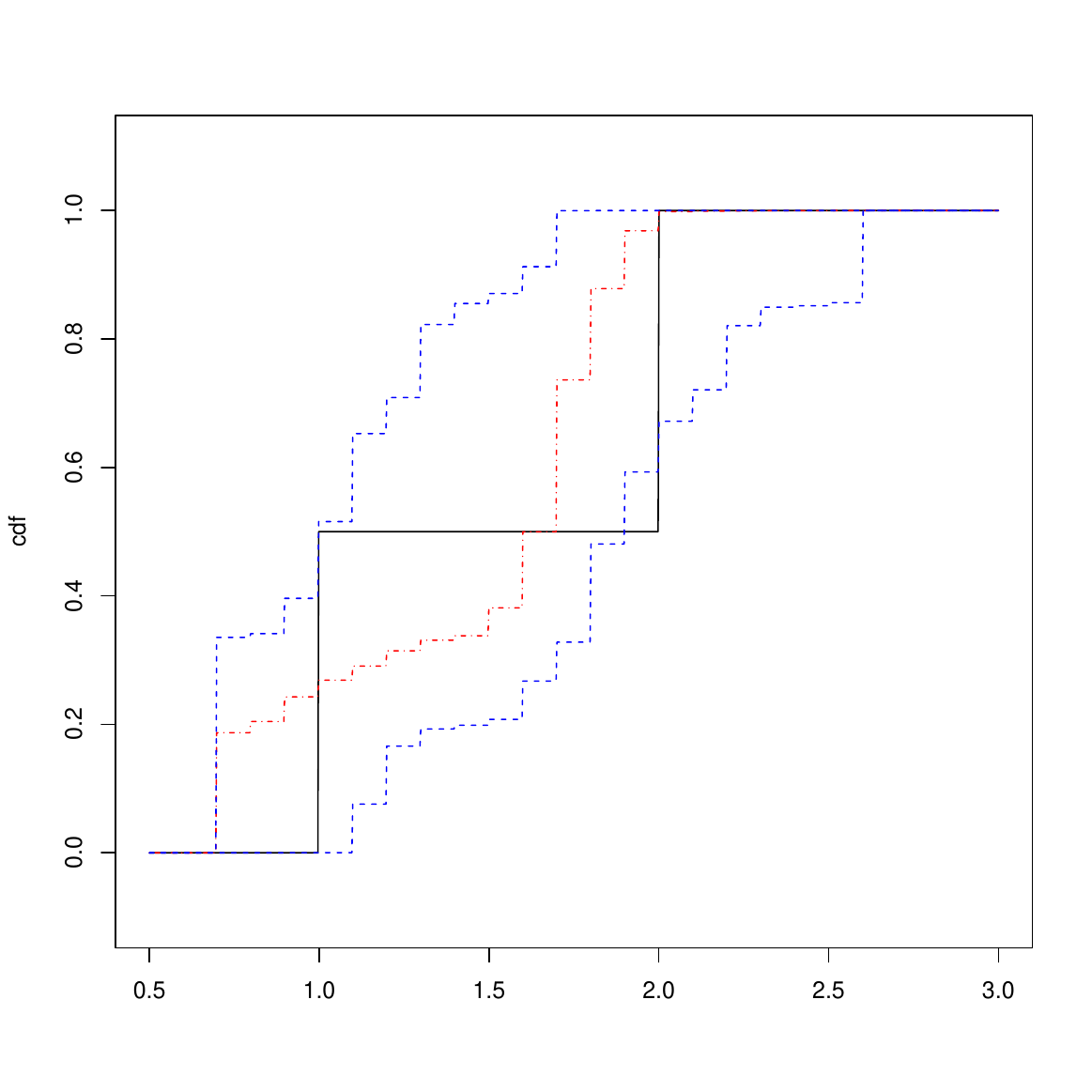}
& \\
\end{tabular}
}
\vspace{25pt}
         \caption{Plot of Median and $90\%$ confidence band for spectral cdf of $\boldsymbol{\Sigma}$ corresponding to the \textbf{case 2.1}. \textbf{Green} \textbf{Dash-Dot Red} curve: median, \textbf{Dashed Blue} curve: $90\%$ confidence band, \textbf{Black Solid} curve: true spectral cdf}
\end{center}
\end{figure}

\begin{figure}[H] \label{fig:AR2_A2}
\begin{center}
\resizebox{\textwidth}{!}{%
\begin{tabular}{cc}
$(p,n)=(400,1600)$ & $(p,n)=(400,800)$ \\ \vspace{-20pt}
\includegraphics[width=2.9in, height=3in, bb = 50 50 650 650]{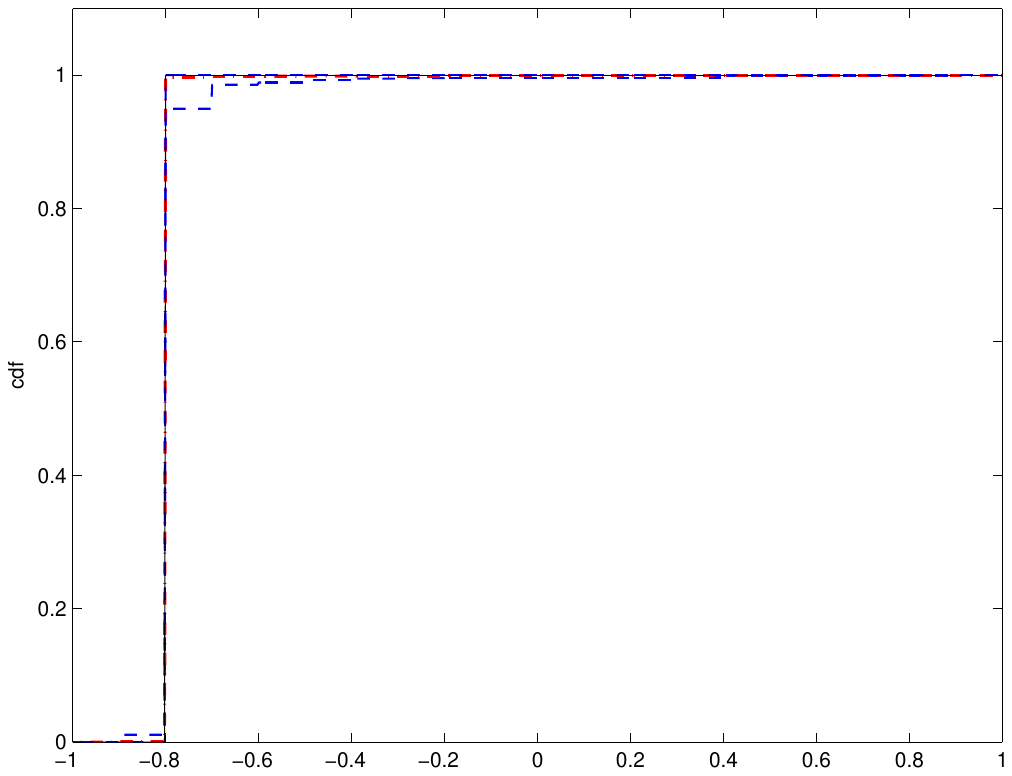}
& \includegraphics[width=2.9in, height=3in, bb = 50 50 650 650]{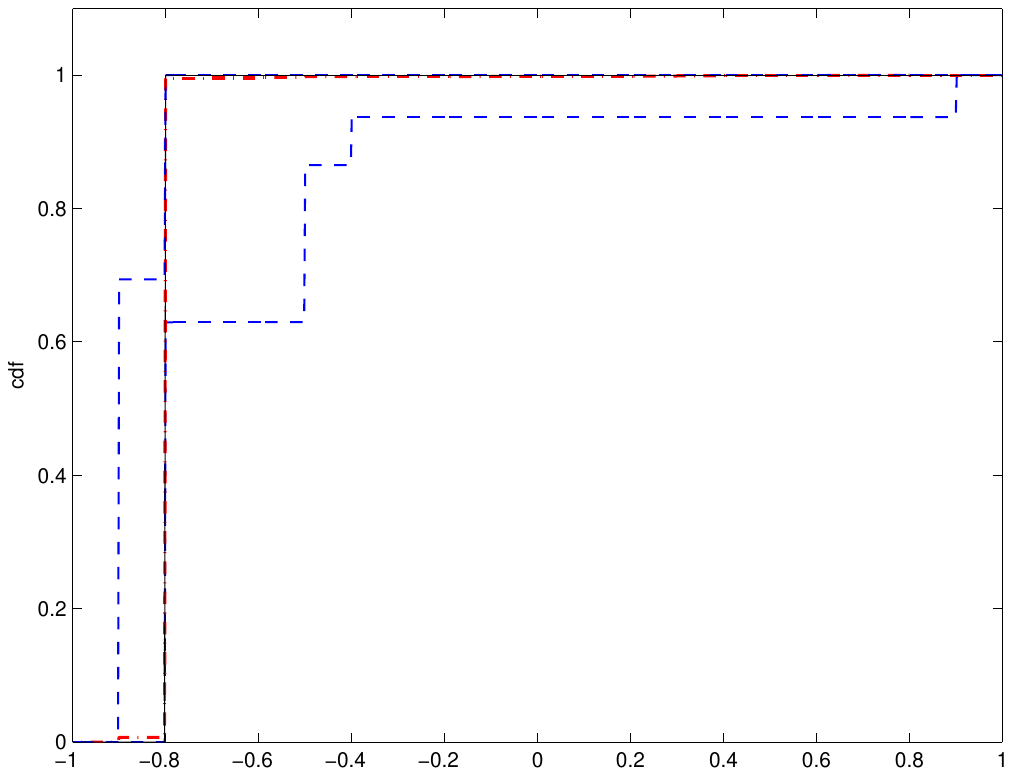}\vspace{-20pt}\\ \vspace{-10pt}
$(p,n)=(200,800)$ & $(p,n)=(200,400)$ \vspace{-5pt}\\ \vspace{-10pt}
\includegraphics[width=2.9in, height=3in, bb = 50 50 650 650]{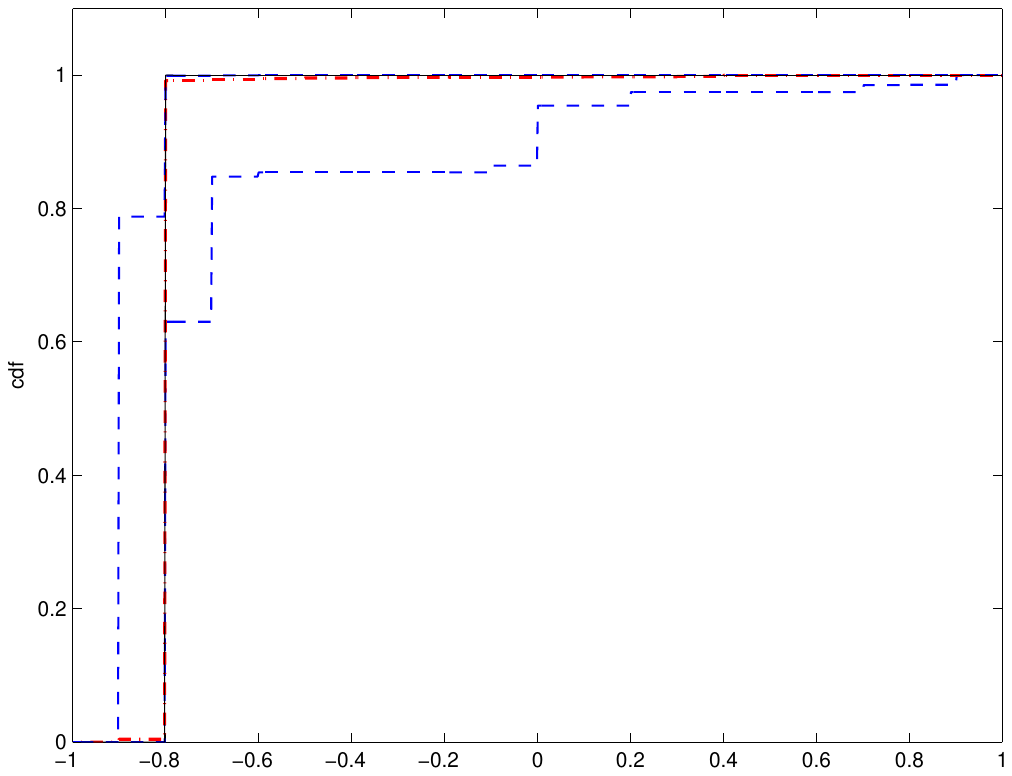} 
& \includegraphics[width=2.9in, height=3in, bb = 50 50 650 650]{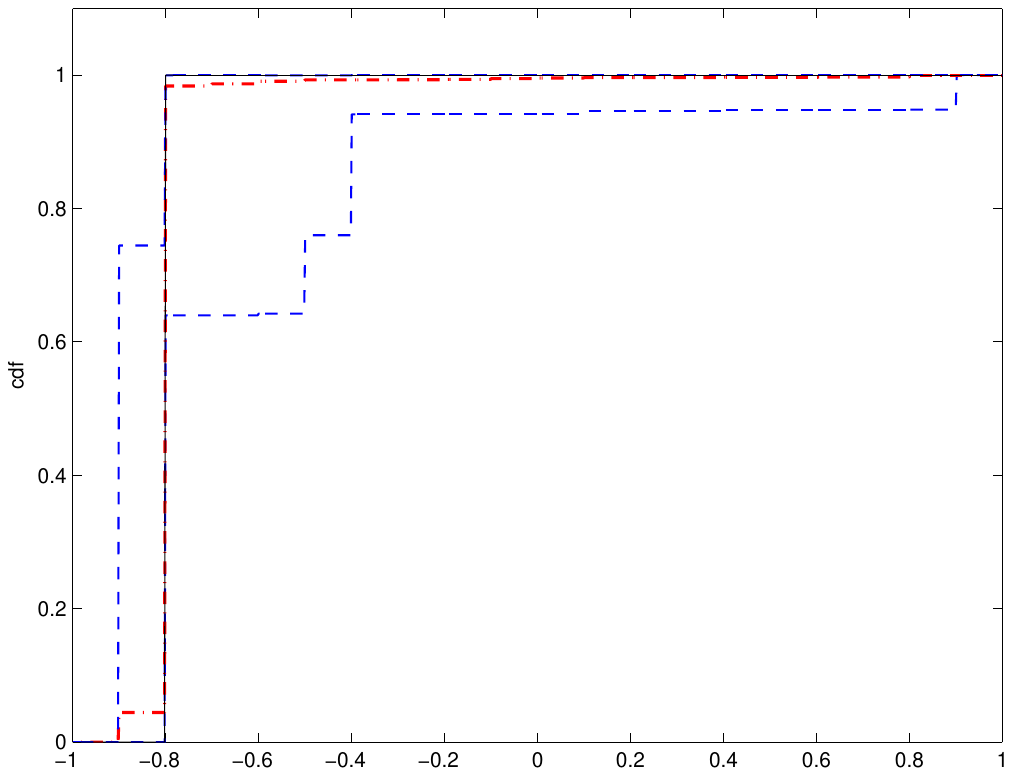} \vspace{-20pt}\\ \vspace{-20pt}
$(p,n)=(200,100)$ & \vspace{-10pt}\\ \vspace{-30pt}
\hspace{10pt}\includegraphics[width=2.8in, height=2.8in, bb = 50 50 650 650]{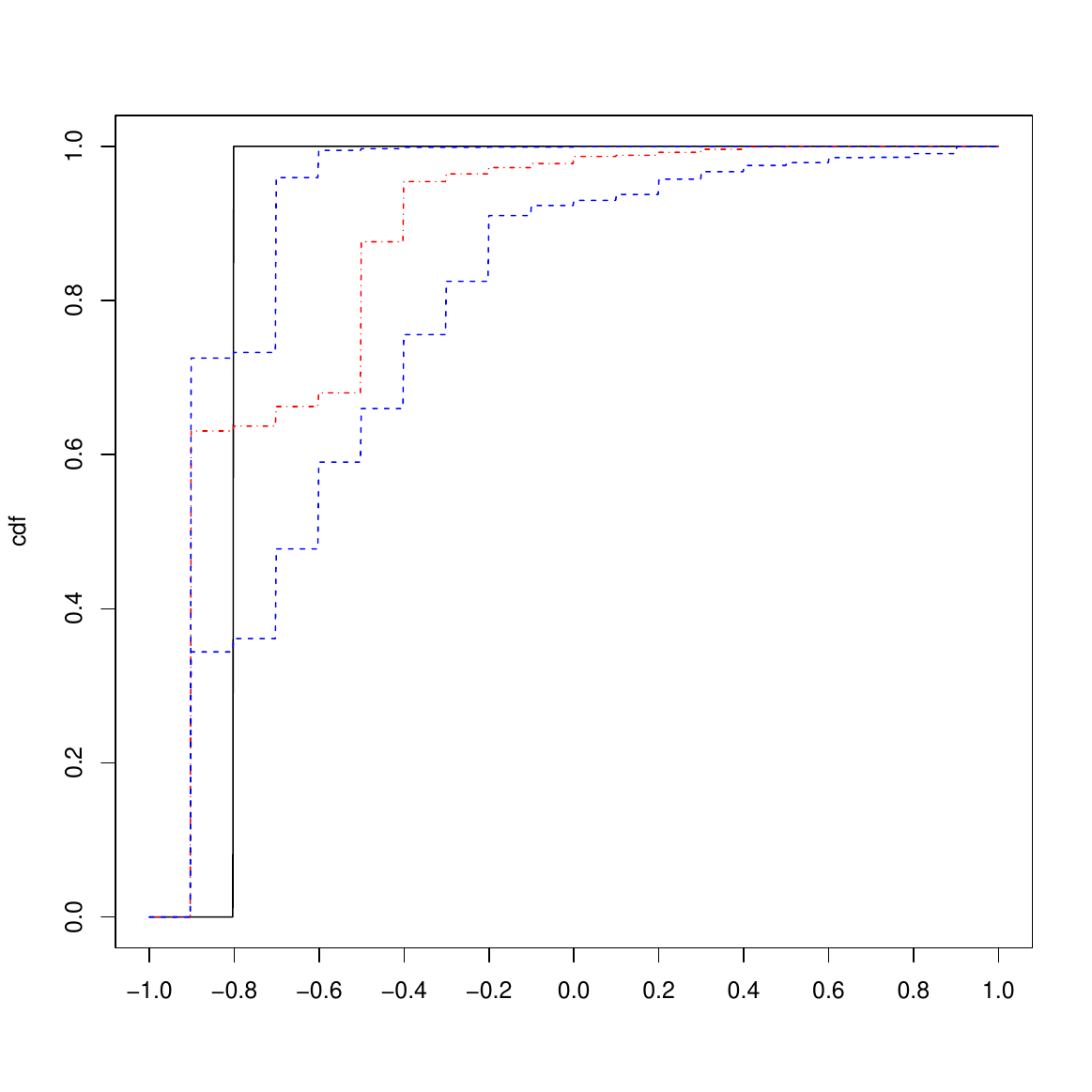}
& \\
\end{tabular}
}
\vspace{25pt}
         \caption{Plot of Median and $90\%$ confidence band for spectral cdf of $\boldsymbol{A_2}$ corresponding to the \textbf{case 2.2}. \textbf{Dash-Dot Red} curve: median, \textbf{Dashed Blue} curve: $90\%$ confidence band, \textbf{Black Solid} curve: true spectral cdf}
\end{center}
\end{figure}

\subsection{Data Analysis}
 Below, we present the plots of the proportion of variation explained (PVE) by leading factors of the log return series and the densities of the pairwise correlation between coordinates of the series $\mathbf{R}$ and $\mathbf{E}_R^{(1)}$.

The following table represents the percentage of times, in 100 times repetition of the simulation study, the model selection procedure selects the correct order of the underlying $AR$ process. $AR(2)\prec AR(1)$ means that the the $L_b$ measure is smaller when the $AR(2)$ model is considered in model fitting v.s. an $AR(1)$ model. 
\begin{center}
\begin{table}[htbp]
\centering
\begin{tabular}{|>{\centering\arraybackslash}m{3cm}|>{\centering\arraybackslash}m{1.2cm}|>{\centering\arraybackslash}m{1.2cm}|>{\centering\arraybackslash}m{1.2cm}|>{\centering\arraybackslash}m{1.2cm}|>{\centering\arraybackslash}m{1.2cm}|>{\centering\arraybackslash}m{1.2cm}|>{\centering\arraybackslash}m{1.2cm}|>{\centering\arraybackslash}m{1.2cm}|}
	\hline
	&\multicolumn{6}{|c|}{lag} \\
	\hline 
	&0&1&2&3&4&5\\
	\hline
	$AR(1)\prec AR(2)$& 0 & 0 & 0 & 0 & 0 & 0 \\
	\hline
	$AR(2)\prec AR(1)$& 100 & 100 & 100 & 100 & 100 & 100 \\
	\hline
\end{tabular}
\caption{Set 3; Percentages}
\end{table}
\end{center}


\begin{thebibliography}{4}

\bibitem[\protect\citeauthoryear{Bhattacharjee and Bose}{2016}]{bb16}
Bhattacharjee, M.\ and Bose, A.\ (2016).
Large sample behaviour of high dimensional autocovariance matrices.
{\it The Annals of Statistics}, {\bf 44}, 598--628.



\bibitem[\protect\citeauthoryear{Brockwell and Davis}{1991}]{bd92}
Brockwell, P.D.\ and Davis, R.A.\ (1991).
{\it Time Series Analysis: Theory and Methods}, 2nd ed.
Springer, New York.

\bibitem[\protect\citeauthoryear{Chatterjee}{2006}]{c06}
Chatterjee, S.\ (2006).
A generalization of the Lindeberg principle.
{\it The Annals of Probability}, {\bf 34}, 2061--2076.

\bibitem[\protect\citeauthoryear{El Karoui}{2008}]{e08}
El Karoui, N.\ (2008).
Spectrum estimation for large dimensional covariance matrices using
random matrix theory.
{\it The Annals of Statistics}, {\bf 36}, 2757--2790.

\bibitem[\protect\citeauthoryear{Flury}{1984}]{flury84}
Flury, B.N.\ (1984).
Common principal components in $k$ groups.
{\it Journal of the American Statistical Association}, {\bf 79},
892--898.


\bibitem[\protect\citeauthoryear{Jin, Wang, Bai, Nair, and Harding}{2014}]{jwbnh14}
Jin, B.,\ Wang, C.,\ Bai, Z.D.,\ Nair, K.K.,\ and Harding, M.\ (2014).
Limiting spectral distribution of a symmetrized auto-cross covariance
matrix.
{\it The Annals of Applied Probability}, {\bf 24}, 1199--1225.

\bibitem[\protect\citeauthoryear{Ledoit and Wolf}{2015}]{lw15}
Ledoit, O.\ and Wolf, M.\ (2015).
Spectrum estimation: A unified framework for covariance matrix
estimation and PCA in large dimensions.
{\it Journal of Multivariate Analysis}, {\bf 139}, 360--384.

\bibitem[\protect\citeauthoryear{Liu, Aue, and Paul}{2015}]{lap15}
Liu, H.,\ Aue, A.,\ and Paul, D.\ (2015).
On the Mar\v{c}enko--Pastur law for linear time series.
{\it The Annals of Statistics}, {\bf 43}, 675--712.


\bibitem[\protect\citeauthoryear{Miettinen et al.}{2016}]{minott16}
Miettinen, J.,\ Illner, K.,\ Nordhausen, K.,\ Oja, H.,\ Taskinen, S.,
and Theis, F.J.\ (2016).
Separation of uncorrelated stationary time series using autocovariance
matrices.
{\it Journal of Time Series Analysis}, {\bf 37}, 337--354.

\bibitem[\protect\citeauthoryear{Miettinen et al.}{2014}]{mnot14}
Miettinen, J.,\ Nordhausen, K.,\ Oja, H., and \ Taskinen, S.,\ (2014). Deflation-based separation of uncorrelated stationary time series. {\it Journal of Time Series Analysis}, {\bf 37}, 337--354.

\bibitem[\protect\citeauthoryear{Namdari}{2018}]{NamdariThesis2018}
Namdari, J.\ (2018).
{\it Estimation of Spectral Distributions of a Class of
High-dimensional Linear Processes}.
Ph.D.\ Thesis, University of California, Davis.

\bibitem[\protect\citeauthoryear{Nordhausen and Ruiz-Gazen}{2022}]{nr22}
Nordhausen, K.\ and Ruiz-Gazen, A.\ (2022).
On the usage of joint diagonalization in multivariate statistics.
{\it Journal of Multivariate Analysis}, {\bf 188}, 104844.

\bibitem[\protect\citeauthoryear{Paul and Aue}{2014}]{pa13}
Paul, D.\ and Aue, A.\ (2014).
Random matrix theory in statistics: A review.
{\it Journal of Statistical Planning and Inference}, {\bf 150}, 1--29.

\bibitem[\protect\citeauthoryear{Wang, Aue, and Paul}{2017}]{wap17}
Wang, L.,\ Aue, A.,\ and Paul, D.\ (2017).
Spectral analysis of sample autocovariance matrices of a class of
linear time series in moderately high dimensions.
{\it Bernoulli}, {\bf 23}, 2181--2209.


\end{thebibliography}
\end{document}